\documentclass[12pt,leqno]{article}
\usepackage{amssymb}
\usepackage[mathscr]{eucal}
\usepackage{amsmath,amssymb,latexsym,theorem,bbm}
\usepackage{amsmath,amssymb,latexsym,theorem}
\usepackage{color,url}
\usepackage{appendix}

\newcommand{\CC}{\mathbb{C}}

\newcommand{\NN}{\mathbb{N}}

\newcommand{\RR}{\mathbb{R}}

\newcommand{\ZZ}{\mathbb{Z}}

\newcommand{\ba}{{\boldsymbol{a}}}
\newcommand{\bA}{{\boldsymbol{A}}}

\newcommand{\tbA}{\widetilde{\bA}}

\newcommand{\bB}{{\boldsymbol{B}}}

\newcommand{\tb}{\widetilde{b}}

\newcommand{\tbB}{{\widetilde{\bB}}}

\newcommand{\bc}{{\boldsymbol{c}}}
\newcommand{\bC}{{\boldsymbol{C}}}

\newcommand{\tC}{\widetilde{C}}

\newcommand{\be}{{\boldsymbol{e}}}

\newcommand{\Bf}{{\boldsymbol{f}}}

\newcommand{\bI}{{\boldsymbol{I}}}
\newcommand{\tI}{\widetilde{I}}

\newcommand{\bM}{{\boldsymbol{M}}}

\newcommand{\bN}{{\boldsymbol{N}}}

\newcommand{\bP}{{\boldsymbol{P}}}

\newcommand{\bQ}{{\boldsymbol{Q}}}

\newcommand{\br}{{\boldsymbol{r}}}

\newcommand{\bu}{{\boldsymbol{u}}}
\newcommand{\tbu}{\widetilde{\bu}}
\newcommand{\bv}{{\boldsymbol{v}}}
\newcommand{\bV}{{\boldsymbol{V}}}

\newcommand{\bx}{{\boldsymbol{x}}}
\newcommand{\bX}{{\boldsymbol{X}}}
\newcommand{\by}{{\boldsymbol{y}}}
\newcommand{\bY}{{\boldsymbol{Y}}}

\newcommand{\tbY}{\widetilde{\bY}}

\newcommand{\bz}{{\boldsymbol{z}}}
\newcommand{\bZ}{{\boldsymbol{Z}}}

\newcommand{\Bbeta}{{\boldsymbol{\beta}}}

\newcommand{\tBbeta}{\widetilde{\Bbeta}}

\newcommand{\bfeta}{{\boldsymbol{\eta}}}

\newcommand{\blambda}{{\boldsymbol{\lambda}}}

\newcommand{\bmu}{{\boldsymbol{\mu}}}

\newcommand{\bbeta}{{\boldsymbol{\beta}}}
\newcommand{\bDelta}{{\boldsymbol{\Delta}}}

\newcommand{\bSigma}{{\boldsymbol{\Sigma}}}

\newcommand{\bzero}{{\boldsymbol{0}}}

\newcommand{\cD}{{\mathcal D}}

\newcommand{\tE}{\widetilde{E}}

\newcommand{\cF}{{\mathcal F}}

\newcommand{\cL}{{\mathcal L}}

\newcommand{\cN}{{\mathcal N}}

\newcommand{\cU}{{\mathcal U}}

\newcommand{\cc}{\mathrm{c}}
\newcommand{\dd}{\mathrm{d}}
\newcommand{\ee}{\mathrm{e}}
\newcommand{\ii}{\mathrm{i}}

\newcommand{\EE}{\operatorname{\mathbb{E}}}

\newcommand{\PP}{\operatorname{\mathbb{P}}}

\renewcommand{\Re}{{\operatorname{Re}}}
\renewcommand{\Im}{{\operatorname{Im}}}

\newcommand{\tM}{\widetilde{M}}
\newcommand{\tN}{\widetilde{N}}

\newcommand{\tY}{\widetilde{Y}}

\newcommand{\vare}{\varepsilon}

\renewcommand{\mid}{\,|\,}

\renewcommand{\leq}{\leqslant}
\renewcommand{\geq}{\geqslant}

\newcommand{\stoch}{\stackrel{\PP}{\longrightarrow}}
\newcommand{\distr}{\stackrel{\cD}{\longrightarrow}}
\newcommand{\distrA}{\stackrel{\cD_A}{\longrightarrow}}
\newcommand{\distrw}{\stackrel{\cD_{\{w_{\bu,\bX_0>0}\}}}{\longrightarrow}}
\newcommand{\distre}{\stackrel{\cD}{=}}
\newcommand{\mean}{\stackrel{L_1}{\longrightarrow}}
\newcommand{\qmean}{\stackrel{L_2}{\longrightarrow}}

\newcommand{\as}{\stackrel{{\mathrm{a.s.}}}{\longrightarrow}}
\newcommand{\ase}{\stackrel{{\mathrm{a.s.}}}{=}}

\newcommand{\bbone}{\mathbbm{1}}

\newcommand{\proofend}{\hfill\mbox{$\Box$}}

\numberwithin{equation}{section}

\theoremstyle{change} \theorembodyfont{\em}
\newtheorem{Lem}{Lemma.}[section]
\newtheorem{Thm}[Lem]{Theorem.}
\newtheorem{Pro}[Lem]{Proposition.}
\newtheorem{Cor}[Lem]{Corollary.}
\newtheorem{Def}[Lem]{Definition.}

\theorembodyfont{\rm}
\newtheorem{Rem}[Lem]{Remark.}

\begin{document}

\begin{center}
 {\bfseries\Large
   Asymptotic behavior of projections of supercritical\\[1mm]
    multi-type continuous state and continuous time\\[3mm]
    branching processes with immigration}

\vspace*{3mm}

 {\sc\large
  M\'aty\'as $\text{Barczy}^{*,\diamond}$,
  \ Sandra $\text{Palau}^{**}$,
  \ Gyula $\text{Pap}^{***}$}

\end{center}

\vskip0.2cm

\noindent
 * MTA-SZTE Analysis and Stochastics Research Group,
   Bolyai Institute, University of Szeged,
   Aradi v\'ertan\'uk tere 1, H--6720 Szeged, Hungary.

\noindent
 ** Department of Statistics and Probability, Instituto de
 Investigaciones en Matem\'aticas Aplicadas y en Sistemas,
 Universidad Nacional Aut\'onoma de M\'exico, M\'exico.

\noindent
 *** Bolyai Institute, University of Szeged,
    Aradi v\'ertan\'uk tere 1, H-6720 Szeged, Hungary.

\noindent e-mail: barczy@math.u-szeged.hu (M. Barczy),
                  sandra@sigma.iimas.unam.mx  (S. Palau).

\noindent $\diamond$ Corresponding author.



\renewcommand{\thefootnote}{}
\vspace*{0.2cm}
\footnote{M\'aty\'as Barczy is supported by the J\'anos Bolyai Research
 Scholarship of the Hungarian Academy of Sciences.
 Sandra Palau is supported by the Royal Society Newton International Fellowship
 and by the EU-funded Hungarian grant EFOP-3.6.1-16-2016-00008.}

\vspace*{-10mm}

\begin{abstract}
Under a fourth order moment condition on the branching and a second order
 moment condition on the immigration mechanisms, we show that an appropriately scaled projection of a supercritical and irreducible continuous state and continuous time branching process with immigration on certain left non-Perron eigenvectors of the branching mean matrix is asymptotically mixed normal.
With an appropriate random scaling, under some conditional probability measure, we prove asymptotic normality as well.
In case of a non-trivial process, under a first order moment condition on the immigration mechanism, we
 also prove the convergence of the relative frequencies of distinct types of individuals on a suitable event; for instance, if the immigration mechanism does not vanish, then this convergence holds almost surely.
\end{abstract}

\section{Introduction}
\label{section_intro}

The asymptotic behavior of multi-type supercritical branching processes
 without or with immigration has been studied for a long time.
Kesten and Stigum \cite[Theorems 2.1, 2.2, 2.3, 2.4]{KesSti1966b} investigated
 the limiting behaviors of the inner products \ $\langle\ba, \bX_n\rangle$ \ as
 \ $n \to \infty$, \ where \ $\bX_n$, \ $n \in \{1, 2, \ldots\}$, \ is a
 supercritical, irreducible and positively regular $d$-type Galton--Watson branching process without immigration and
 \ $\ba \in \RR^d \setminus \{\bzero\}$ \ is orthogonal to the left Perron eigenvector of the branching mean matrix
 \ $\bM := (\EE(\langle\be_j, \bX_1\rangle \mid \bX_0 = \be_i))_{i,j\in\{1,\ldots,d\}}$ \ of the process, where \ $\be_1, \ldots, \be_d$ \ denotes the natural basis in \ $\RR^d$.
\ Of course, this can arise only if \ $d \in \{2, 3, \ldots\}$.
\ It is enough to consider the case of \ $\|\ba\| = 1$, \ when \ $\langle\ba, \bX_n\rangle$ \ is the scalar projection of \ $\bX_n$ \ on \ $\ba$.
\ The appropriate scaling factor of \ $\langle\ba, \bX_n\rangle$, \ $n \in \{1, 2, \ldots\}$, \ depends not only on the Perron eigenvalue \ $r(\bM)$
 \ (which is the spectral radius of \ $\bM$) and on the left and right Perron eigenvectors of \ $\bM$, \ but also  on the full spectral representation of \ $\bM$.
\ Badalbaev and Mukhitdinov \cite[Theorems 1 and 2]{BadMuk} extended these results of Kesten and Stigum \cite{KesSti1966b}, namely, they described in a more explicit way the asymptotic behavior of
 \ $(\langle\ba^{(1)}, \bX_n\rangle, \ldots, \langle\ba^{(d-1)}, \bX_n\rangle)$
 \ as \ $n \to \infty$, \ where \ $\{\ba^{(1)}, \ldots, \ba^{(d-1)}\}$ \ is a
 basis of the hyperplane in \ $\RR^d$ \ orthogonal to the left Perron eigenvector of \ $\bM$.
\ They also pointed out the necessity of considering the functionals above originated in statistical investigations
 for \ $\bX_n$, \ $n\in\{1,2,\ldots\}$.

Athreya \cite{Ath0, Ath1} investigated the limiting behavior of \ $\bX_t$ \ and the inner products
 \ $\langle\bv, \bX_t\rangle$ \ as \ $t \to \infty$, \ where
 \ $(\bX_t)_{t\in[0,\infty)}$ \ is a supercritical, positively regular and non-singular $d$-type continuous time Galton--Watson branching process without immigration and
 \ $\bv \in \CC^d$ \ is a right eigenvector corresponding to an eigenvalue \ $\lambda \in \CC$ \ of the infinitesimal generator \ $\bA$ \ of the branching mean matrix semigroup
 \ $\bM(t) := (\EE( \langle \be_j, \bX_t\rangle \mid \bX_0 = \be_i))_{i,j\in\{1,\ldots,d\}} = \ee^{t\bA}$, \ $t \in [0, \infty)$, \ of the process.
Under a first order moment condition on the branching distributions, denoting by \ $s(\bA)$ \ the maximum of the real parts of the eigenvalues of \ $\bA$, \
 it was shown that there exists a non-negative random variable \ $w_{\bu,\bX_0}$ \ such that
 \ $\ee^{-s(\bA)t}\bX_t$ \ converges to \ $w_{\bu,\bX_0}\bu$ \ almost surely  as \ $t\to\infty$, \ where \ $\bu$ \ denotes
 the left Perron eigenvector of the branching mean matrix \ $\bM(1)$.
\ Under a second order moment condition on the branching distributions, it was shown that if
\ $\Re(\lambda) \in \big(\frac{1}{2} s(\bA), s(\bA)\big]$, \ then \ $\ee^{-\lambda t} \langle\bv, \bX_t\rangle$ \ converges almost surely and in $L_2$ \ to a (complex) random
 variable as \ $t \to \infty$, \ and if \ $\Re(\lambda) \in \big(-\infty, \frac{1}{2} s(\bA)\big]$ \ and \ $\PP(w_{\bu,\bX_0} >0) > 0$, \
 then, under the conditional probability measure \ $\PP(\cdot\mid w_{\bu,\bX_0} >0)$, \ the limit distribution of
 \ $t^{-\theta} \ee^{-s(\bA)t/2} \langle\bv, \bX_t\rangle$ \ as \ $t \to \infty$ \ is mixed normal, where
 \ $\theta = \frac{1}{2}$ \ if \ $\Re(\lambda) = \frac{1}{2} s(\bA)$ \ and
 \ $\theta = 0$ \ if \ $\Re(\lambda) \in \big(-\infty, \frac{1}{2} s(\bA)\big)$.
\ Further, in case of \ $\Re(\lambda) \in \big(-\infty, \frac{1}{2} s(\bA)\big]$,
 \ under the conditional probability measure \ $\PP(\cdot\mid w_{\bu,\bX_0} >0)$,
 \ with an appropriate random scaling, asymptotic normality has been derived as well
 with an advantage that the limit laws do not depend on the initial value \ $\bX_0$.
\ We also recall that Athreya \cite{Ath0} described the asymptotic behaviour of \ $\EE(\vert \langle\bv, \bX_t\rangle^2\vert)$ \
 as \ $t\to\infty$ \ under a second order moment condition on the branching distributions.
These results have been extended by Athreya \cite{Ath2} for the inner products
 \ $\langle\ba, \bX_t\rangle$, \ $t \in [0, \infty)$, \ with arbitrary \ $\ba \in \CC^d$. \
Janson \cite[Theorem 3.1]{Janson} gave a functional version of  Athreya's above mentioned
 results in \cite{Ath0, Ath1}.
Under some weaker conditions than Athreya \cite{Ath0, Ath1}, Janson \cite{Janson} described the asymptotic behaviour
 of \ $(\langle\bv, \bX_{t+s}\rangle)_{s \in [0, \infty)}$ \ as  \ $t \to \infty$ \ by giving more explicit
 formulas for the asymptotic variances and covariances as well.
For a more detailed comparison of Athreya's and Janson's results, see Janson \cite[Section 6]{Janson}.

Kyprianou et al.\ \cite{KypPalRen} described the limit behavior of the inner
 product \ $\langle\bu, \bX_t\rangle$ \ as \ $t \to \infty$ \ for supercritical and
 irreducible $d$-type continuous state and continuous time branching processes
 (without immigration),
 \ where \ $\bu$ \ denotes the left Perron vector
 	of the branching mean matrix of \ $(\bX_t)_{t\in[0,\infty)}$.
Barczy et al.\ \cite{BarPalPap} started to investigate the limiting behavior of
 the inner products \ $\langle\bv, \bX_t\rangle$ \ as \ $t \to \infty$, \ where
 \ $(\bX_t)_{t\in[0,\infty)}$ \ is a supercritical and irreducible $d$-type
 continuous state and continuous time branching process with immigration (CBI process)
 and \ $\bv \in \CC^d$ \ is a left eigenvector corresponding to an eigenvalue \ $\lambda \in \CC$ \
 of the infinitesimal generator \ $\tbB$ \ of the branching mean matrix semigroup \ $\ee^{t\tbB}$,
 \ $t \in [0, \infty)$, \ of the process.
 Note that for each \ $t \in [0, \infty)$ \ and \ $i, j \in \{1, \ldots, d\}$, \ we have
 \ $\langle \be_i,  \ee^{t\tbB} \be_j \rangle
    = \EE( \langle \be_i, \bY_t\rangle \mid \bY_0 = \be_j)$,
 \ where \ $(\bY_t)_{t\in[0,\infty)}$ \ is a multi-type  continuous state and
 continuous time branching process without immigration and with the same branching mechanism as \ $(\bX_t)_{t\in[0,\infty)}$, \ so \ $\tbB$ \ plays the role of \ $\bA^\top$ \ in Athreya \cite{Ath1}, hence in our results the right and left eigenvectors are interchanged compared to Athreya \cite{Ath1}.
Under first order moment conditions on the branching and immigration mechanisms,
 it was shown that there exists a non-negative random variable \ $w_{\bu,\bX_0}$ \ such that
 \ $\ee^{-s(\tbB)t}\bX_t$ \ converges to \ $w_{\bu,\bX_0}\tbu$ \ almost surely as \ $t\to\infty$,
 \ where \ $\tbu$ \ is the right Perron vector of \ $\ee^{\tbB}$, \ see Barczy et al.\ \cite[Theorem 3.3]{BarPalPap}.
\ If \ $\bv$ \ is a left non-Perron eigenvector of the branching mean matrix
 \ $\ee^{\tbB}$, \ then this result implies that
 \ $\ee^{-s(\tbB)t} \langle\bv, \bX_t\rangle \to w_{\bu,\bX_0} \langle\bv, \tbu\rangle = 0$ \
 almost surely as \ $t \to \infty$, \ since \ $\langle\bv, \tbu\rangle = 0$ \ due
 to the so-called principle of biorthogonality (see, e.g., Horn and Johnson
 \cite[Theorem 1.4.7(a)]{HorJoh}), consequently, the scaling factor
 \ $\ee^{-s(\tbB)t}$ \ is not appropriate for describing the asymptotic behavior of
 the projection \ $\langle\bv, \bX_t\rangle$ \ as \ $t \to \infty$.
\ Under suitable moment conditions on the branching and immigration mechanisms, it was shown that if
\ $\Re(\lambda) \in \big(\frac{1}{2} s(\tbB), s(\tbB)\big]$, \ then \ $\ee^{-\lambda t} \langle\bv, \bX_t\rangle$ \ converges almost surely
 and in $L_1$ \ (in \ $L_2$) \ to a (complex) random variable as \ $t \to \infty$, \
 see Barczy et al.\ \cite[Theorems 3.1 and 3.4]{BarPalPap}.

The aim of the present paper is to continue the investigations of Barczy et al.\ \cite{BarPalPap}.
We will prove that under a fourth order moment condition on the branching mechanism and a second order
 moment condition on the immigration mechanism, if
\ $\Re(\lambda) \in \big(-\infty, \frac{1}{2} s(\tbB)\big]$, \ then the limit distribution of
\ $t^{-\theta} \ee^{-s(\tbB)t/2} \langle\bv, \bX_t\rangle$ \ as \ $t \to \infty$ \ is mixed normal, where
\ $\theta = \frac{1}{2}$ \ if \ $\Re(\lambda) = \frac{1}{2} s(\tbB)$ \ and
\ $\theta = 0$ \ if \ $\Re(\lambda) \in \big(-\infty, \frac{1}{2} s(\tbB)\big)$,
 \ see  parts (ii) and (iii) of Theorem \ref{convCBIweak1}.
If \ $\Re(\lambda) \in \big(-\infty, \frac{1}{2} s(\tbB)\big]$ \  and
 \ $(\bX_t)_{t\in[0,\infty)}$ \ is non-trivial (equivalently, \ $\PP(w_{\bu,\bX_0} > 0) > 0$, \ see Lemma \ref{trivial}),
 then under the conditional probability measure \ $\PP(\cdot\mid w_{\bu,\bX_0} > 0)$,
 \ with an appropriate random scaling, we prove asymptotic normality as well
 with an advantage that the limit laws do not depend on the initial value \ $\bX_0$,
 \ see Theorem \ref{convCBIweakr}.
For the asymptotic variances, explicit formulas are presented.
 In case of a non-trivial process, under a first order moment condition on the immigration mechanism, we
 also prove the convergence of the relative frequencies of distinct types of individuals on the event \ $\{w_{\bu,\bX_0} > 0\}$ \ (see Proposition \ref{Cor_relative_frequency}); for instance, if the immigration mechanism does not vanish, then this convergence holds almost surely (see Theorem \ref{convCBIasL1+}).

Now, we summary the novelties of our paper. We point out that we investigate the asymptotic
 behavior of the projections of a multi-type CBI process on certain left non-Perron eigenvectors of its
 branching mean matrix.
Our approach is based on a decomposition of the process
 \ $(\ee^{-\lambda t} \langle\bv, \bX_t\rangle)_{t\in[0,\infty)}$ \ as the sum
 of a deterministic process and three square-integrable martingales, see the beginning of the proof of part (iii) of Theorem \ref{convCBIweak1}.
For proving asymptotic normality of the martingales in question, we use a result due to Crimaldi and Pratelli \cite[Theorem 2.2]{CriPra}
 (see also Theorem \ref{THM_Cri_Pra}), which provides a set of sufficient conditions for the
 asymptotic normality of multivariate martingales.
These sufficient conditions are about the quadratic variation process and the jumps of the
 multivariate martingale in question.
In the course of checking the conditions of Theorem \ref{THM_Cri_Pra}, we need to study the asymptotic behaviour of the expectation of the running supremum of the jumps of a
compensated Poisson integral process having time dependent integrand
over an interval \ $[0,t]$ \ as \ $t \to \infty$.
\ There is a new interest in this type of questions, see, e.g., the paper of He and Li \cite{HeLi}
on the distributions of jumps of a single-type CBI process.

Next, we compare our methodology with the discrete-valued settings.
Athreya \cite{Ath1} decomposed \ $\ee^{-\lambda t}\langle\bv, \bX_t\rangle$ \ in three terms, where
 \ $(\bX_t)_{t\in[0,\infty)}$ \ is a supercritical, positively regular and non-singular d-type
 continuous time Galton--Watson branching process without immigration, and
 he showed that two of them are small in probability and, using the central limit theorem, the third one converges to the desired normal distribution.
Janson's proof \cite[Theorem 3.1]{Janson} for a functional extension of Athreya's results is based on a martingale convergence theorem
 (see \cite[Proposition 9.1]{Janson}) that relies on the convergence of the quadratic variation of an  \ $L_2$-locally bounded
 (see \cite[condition (9.2)]{Janson}) martingale sequence.
Then, he needed to define a suitable martingale sequence, and estimate its quadratic variation.
Observe that he asked for a finite second moment for the branching mechanism in order to have an \ $L_2$-locally bounded martingale
 (see \cite[assumption (A.2)]{Janson}).
In our case, where \ $(\bX_t)_{t\in\RR_+}$ \ is a supercritical and irreducible d-type CBI process,
 the three martingales appearing in the previously mentioned decomposition of
 \ $(\ee^{-\lambda t}\langle\bv, \bX_t\rangle)_{t\in[0,\infty)}$ \ turn out to be square-integrable
  under our moment assumptions of the branching and immigration mechanisms.
One of the three martingales in question is an integral with respect to a standard Wiener process,
 and the other two are integrals with respect to compensated Poisson measures.
The decomposition in question was derived using an SDE representation of \ $(\bX_t)_{t\in[0,\infty)}$ \
 together with an application of the multidimensional It\^{o}'s formula, see Barczy et al. \cite[Lemma 4.1]{BarLiPap3}.
Concerning our moment assumptions, in order to be able to check the conditions of the previously mentioned
 Theorem 2.2 in Crimaldi and Pratelli \cite{CriPra} (see also Theorem \ref{THM_Cri_Pra}) we need a fourth order moment condition
 on the branching mechanism and a second order moment condition on the immigration mechanism.
So our proof technique can not be considered as an easy adaption of that of Athreya's \cite{Ath0, Ath1} or that of Janson \cite[Theorem 3.1]{Janson}.

The paper is structured as follows.
In Section \ref{section_CBI}, we recall the definition of multi-type CBI processes together with the notion of irreducibility,
 and we introduce a classification of multi-type CBI processes as well.
Sections \ref{section_results} and \ref{section_proofs} contain our results and their proofs, respectively.
We close the paper with five appendices.
In Appendix \ref{deco_CBI} we recall a decomposition of multi-type CBI processes,
 Appendix \ref{det_CBI} is devoted to a description of deterministic projections of multi-type CBI processes (i.e., projections that are deterministic).
In Appendix \ref{section_stoch_fixed_point_eq}, based on Buraczewski et al. \cite[Proposition 4.3.2]{BurDamMik},
 we recall some mild conditions under which the solution of a stochastic fixed point equation is atomless.
Appendix \ref{second_moment_CBI} is devoted to the description of the asymptotic behaviour of the second moment of projections
 of multi-type CBI processes.
In Appendix \ref{ltm} we recall a result on the asymptotic behavior of multivariate martingales
 due to Crimaldi and Pratelli \cite[Theorem 2.2]{CriPra}, which serves us as a key tool for proving our results,
 see Theorem \ref{THM_Cri_Pra}.

\section{Preliminaries}
\label{section_CBI}

Let \ $\ZZ_+$, \ $\NN$, \ $\RR$, \ $\RR_+$, \ $\RR_{++}$ \ and \ $\CC$ \ denote
 the set of non-negative integers, positive integers, real numbers,
 non-negative real numbers, positive real numbers and complex numbers, respectively.
For \ $x , y \in \RR$, \ we will use the notations
 \ $x \land y := \min \{x, y\}$, \ $x \lor y := \max \{x, y\}$ \ and
 \ $x^+ := \max \{0, x\}$.
\ By \ $\langle\bx, \by\rangle := \sum_{j=1}^d x_j \overline{y_j}$, \ we denote
 the Euclidean inner product of \ $\bx = (x_1, \ldots, x_d)^\top \in \CC^d$
 \ and \ $\by = (y_1, \ldots, y_d)^\top \in \CC^d$, \ and by \ $\|\bx\|$ \ and
 \ $\|\bA\|$, \ we denote the induced norm of \ $\bx \in \CC^d$ \ and
 \ $\bA \in \CC^{d\times d}$, \ respectively.
By \ $r(\bA)$, \ we denote the spectral radius of \ $\bA \in \CC^{d\times d}$.
\ The null vector and the null matrix will be denoted by \ $\bzero$.
\ Moreover, \ $\bI_d \in \RR^{d\times d}$ \ denotes the identity matrix.
If \ $\bA \in \RR^{d\times d}$ \ is positive semidefinite, then \ $\bA^{1/2}$
 \ denotes the unique positive semidefinite square root of \ $\bA$.
\ If \ $\bA \in \RR^{d\times d}$ \ is strictly positive definite, then
 \ $\bA^{1/2}$ \ is strictly positive definite and \ $\bA^{-1/2}$ \ denotes the
 inverse of \ $\bA^{1/2}$.
\ The set of \ $d \times d$ \ matrices with non-negative off-diagonal entries (also called essentially non-negative matrices)
 is denoted by \ $\RR^{d\times d}_{(+)}$.
\ By \ $C^2_\cc(\RR_+^d,\RR)$, \ we denote the set of twice continuously
 differentiable real-valued functions on \ $\RR_+^d$ \ with compact support. By \ $B(\RR_+^d,\RR)$, \ we denote the Banach space (endowed with the supremum norm) of real-valued
bounded Borel functions on \ $\RR_+^d$.
\  Convergence almost surely, in $L_1$, in $L_2$, in probability and in
 distribution will be denoted by \ $\as$, \ $\mean$, \ $\qmean$, \ $\stoch$ \ and
 \ $\distr$, \ respectively.
 For an event \ $A$ \ with \ $\PP(A) > 0$,
 \ let \ $\PP_A(\cdot) := \PP(\cdot\mid A) = \PP(\cdot \cap A) / \PP(A)$ \ denote
 the conditional probability measure given \ $A$, \ and let \ $\distrA$ \ denote
 convergence in distribution under the conditional probability measure \ $\PP_A$.
\ Almost sure equality and equality in distribution will be denoted by \ $\ase$
 \ and \ $\distre$, \ respectively.
If \ $\bV \in \RR^{d\times d}$ \ is symmetric and positive semidefinite, then
 \ $\cN_d(\bzero, \bV)$ \ denotes the $d$-dimensional normal distribution with
 zero mean and variance matrix \ $\bV$.
\ Throughout this paper, we make the conventions \ $\int_a^b := \int_{(a,b]}$
 \ and \ $\int_a^\infty := \int_{(a,\infty)}$ \ for any \ $a, b \in \RR$ \ with
 \ $a < b$.

\begin{Def}\label{Def_admissible}
A tuple \ $(d, \bc, \Bbeta, \bB, \nu, \bmu)$ \ is called a set of admissible
 parameters if
 \renewcommand{\labelenumi}{{\rm(\roman{enumi})}}
 \begin{enumerate}
  \item
   $d \in \NN$,
  \item
   $\bc = (c_i)_{i\in\{1,\ldots,d\}} \in \RR_+^d$,
  \item
   $\Bbeta = (\beta_i)_{i\in\{1,\ldots,d\}} \in \RR_+^d$,
  \item
   $\bB = (b_{i,j})_{i,j\in\{1,\ldots,d\}} \in \RR^{d \times d}_{(+)}$,
  \item
   $\nu$ \ is a Borel measure on \ $\cU_d := \RR_+^d \setminus \{\bzero\}$
    \ satisfying \ $\int_{\cU_d} (1 \land \|\br\|) \, \nu(\dd\br) < \infty$,
  \item
   $\bmu = (\mu_1, \ldots, \mu_d)$, \ where, for each
    \ $i \in \{1, \ldots, d\}$, \ $\mu_i$ \ is a Borel measure on
    \ $\cU_d$ \ satisfying
    \[
      \int_{\cU_d}
       \biggl[\|\bz\| \land \|\bz\|^2
              + \sum_{j \in \{1, \ldots, d\} \setminus \{i\}} (1 \land z_j)\biggr]
       \mu_i(\dd\bz)
      < \infty .
    \]
  \end{enumerate}
\end{Def}

\begin{Thm}\label{CBI_exists}
Let \ $(d, \bc, \Bbeta, \bB, \nu, \bmu)$ \ be a set of admissible parameters.
Then there exists a unique conservative transition semigroup \ $(P_t)_{t\in\RR_+}$
 \ acting on  \ $B(\RR_+^d,\RR)$ \ such that its
 Laplace transform has a representation
 \[
  \int_{\RR_+^d} \ee^{- \langle \blambda, \by \rangle} P_t(\bx, \dd \by)
  = \ee^{- \langle \bx, \bv(t, \blambda) \rangle
         - \int_0^t \psi(\bv(s, \blambda)) \, \dd s} , \qquad
  \bx \in \RR_+^d, \quad \blambda \in \RR_+^d , \quad t \in \RR_+ ,
 \]
 where, for any \ $\blambda \in \RR_+^d$, \ the continuously differentiable
 function
 \ $\RR_+ \ni t \mapsto \bv(t, \blambda)
    = (v_1(t, \blambda), \ldots, v_d(t, \blambda))^\top \in \RR_+^d$
 \ is the unique locally bounded solution to the system of differential equations
 \[
   \partial_t v_i(t, \blambda) = - \varphi_i(\bv(t, \blambda)) , \qquad
   v_i(0, \blambda) = \lambda_i , \qquad i \in \{1, \ldots, d\} ,
 \]
 with
 \[
   \varphi_i(\blambda)
   := c_i \lambda_i^2 -  \langle \bB \be_i, \blambda \rangle
      + \int_{\cU_d}
         \bigl( \ee^{- \langle \blambda, \bz \rangle} - 1
                + \lambda_i (1 \land z_i) \bigr)
         \, \mu_i(\dd \bz)
 \]
 for \ $\blambda \in \RR_+^d$, \ $i \in \{1, \ldots, d\}$, \ and
 \[
   \psi(\blambda)
   := \langle \bbeta, \blambda \rangle
      + \int_{\cU_d}
         \bigl( 1 - \ee^{- \langle\blambda, \br\rangle} \bigr)
         \, \nu(\dd\br) , \qquad
   \blambda \in \RR_+^d .
 \]
\end{Thm}

Theorem \ref{CBI_exists} is a special case of Theorem 2.7 of Duffie et al.\
\cite{DufFilSch} with \ $m = d$, \ $n = 0$ \ and zero killing rate.
For more details, see Remark 2.5 in Barczy et al.\ \cite{BarLiPap2}.
\begin{Def}\label{Def_CBI}
A conservative Markov process with state space \ $\RR_+^d$ \ and with transition
 semigroup \ $(P_t)_{t\in\RR_+}$ \ given in Theorem \ref{CBI_exists} is called a
 multi-type CBI process with parameters \ $(d, \bc, \Bbeta, \bB, \nu, \bmu)$.
\ The function
 \ $\RR_+^d \ni \blambda
    \mapsto (\varphi_1(\blambda), \ldots, \varphi_d(\blambda))^\top \in \RR^d$
 \ is called its branching mechanism, and the function
 \ $\RR_+^d \ni \blambda \mapsto \psi(\blambda) \in \RR_+$ \ is called its
 immigration mechanism.
A multi-type CBI process with parameters \ $(d, \bc, \Bbeta, \bB, \nu, \bmu)$ \ is
 called a CB process (a continuous state and continuous time branching process
 without immigration) if \ $\Bbeta = \bzero$ \ and \ $\nu = 0$ \  (equivalently, \ $\psi=0$).
\end{Def}

Let \ $(\bX_t)_{t\in\RR_+}$ \ be a multi-type CBI process with parameters
 \ $(d, \bc, \Bbeta, \bB, \nu, \bmu)$ \ such that \ $\EE(\|\bX_0\|) < \infty$ \ and
 the moment condition
 \begin{equation}\label{moment_condition_m_new}
  \int_{\cU_d} \|\br\| \bbone_{\{\|\br\|\geq1\}} \, \nu(\dd\br) < \infty
 \end{equation}
 holds.
Then, by formula (3.4) in Barczy et al. \cite{BarLiPap2},
 \begin{equation}\label{EXcond}
  \EE(\bX_t \mid \bX_0 = \bx)
  = \ee^{t\tbB} \bx + \int_0^t \ee^{u\tbB} \tBbeta \, \dd u ,
  \qquad \bx \in \RR_+^d , \quad t \in \RR_+ ,
 \end{equation}
 where
 \begin{gather*}
  \tbB := (\tb_{i,j})_{i,j\in\{1,\ldots,d\}} , \qquad
  \tb_{i,j}
  := b_{i,j} + \int_{\cU_d} (z_i - \delta_{i,j})^+ \, \mu_j(\dd\bz) , \qquad
  \tBbeta := \Bbeta + \int_{\cU_d} \br \, \nu(\dd\br) ,
 \end{gather*}
 with \ $\delta_{i,j}:=1$ \ if \ $i = j$, \ and \ $\delta_{i,j} := 0$ \ if
 \ $i \ne j$.
\ Note that, for each \ $\bx \in \RR^d_+$, \ the function
 \ $\RR_+ \ni t \mapsto \EE(\bX_t \mid \bX_0 = \bx)$ \ is continuous, and
 \ $\tbB \in \RR^{d \times d}_{(+)}$ \ and \ $\tBbeta \in \RR_+^d$, \ since
 \[
   \int_{\cU_d} \|\br\| \, \nu(\dd\br) < \infty , \qquad
   \int_{\cU_d} (z_i - \delta_{i,j})^+ \, \mu_j(\dd \bz) < \infty , \quad
   i, j \in \{1, \ldots, d\} ,
 \]
 see Barczy et al. \cite[Section 2]{BarLiPap2}.
Further, \ $\EE(\bX_t \mid \bX_0 = \bx)$, \ $\bx \in \RR_+^d$, \ does not depend on
 the parameter \ $\bc$.
\ One can give probabilistic interpretations of the modified parameters \ $\tbB$
 \ and \ $\tBbeta$, \ namely,  for each \ $t \in \RR_+$, \ we have \ $\ee^{t\tbB} \be_j = \EE(\bY_t \mid \bY_0 = \be_j)$,
 \ $j \in \{1, \ldots, d\}$, \ and \ $t \tBbeta = \EE(\bZ_t \mid \bZ_0 = \bzero)$,
 \ where \ $(\bY_t)_{t\in\RR_+}$ \ and \ $(\bZ_t)_{t\in\RR_+}$ \ are multi-type CBI
 processes with parameters \ $(d, \bc, \bzero, \bB, 0, \bmu)$ \ and
 \ $(d, \bzero, \Bbeta, \bzero, \nu, \bzero)$, \ respectively, see formula
 \eqref{EXcond}.
The processes \ $(\bY_t)_{t\in\RR_+}$ \ and \ $(\bZ_t)_{t\in\RR_+}$ \ can be
 considered as pure branching (without immigration) and pure immigration (without
 branching) processes, respectively.
Consequently, \ $\ee^\tbB$ \ and \ $\tBbeta$ \ may be called the branching mean
 matrix and the immigration mean vector, respectively.
Note that the branching mechanism depends only on the parameters \ $\bc$, \ $\bB$
 \ and \ $\bmu$, \ while the immigration mechanism depends only on the parameters
 \ $\Bbeta$ \ and \ $\nu$.

If \ $(d, \bc, \Bbeta, \bB, \nu, \bmu)$ \ is a set of admissible parameters,
 \ $\EE(\|\bX_0\|) < \infty$ \ and the moment condition
 \eqref{moment_condition_m_new} holds, then the multi-type CBI process with parameters
 \ $(d, \bc, \Bbeta, \bB, \nu, \bmu)$ \ can be represented as a pathwise unique strong solution of the stochastic differential equation (SDE)
 \begin{align}\label{SDE_atirasa_dimd}
  \begin{split}
   \bX_t
   &=\bX_0
     + \int_0^t (\Bbeta + \tbB \bX_u) \, \dd u
     + \sum_{\ell=1}^d
        \int_0^t \sqrt{2 c_\ell \max \{0, X_{u,\ell}\}} \, \dd W_{u,\ell}
        \, \be_\ell \\
   &\quad
      + \sum_{\ell=1}^d
         \int_0^t \int_{\cU_d} \int_{\cU_1}
          \bz \bbone_{\{w\leq X_{u-,\ell}\}} \, \tN_\ell(\dd u, \dd\bz, \dd w)
      + \int_0^t \int_{\cU_d} \br \, M(\dd u, \dd\br)
  \end{split}
 \end{align}
 for \ $t \in\RR_+$, \ see, Theorem 4.6 and Section 5 in
 Barczy et al.~\cite{BarLiPap2}, where \eqref{SDE_atirasa_dimd} was proved only for
 \ $d \in \{1, 2\}$, \ but their method clearly works for all
 \ $d \in \NN$.
\  Here \ $X_{t,\ell}$, \ $\ell\in\{1,\ldots,d\}$, \ denotes
 the \ $\ell^{\mathrm{th}}$ \ coordinate of \ $\bX_t$, \ $\PP(\bX_0\in\RR_+^d)=1$,
 \ $(W_{t,1})_{t\in\RR_+}$, \ \ldots, \ $(W_{t,d})_{t\in\RR_+}$ \ are standard
 Wiener processes, \ $N_\ell$, \ $\ell \in \{1, \ldots, d\}$, \ and \ $M$ \ are Poisson random measures on
 \ $\RR_{++} \times \cU_d \times \RR_{++}$
 \ and on \ $\RR_{++} \times \cU_d$ \ with intensity measures
 \ $\dd u \, \mu_\ell(\dd\bz) \, \dd w$, \ $\ell \in \{1, \ldots, d\}$, \ and
 \ $\dd u \, \nu(\dd\br)$, \ respectively, such that \ $\bX_0$,
 \ $(W_{t,1})_{t\in\RR_+}$, \ldots, $(W_{t,d})_{t\in\RR_+}$,
 \ $N_1,\dots,N_d$ \ and \ $M$ \ are independent, and
 \ $\tN_\ell(\dd u, \dd\bz, \dd w)
    := N_\ell(\dd u, \dd\bz, \dd w) - \dd u \, \mu_\ell(\dd\bz) \, \dd w$,
 \ $\ell \in \{1, \ldots, d\}$.

Next we recall a classification of multi-type CBI processes.
For a matrix \ $\bA \in \RR^{d \times d}$, \ $\sigma(\bA)$ \ will denote the
 spectrum of \ $\bA$, \ that is, the set of all \ $\lambda \in \CC$ \ that are
 eigenvalues of \ $\bA$.
\ Then \ $r(\bA) = \max_{\lambda \in \sigma(\bA)} |\lambda|$ \ is the spectral
 radius of \ $\bA$.
\ Moreover, we will use the notation
 \[
   s(\bA) := \max_{\lambda \in \sigma(\bA)} \Re(\lambda) .
 \]
A matrix \ $\bA \in \RR^{d\times d}$ \ is called reducible if there exist a
 permutation matrix \ $\bP \in \RR^{ d \times d}$ \ and an integer \ $r$ \ with
 \ $1 \leq r \leq d-1$ \ such that
 \[
  \bP^\top \bA \bP
   = \begin{pmatrix} \bA_1 & \bA_2 \\ \bzero & \bA_3 \end{pmatrix},
 \]
 where \ $\bA_1 \in \RR^{r\times r}$, \ $\bA_3 \in \RR^{ (d-r) \times (d-r) }$,
 \ $\bA_2 \in \RR^{r \times (d-r) }$, \ and \ $\bzero \in \RR^{(d-r)\times r}$ \ is
 a null matrix.
A matrix \ $\bA \in \RR^{d\times d}$ \ is called irreducible if it is not
 reducible, see, e.g., Horn and Johnson
 \cite[Definitions 6.2.21 and 6.2.22]{HorJoh}.
We do emphasize that no 1-by-1 matrix is reducible.

\begin{Def}\label{Def_irreducible}
Let \ $(\bX_t)_{t\in\RR_+}$ \ be a multi-type CBI process with parameters
 \ $(d, \bc, \Bbeta, \bB, \nu, \bmu)$ \ such that the moment condition
 \eqref{moment_condition_m_new} holds.
Then \ $(\bX_t)_{t\in\RR_+}$ \ is called irreducible if \ $\tbB$ \ is irreducible.
\end{Def}

Recall that if \ $\tbB \in \RR^{d\times d}_{(+)}$ \ is irreducible, then
 \ $\ee^{t\tbB} \in \RR^{d \times d}_{++}$ \ for all \ $t \in \RR_{++}$, \ and
 \ $s(\tbB)$ \ is a real eigenvalue of \ $\tbB$, \ the algebraic and geometric
 multiplicities of \ $s(\tbB)$ \ is 1, and the real parts of the other eigenvalues
 of \ $\tbB$ \ are less than \ $s(\tbB)$.
\ Moreover, corresponding to the eigenvalue \ $s(\tbB)$ \ there exists a unique
 (right) eigenvector \ $\tbu \in \RR^d_{++}$ \ of \ $\tbB$ \  such that the sum of
 its coordinates is 1 which is also the unique (right) eigenvector of \ $\ee^\tbB$,
 \ called the right Perron vector of \ $\ee^\tbB$, \ corresponding to the
 eigenvalue \ $r(\ee^\tbB) = \ee^{s(\tbB)}$ \ of \ $\ee^\tbB$ \ such that the sum
 of its coordinates is 1.
\ Further, there exists a unique left eigenvector \ $\bu \in \RR^d_{++}$ \ of
 \ $\tbB$ \ corresponding to the eigenvalue \ $s(\tbB)$ \ with
 \ $\tbu^\top \bu = 1$, \ which is also the unique (left) eigenvector of
 \ $\ee^\tbB$, \ called the left Perron vector of \ $\ee^\tbB$, \ corresponding to
 the eigenvalue \ $r(\ee^\tbB) = \ee^{s(\tbB)}$ \ of \ $\ee^\tbB$ \ such that
 \ $\tbu^\top \bu = 1$.
\ Moreover, there exist \ $C_1, C_2, C_3, C_4 \in \RR_{++}$ \ such that
 \begin{gather}\label{C}
  \|\ee^{-s(\tbB)t} \ee^{t\tbB} - \tbu \bu^\top\| \leq C_1 \ee^{-C_2 t} , \qquad
  \|\ee^{t\tbB}\| \leq C_3 \ee^{s(\tbB)t} , \qquad t \in \RR_+ , \\
  \label{EX}
  \EE(\|\bX_t\|) \leq C_4 \ee^{s(\tbB)t} , \qquad t \in \RR_+ .
 \end{gather}
These Frobenius and Perron type results can be found, e.g., in Barczy and Pap
 \cite[Appendix A]{BarPap} and Barczy et al.\ \cite[(3.8)]{BarPalPap}.

We will need the following dichotomy of the expectation of an irreducible
 multi-type CBI process.

\begin{Lem}\label{expectation_CBI}
Let \ $(\bX_t)_{t\in\RR_+}$ \ be an irreducible multi-type CBI process with
 parameters \ $(d, \bc, \Bbeta, \bB, \nu, \bmu)$ \ such that
 \ $\EE(\|\bX_0\|) < \infty$ \ and the moment condition
 \eqref{moment_condition_m_new} holds.
Then either \ $\EE(\bX_t) = \bzero$ \ for all \ $t \in \RR_+$, \ or
 \ $\EE(\bX_t) \in \RR_{++}^d$ \ for all \ $t \in \RR_{++}$.
\ Namely, if \ $\PP(\bX_0 = \bzero) = 1$, \ $\Bbeta = \bzero$ \ and \ $\nu = 0$,
 \ then \ $\EE(\bX_t) = \bzero$ \ for all \ $t \in \RR_+$, \ and hence
 \ $\PP(\bX_t = \bzero) = 1$ \ for all \ $t \in \RR_+$, \ otherwise
 \ $\EE(\bX_t) \in \RR_{++}^d$ \ for all \ $t \in \RR_{++}$.
\end{Lem}

\noindent
\textbf{Proof.}
For each \ $t \in \RR_+$, \ by \eqref{EXcond}, we have
 \[
   \EE(\bX_t) = \ee^{t\tbB} \EE(\bX_0) + \int_0^t \ee^{u\tbB} \tBbeta \, \dd u ,
   \qquad t \in \RR_+^d .
 \]
Since \ $\ee^{u\tbB} \in \RR_{++}^{d\times d}$ \ for all \ $u \in \RR_{++}$,
 \ $\EE(\bX_0) \in \RR_+^d$ \ and \ $\tBbeta \in \RR_+^d$, \ we obtain the
 assertions.
\proofend

\begin{Def}\label{Def_nontrivial}
Let \ $(\bX_t)_{t\in\RR_+}$ \ be an irreducible multi-type CBI process with parameters
 \ $(d, \bc, \Bbeta, \bB, \nu, \bmu)$.
\ Then \ $(\bX_t)_{t\in\RR_+}$ \ is called trivial if \ $\PP(\bX_0 = \bzero) = 1$,
 \ $\Bbeta = \bzero$ \ and \ $\nu = 0$, \ equivalently, if
 \ $\PP(\bX_t = \bzero) = 1$ \ for all \ $t \in \RR_+$.
\ Otherwise \ $(\bX_t)_{t\in\RR_+}$ \ is called non-trivial.
\end{Def}

We do recall the attention that if \ $(\bX_t^{(1)})_{t\in\RR_+}$ \ and
 \ $(\bX_t^{(2)})_{t\in\RR_+}$ \ are multi-type CBI processes with parameters
 \ $(d, \bc^{(1)}, \Bbeta, \bB^{(1)}, \nu, \bmu^{(1)})$ \ and
 \ $(d, \bc^{(2)}, \Bbeta, \bB^{(2)}, \nu, \bmu^{(2)})$, \ respectively, \ $\bX_0^{(1)} \ase \bX_0^{(2)}$ \ and
 \ $(\bX_t^{(1)})_{t\in\RR_+}$ \ is trivial, then
 \ $(\bX_t^{(2)})_{t\in\RR_+}$ \ is also trivial.

\begin{Def}\label{Def_indecomposable_crit}
Let \ $(\bX_t)_{t\in\RR_+}$ \ be an irreducible multi-type CBI process with
 parameters \ $(d, \bc, \Bbeta, \bB, \nu, \bmu)$ \ such that
 \ $\EE(\|\bX_0\|) < \infty$ \ and the moment condition
 \eqref{moment_condition_m_new} holds.
Then \ $(\bX_t)_{t\in\RR_+}$ \ is called
 \[
   \begin{cases}
    subcritical & \text{if \ $s(\tbB) < 0$,} \\
    critical & \text{if \ $s(\tbB) = 0$,} \\
    supercritical & \text{if \ $s(\tbB) > 0$.}
   \end{cases}
 \]
\end{Def}
For motivations of Definitions \ref{Def_irreducible} and
 \ref{Def_indecomposable_crit}, see Barczy and Pap \cite[Section 3]{BarPap}.

\section{Results}
\label{section_results}

Now we present the main result of this paper. Recall that \ $\bu \in \RR^d_{++}$ \ is the left Perron vector of \ $\ee^{\tbB}$ \ corresponding to the eigenvalue \ $\ee^{s(\tbB)}$.

\begin{Thm}\label{convCBIweak1}
Let \ $(\bX_t)_{t\in\RR_+}$ \ be a supercritical and irreducible multi-type CBI
 process with parameters \ $(d, \bc, \Bbeta, \bB, \nu, \bmu)$ \ such that
 \ $\EE(\|\bX_0\|) < \infty$ \ and the moment condition \eqref{moment_condition_m_new} holds.
Let \ $\lambda \in \sigma(\tbB)$ \ and let \ $\bv \in \CC^d$ \ be a left
 eigenvector \ of \ $\tbB$ \ corresponding to the eigenvalue \ $\lambda$.
\renewcommand{\labelenumi}{{\rm(\roman{enumi})}}
\begin{enumerate}
 \item
  If \ $\Re(\lambda) \in \bigl(\frac{1}{2} s(\tbB), s(\tbB)\bigr]$ \ and the moment
   condition
   \begin{equation}\label{moment_condition_xlogx}
    \sum_{\ell=1}^d
     \int_{\cU_d} g(\|\bz\|) \bbone_{\{\|\bz\|\geq1\}} \, \mu_\ell(\dd \bz)
    < \infty
   \end{equation}
   with
   \[
     g(x)
     := \begin{cases}
         x^{\frac{s(\tbB)}{\Re(\lambda)}}
          & \text{if \ $\Re(\lambda)
                        \in \bigl(\frac{1}{2} s(\tbB), s(\tbB)\bigr)$,} \\
         x \log(x)
          & \text{if \ $\Re(\lambda) = s(\tbB)$
                  \ ($\Longleftrightarrow$ $\lambda = s(\tbB)$),}
         \end{cases}
     \qquad x \in \RR_{++}
   \]
   holds, then there exists a complex random variable \ $w_{\bv,\bX_0}$ \ with
   \ $\EE(|w_{\bv,\bX_0}|) < \infty$ \ such that
   \begin{equation}\label{convwv}
    \ee^{-\lambda t} \langle \bv, \bX_t \rangle \to w_{\bv,\bX_0} \qquad
    \text{as \ $t \to \infty$ \ in \ $L_1$ \ and almost surely.}
   \end{equation}
 \item
  If \ $\Re(\lambda) = \frac{1}{2} s(\tbB)$ \ and the moment condition
   \begin{equation}\label{moment_4_2}
    \sum_{\ell=1}^d
     \int_{\cU_d} \|\bz\|^4
      \bbone_{\{\|\bz\|\geq1\}} \, \mu_\ell(\dd\bz) < \infty , \qquad
    \int_{\cU_d} \|\br\|^2 \bbone_{\{\|\br\|\geq1\}} \, \nu(\dd \bz) < \infty
   \end{equation}
   holds, then
   \begin{equation}\label{convvweak2}
    t^{-1/2} \ee^{-s(\tbB)t/2}
    \begin{pmatrix}
     \Re(\langle \bv, \bX_t \rangle) \\
     \Im(\langle \bv, \bX_t \rangle)
    \end{pmatrix}
    \distr \sqrt{w_{\bu,\bX_0}} \, \bZ_\bv \qquad \text{as \ $t \to \infty$,}
   \end{equation}
   where \ $\bZ_\bv$ \ is a $2$-dimensional random vector such that
   \ $\bZ_\bv \distre \cN_2(\bzero, \bSigma_\bv)$ \ independent of
   \ $w_{\bu,\bX_0}$, \ where
   \begin{align}\label{help15_Sigma_v2}
     \bSigma_\bv
     := \frac{1}{2}
         \sum\limits_{\ell=1}^d
          \langle\be_\ell, \tbu\rangle
          \left(C_{\bv,\ell} \bI_2
                + \begin{pmatrix}
                   \Re(\tC_{\bv,\ell}) & \Im(\tC_{\bv,\ell}) \\
                   \Im(\tC_{\bv,\ell}) & -\Re(\tC_{\bv,\ell})
                  \end{pmatrix}
                  \bbone_{\{\Im(\lambda)=0\}}\right)
   \end{align}
   with
   \begin{align*}
    C_{\bv,\ell}
    &:= 2 |\langle\bv, \be_\ell\rangle|^2 c_\ell
        + \int_{\cU_d} |\langle\bv, \bz\rangle|^2 \, \mu_\ell(\dd\bz) , \qquad
     \ell \in \{1, \ldots, d\} , \\
    \tC_{\bv,\ell}
    &:= 2 \langle\bv, \be_\ell\rangle^2 c_\ell
        + \int_{\cU_d} \langle\bv, \bz\rangle^2 \, \mu_\ell(\dd\bz) , \qquad
     \ell \in \{1, \ldots, d\} .
   \end{align*}
 \item
  If \ $\Re(\lambda) \in \bigl(-\infty, \frac{1}{2} s(\tbB)\bigr)$ \ and the moment
   condition \eqref{moment_4_2} holds, then
   \begin{equation}\label{convvweak1}
    \ee^{-s(\tbB)t/2}
    \begin{pmatrix}
     \Re(\langle \bv, \bX_t \rangle) \\
     \Im(\langle \bv, \bX_t \rangle)
    \end{pmatrix}
    \distr \sqrt{w_{\bu,\bX_0}} \, \bZ_\bv \qquad \text{as \ $t \to \infty$,}
   \end{equation}
   where \ $\bZ_\bv$ \ is a $2$-dimensional random vector such that
   \ $\bZ_\bv \distre \cN_2(\bzero, \bSigma_\bv)$ \ independent of
   \ $w_{\bu,\bX_0}$, \ where
   \begin{align}\label{help15_Sigma_v}
    \bSigma_\bv
    := \frac{1}{2}
       \sum\limits_{\ell=1}^d
       \langle\be_\ell, \tbu\rangle
        \left\{\frac{C_{\bv,\ell}}{s(\tbB)-2\Re(\lambda)} \bI_2
              + \begin{pmatrix}
                \Re\Bigl(\frac{\tC_{\bv,\ell}}{s(\tbB)-2\lambda}\Bigr)
                 & \Im\Bigl(\frac{\tC_{\bv,\ell}}{s(\tbB)-2\lambda}\Bigr) \\[1.5mm]
                \Im\Bigl(\frac{\tC_{\bv,\ell}}{s(\tbB)-2\lambda}\Bigr)
                 & -\Re\Bigl(\frac{\tC_{\bv,\ell}}{s(\tbB)-2\lambda}\Bigr)
               \end{pmatrix}\right\}
   \end{align}
   with \ $C_{\bv,\ell}$, \ $\ell \in \{1, \ldots, d\}$, \ and \ $\tC_{\bv,\ell}$, \ $\ell \in \{1, \ldots, d\}$, \ defined in part \upshape{(ii)}.
\end{enumerate}
\end{Thm}

 First we have some remarks concerning the limit distributions in parts (ii) and (iii) of Theorem \ref{convCBIweak1}.
Note that under the moment condition \eqref{moment_4_2}, the moment condition \eqref{moment_condition_xlogx}
 holds for \ $\lambda=s(\tbB)$ \ and hence there exists a non-negative random variable \ $w_{\bu,\bX_0}$ \
 with \ $\EE(w_{\bu,\bX_0})<\infty$ \ such that \ $\ee^{-s(\tbB)t}\langle\bu,\bX_t\rangle \to w_{\bu,\bX_0}$ \
 as \ $t\to\infty$ \ in \ $L_1$ \ and almost surely.
  Observe that if \ $(\bX_t)_{t\in\RR_+}$ \ is not a trivial process (see Definition \ref{Def_nontrivial}) and \ $\bSigma_\bv\ne \bzero$,
 \ then the scaling factors \ $t^{-1/2}\ee^{-s(\tbB)t/2}$ \ and \ $\ee^{-s(\tbB)t/2}$ \ in parts (ii) and (iii)
 of Theorem \ref{convCBIweak1} are correct in the sense that the corresponding limits are non-degenerate random variables,
 since \ $\PP(w_{\bu,\bX_0} = 0)<1$ \ due to Theorem 3.1 in Barczy et al.\ \cite{BarPalPap} or to Lemma \ref{trivial}.
The correctness of the scaling factor in part (i) of Theorem \ref{convCBIweak1} will be studied later on, this motivates the forthcoming Theorem \ref{convCBIasL1+}.
Note also that Theorem \ref{convCBIweak1} is valid even if \ $\bSigma_\bv$ \ is not invertible.
In Proposition \ref{Pro_Sigma_invertible}, necessary and sufficient conditions are
 given for the invertibility of \ $\bSigma_\bv$ \ provided that \ $\EE(\Vert\bX_0\Vert^2) < \infty$,
 \ $\Im(\lambda) \ne 0$, \ and the moment condition
 \begin{equation}\label{moment_condition_CBI2}
  \sum_{\ell=1}^d
   \int_{\cU_d} \|\bz\|^2 \bbone_{\{\|\bz\|\geq1\}} \, \mu_\ell(\dd\bz) < \infty ,
  \qquad
  \int_{\cU_d} \|\br\|^2 \bbone_{\{\|\br\|\geq1\}} \, \nu(\dd \br) < \infty
 \end{equation}
  holds.

Moreover, in Proposition \ref{variance_asymptotics_CBI} under the moment condition \eqref{moment_condition_CBI2} together with
 \ $\EE(\Vert \bX_0\Vert^2)<\infty$ \ we describe the asymptotic behavior of the variance matrix of the real and imaginary
 parts of \ $\langle \bv, \bX_t \rangle$ \ as \ $t \to \infty$, \ which explains the phase
 transition at \ $\Re(\lambda) = \frac{1}{2} s(\tbB)$ \ in Theorem \ref{convCBIweak1}.
This result can be considered as an extension of Proposition B.1 in Barczy et al.\ \cite{BarPalPap} (see also Proposition \ref{second_moment_asymptotics_CBI}),
 where the asymptotic behaviour of the second absolute moment
 \ $\EE(\vert \langle \bv, \bX_t \rangle\vert^2)$ \  of \ $\langle \bv, \bX_t \rangle$ \ has been described as \ $t\to\infty$.
\ The proof of Proposition \ref{variance_asymptotics_CBI} is based on the decomposition of \ $\ee^{-\lambda t}\langle \bv, \bX_t \rangle$ \
 mentioned in the Introduction (see the beginning of the proof of part (iii) of Theorem \ref{convCBIweak1})
 yielding an appropriate decomposition of \ $\EE(\langle \bv, \bX_t \rangle^2)$ \ containing \ $\EE(\langle \bv, \bX_0 \rangle^2)$, \ $\EE(\langle \bv, \bX_0 \rangle)$
 \ and \ $\EE(X_{u,\ell})$, \ $u\in[0,t]$, \ $\ell\in\{1,\ldots,d\}$.
So the proof of Proposition \ref{variance_asymptotics_CBI} can be finished by delicate estimations
 using the explicit form of \ $\EE(\bX_t \mid \bX_0=\bx)$, \ $\bx\in\RR_+^d$, \ $t\in\RR_+$, \ given in \eqref{EXcond}.

 In the next statement, sufficient conditions are derived for \ $\PP(w_{\bv,\bX_0} = 0) = 0$.
\ Note that in case of \ $\PP(w_{\bv,\bX_0} = 0) = 0$, \ the scaling factor
 \ $\ee^{-\lambda t}$ \ is correct in part (i) of Theorem \ref{convCBIweak1} in the sense that the limit is a non-degenerate random variable.

\begin{Thm}\label{convCBIasL1+}
Let \ $(\bX_t)_{t\in\RR_+}$ \ be a supercritical and irreducible multi-type CBI
 process with parameters \ $(d, \bc, \Bbeta, \bB, \nu, \bmu)$ \ such that
 \ $\EE(\|\bX_0\|) < \infty$ \ and the moment conditions
 \eqref{moment_condition_m_new} and \eqref{moment_condition_CBI2} hold.
Let \ $\lambda \in \sigma(\tbB)$ \ be such that
 \ $\Re(\lambda) \in \bigl(\frac{1}{2} s(\tbB), s(\tbB)\bigr]$, \ and let \ $\bv \in \CC^d$ \ be a
 left eigenvector of \ $\tbB$ \ corresponding to the eigenvalue \ $\lambda$.

If the conditions
 \renewcommand{\labelenumi}{{\rm(\roman{enumi})}}
 \begin{enumerate}
  \item
   $\tBbeta \ne \bzero$, \ i.e., \ $\Bbeta \ne \bzero$ \ or \ $\nu \ne 0$,
  \item
   $\nu(\{\br \in \cU_d : \langle\bv, \br\rangle \ne 0\}) > 0$, \ or \ there exists
    \ $\ell \in \{1, \ldots, d\}$ \ such that
    \ $\langle\bv, \be_\ell\rangle c_\ell \ne 0$ \ or
    \ $\mu_\ell(\{\bz \in \cU_d : \langle\bv, \bz\rangle \ne 0\}) > 0$
 \end{enumerate}
 hold, then the law of \ $w_{\bv,\bX_0}$ \ does not have atoms, where
 \ $w_{\bv,\bX_0}$ \ is given in part {\upshape (i)} of Theorem \ref{convCBIweak1}.
In particular, \ $\PP(w_{\bv,\bX_0} = 0) = 0$.

If the condition \textup{(ii)} does not hold, then
 \ $\PP(w_{\bv,\bX_0} = \langle\bv, \bX_0 + \lambda^{-1} \tBbeta\rangle ) = 1$,
 \ and in particular,
 \ $\PP(w_{\bv,\bX_0} = 0)
    = \PP(\langle\bv, \bX_0 + \lambda^{-1} \tBbeta\rangle = 0)$.

If \ $\lambda = s(\tbB)$, \ $\bv = \bu$ \ and the condition \textup{(i)} holds,
 then \ $\PP(w_{\bu,\bX_0} = 0) = 0$.

If \ $\lambda = s(\tbB)$, \ $\bv = \bu$, \ and the conditions \textup{(i)} and
 \textup{(ii)} do not hold, then \ $\PP(w_{\bu,\bX_0} = 0) = \PP(\bX_0 = \bzero)$.
\end{Thm}

Next, we show that with an appropriate random scaling in parts (ii) and (iii)
 in Theorem \ref{convCBIweak1}, \ $\langle \bv, \bX_t \rangle$ \ is asymptotically
 normal as \ $t \to \infty$ \ under the conditional probability measure
 \ $\PP(\cdot\mid w_{\bu,\bX_0} > 0)$, \  provided that
 \ $\PP(w_{\bu,\bX_0} > 0) > 0$.
\ Parts  (ii) and (iii) of the forthcoming Theorem \ref{convCBIweakr} are
 analogous to Theorems 1 and 2 and part 5 of Corollary 5 in Athreya
 \cite{Ath1}.
 First we give a necessary and sufficient condition for \ $w_{\bu,\bX_0} \ase 0$.

\begin{Lem}\label{trivial}
Suppose that \ $(\bX_t)_{t\in\RR_+}$ \ is a supercritical and irreducible multi-type
 CBI process with parameters \ $(d, \bc, \Bbeta, \bB, \nu, \bmu)$ \ such that
 \ $\EE(\|\bX_0\|) < \infty$, \ the moment condition \eqref{moment_condition_m_new}
 holds, and the moment condition \eqref{moment_condition_xlogx} holds for
 \ $\lambda = s(\tbB)$.
\ Then \ $w_{\bu,\bX_0} \ase 0$ \ if and only if \ $(\bX_t)_{t\in\RR_+}$ \ is a
 trivial process (equivalently, \  $\bX_0 \ase \bzero$ \ and \ $\tBbeta = \bzero$, \ see Lemma
 \ref{expectation_CBI} and Definition \ref{Def_nontrivial}).
\end{Lem}

\begin{Thm}\label{convCBIweakr}
Suppose that \ $(\bX_t)_{t\in\RR_+}$ \ is a supercritical, irreducible and
 non-trivial multi-type CBI
 process with parameters \ $(d, \bc, \Bbeta, \bB, \nu, \bmu)$ \ such that
 \ $\EE(\|\bX_0\|) < \infty$ \ and the moment condition \eqref{moment_condition_m_new} holds.
\renewcommand{\labelenumi}{{\rm(\roman{enumi})}}
\begin{enumerate}
 \item
  If \ $\Re(\lambda) \in \bigl(\frac{1}{2} s(\tbB), s(\tbB)\bigr]$ \ and the moment condition \eqref{moment_condition_xlogx} holds, then
   \begin{align*}
    &\bbone_{\{\bX_t\ne\bzero\}}
     \frac{1}{\langle \bu, \bX_t \rangle^{\Re(\lambda)/s(\tbB)}}
     \begin{pmatrix}
      \cos(\Im(\lambda)t) & \sin(\Im(\lambda)t) \\
      - \sin(\Im(\lambda)t) & \cos(\Im(\lambda)t)
     \end{pmatrix}
     \begin{pmatrix}
      \Re(\langle \bv, \bX_t \rangle) \\
      \Im(\langle \bv, \bX_t \rangle)
     \end{pmatrix} \\
    &\to \frac{1}{w_{\bu,\bX_0}^{\Re(\lambda)/s(\tbB)}}
         \begin{pmatrix}
          \Re(w_{\bv,\bX_0}) \\
          \Im(w_{\bv,\bX_0})
         \end{pmatrix}
     \qquad \text{as \ $t \to \infty$}
   \end{align*}
   on the event \ $\{w_{\bu,\bX_0} > 0\}$.
 \item
  If \ $\Re(\lambda) = \frac{1}{2} s(\tbB)$ \ and the moment condition \eqref{moment_4_2} holds,
  then, under the conditional probability measure \ $\PP(\cdot\mid w_{\bu,\bX_0} > 0)$, \ we have
 \[
  \bbone_{\{\langle\bu, \bX_t\rangle>1\}}
  \frac{1}{\sqrt{\langle\bu, \bX_t\rangle\log(\langle\bu, \bX_t\rangle)}}
  \begin{pmatrix}
   \Re(\langle\bv, \bX_t\rangle) \\
   \Im(\langle\bv, \bX_t\rangle)
  \end{pmatrix}
  \distrw \cN_2\biggl(\bzero, \frac{1}{s(\tbB)}\bSigma_\bv\biggr)
 \]
 as \ $t \to \infty$.
 \item
  If \ $\Re(\lambda) \in \bigl(-\infty, \frac{1}{2} s(\tbB)\bigr)$ \ and the moment condition \eqref{moment_4_2} holds, then,
  under the conditional probability measure \ $\PP(\cdot\mid w_{\bu,\bX_0} > 0)$, \ we have
 \[
  \bbone_{\{\bX_t\ne\bzero\}} \frac{1}{\sqrt{\langle \bu, \bX_t \rangle}}
  \begin{pmatrix}
   \Re(\langle \bv, \bX_t \rangle) \\
   \Im(\langle \bv, \bX_t \rangle)
  \end{pmatrix}
  \distrw \cN_2(\bzero, \bSigma_\bv) \qquad \text{as \ $t \to \infty$.}
 \]
\end{enumerate}
\end{Thm}

\begin{Rem}\label{Rem_relative_frequency_1}
The indicator function \ $\bbone_{\{\bX_t\ne\bzero\}}$ \ are needed in parts  (i) and
 (iii) of Theorem \ref{convCBIweakr}, and the indicator function \ $\bbone_{\{\langle\bu, \bX_t\rangle>1\}}$ \
  is needed in part  (ii) of Theorem \ref{convCBIweakr},
  since it can happen that \ $\PP(\bX_t = \bzero) > 0$, \ $t \in \RR_{++}$,
  \ even if \ $\tBbeta \ne \bzero$.
\ For example, if \ $(\bX_t)_{t\in\RR_+}$ \ is a multi-type CBI process with
 parameters \ $(d, \bc, \bzero, \bB, \nu, \bzero)$ \ such that
 \ $\bX_0 = \bzero$, \ $\bB$ \ is irreducible with \ $s(\bB) > 0$ \ and
 \ $\nu \ne 0$ \ with \ $\int_{U_d} (1 \lor \|\br\|) \, \nu(\dd \br) < \infty$,
 \ then \ $\tbB = \bB$, \ thus \ $(\bX_t)_{t\in\RR_+}$ \ is irreducible and
 supercritical.
One can choose, for instance, \ $d = 2$ \ and
 \[
   \bB = \begin{pmatrix} 1 & 1 \\ 1 & 1 \end{pmatrix} \in \RR_{(+)}^{2\times2} ,
 \]
 yielding that \ $\sigma(\bB)=\sigma(\tbB) = \{0,2\}$ \ and \ $s(\bB)=s(\tbB)=2$, \
  hence, by choosing \ $\lambda =0\in\sigma(\tbB)$, \ we have \ $\Re(\lambda) = 0\in(-\infty,1)=(-\infty,\frac{1}{2}s(\tbB))$,
  \ and, by choosing \ $\lambda = 2 \in \sigma(\tbB)$, \ we have \ $\Re(\lambda) = 2 \in (1, 2] = (\frac{1}{2}s(\tbB), s(\tbB)]$.
 If \ $d = 2$ \ and we choose
 \[
   \bB = \begin{pmatrix}
           3 & 1 \\
           1 & 3
         \end{pmatrix}
   \in \RR_{(+)}^{2\times2} ,
 \]
  then \ $\sigma(\bB) = \sigma(\tbB) = \{2,4\}$, \ $s(\bB) = s(\tbB)= 4$, \ and with \ $\lambda = 2$ \ we have \ $\Re(\lambda) = \frac{1}{2}s(\tbB)$.
\ Further, using \ $\tBbeta = \int_{\cU_d} \br \, \nu(\dd\br)$, \ by Lemma 4.1 in
 Barczy et al.\ \cite{BarLiPap3},
 \begin{align*}
  \bX_t &= \int_0^t \ee^{(t-u)\tbB} \tBbeta \, \dd u
           + \sum_{\ell=1}^d
              \int_0^t
               \ee^{(t-u)\tbB} \be_\ell \sqrt{2 c_\ell X_{u,\ell}}
               \, \dd W_{u,\ell}
           + \int_0^t \int_{\cU_d} \ee^{(t-u)\tbB} \br \, \tM(\dd u, \dd\br) \\
        &= \sum_{\ell=1}^d
            \int_0^t
             \ee^{(t-u)\tbB} \be_\ell \sqrt{2 c_\ell X_{u,\ell}} \, \dd W_{u,\ell}
           + \int_0^t \int_{\cU_d} \ee^{(t-u)\tbB} \br \, M(\dd u, \dd\br)
 \end{align*}
 for all \ $t \in \RR_+$, \ where \ $\tM(\dd u, \dd\br) := M(\dd u, \dd\br) - \dd u\, \nu(\dd \br)$.
\ Note that until the first jump of \ $M$ \ in \ $\RR_+ \times \cU_d$, \ the
 pathwise unique solution of this SDE is the identically zero process.
Hence, using that
 \ $\ee^{(t-u)\tbB} \in \RR_{++}^{d\times d}$ \ and \ $\ee^{(t-u)\tbB}$ \ is
 invertible \ for all \ $t \in \RR_{++}$ \ and \ $u \in  [0, t]$, \ we have
 \begin{align*}
  &\PP(\text{$\bX_s = \bzero$ \ for each \ $s \in [0, t]$})
    \geq  \PP(\text{$M$ \ has no point in
               \ $\{(u, \br) \in (0, t] \times \cU_d
                    : \ee^{(t-u)\tbB} \br \ne \bzero\}$}) \\
  &= \PP(\text{$M$ \ has no point in
               \ $\{(u, \br) \in (0, t] \times \cU_d : \br \ne \bzero\}$})
   = \ee^{-\int_0^t\int_{\cU_d}\bbone_{\{\br\ne\bzero\}}\,\dd u\,\nu(\dd\br)}
   = \ee^{-t\nu(\cU_d)}
 \end{align*}
 for all \ $t \in \RR_{++}$.
\ Consequently, since \ $\nu(\cU_d) < \infty$, \ we obtain
 \ $\PP(\bX_t = \bzero) > 0$, \ $t \in \RR_{++}$.
\proofend
\end{Rem}

Next we describe the asymptotic behavior of the relative frequencies of distinct types of individuals on the event \ $\{w_{\bu,\bX_0} > 0\}$.
\ For different models, one can find similar results in
 Jagers \cite[Corollary 1]{Jag} and Yakovlev and Yanev \cite[Theorem 2]{YakYan}.
For critical and irreducible multi-type CBI processes, see Barczy and Pap \cite[Corollary 4.1]{BarPap}.

\begin{Pro}\label{Cor_relative_frequency}
If \ $(\bX_t)_{t\in\RR_+}$ \ is a non-trivial, supercritical and irreducible multi-type CBI
 process with parameters \ $(d, \bc, \Bbeta, \bB, \nu, \bmu)$ \ such that
 \ $\EE(\|\bX_0\|) < \infty$ \ and the moment condition
 \eqref{moment_condition_m_new} holds, then for each \ $i, j \in \{1, \ldots, d\}$, \ we have
 \[
   \bbone_{\{\langle\be_j,\bX_t\rangle\ne0\}}
   \frac{\langle\be_i,\bX_t\rangle}{\langle\be_j,\bX_t\rangle}
   \to \frac{\langle\be_i,\tbu\rangle}{\langle\be_j,\tbu\rangle}
   \qquad \text{and} \qquad
   \bbone_{\{\bX_t\ne\bzero\}}
   \frac{\langle\be_i,\bX_t\rangle}{\sum_{k=1}^d\langle\be_k,\bX_t\rangle}
   \to {\langle\be_i, \tbu\rangle}
   \qquad \text{as \ $t \to \infty$}
 \]
 on the event \ $\{w_{\bu,\bX_0} > 0\}$.
\end{Pro}

\begin{Rem}
The indicator functions \ $\bbone_{\{\be_j^\top\bX_t\ne0\}}$ \ and
 \ $\bbone_{\{\bX_t\ne\bzero\}}$ \ are needed in Proposition
 \ref{Cor_relative_frequency}, since it can happen that \ $\PP(\bX_t = \bzero) > 0$, \ $t \in \RR_{++}$,
 \ see Remark \ref{Rem_relative_frequency_1}.
\proofend
\end{Rem}

\begin{Rem}
If \ $\PP(w_{\bu,\bX_0} = 0) = 0$, \ then the convergence in part (i) of  Theorem \ref{convCBIweakr} and in Proposition
 \ref{Cor_relative_frequency} holds almost surely, and the convergences in parts (ii) and (iii) hold
 under the unconditional probability measure \ $\PP$.
\proofend
\end{Rem}

\section{Proofs}
\label{section_proofs}

\noindent
\textbf{Proof of part (i) of Theorem \ref{convCBIweak1}.}
 This statement has been proved in Barczy et al.\ \cite[Theorem 3.1]{BarPalPap}.

\noindent
\textbf{Proof of part (iii) of Theorem \ref{convCBIweak1}.}
The proof is divided into three main steps.
First, we  decompose the process \ $(\ee^{-\lambda t} \langle\bv, \bX_t\rangle)_{t\in\RR_+}$ \ as the sum of
 a deterministic process and three square-integrable martingales.
We show that the deterministic process goes to zero as \ $t\rightarrow\infty$.
\ For proving asymptotic normality of the martingales in question, we use Theorem \ref{THM_Cri_Pra}
 due to Crimaldi and Pratelli \cite[Theorem 2.2]{CriPra} which provides a set of sufficient conditions for
 the asymptotic normality of multivariate martingales.
Then, the proof is complete as soon as we show that the conditions \eqref{ltm_cond1} and \eqref{ltm_cond2}
 of Theorem \ref{THM_Cri_Pra} are satisfied.
In the second and third steps, we prove that \eqref{ltm_cond1} and \eqref{ltm_cond2} are satisfied, respectively.

\textbf{Step 1.}
For each \ $t \in \RR_+$, \ we have the representation
 \ $\ee^{-\lambda t} \langle\bv, \bX_t\rangle
    = Z_t^{(0,1)} + Z_t^{(2)} + Z_t^{(3,4)} + Z_t^{(5)}$
 \ with
 \begin{align*}
  Z_t^{(0,1)}
  &:= \langle\bv, \bX_0\rangle
      + \langle\bv, \tBbeta\rangle \int_0^t \ee^{-\lambda u} \, \dd u , \\
  Z_t^{(2)}
  &:= \sum_{\ell=1}^d
       \langle\bv, \be_\ell\rangle
       \int_0^t \ee^{-\lambda u} \sqrt{2 c_\ell X_{u,\ell}} \, \dd W_{u,\ell} , \\
  Z_t^{(3,4)}
  &:= \sum_{\ell=1}^d
       \int_0^t \int_{\cU_d} \int_{\cU_1}
        \ee^{-\lambda u} \langle\bv, \bz\rangle
        \bbone_{\{w\leq X_{u-,\ell}\}}
        \, \tN_\ell(\dd u, \dd\bz, \dd w) , \\
  Z_t^{(5)}
  &:= \int_0^t \int_{\cU_d}
       \ee^{-\lambda u} \langle\bv, \br\rangle \, \tM(\dd u, \dd\br) ,
 \end{align*}
 see Barczy et al.\ \cite[Lemma 4.1]{BarLiPap3}.
Thus for each \ $t \in \RR_+$, \ we have
 \[
   \ee^{-s(\tbB)t/2} \langle \bv, \bX_t \rangle
   = \ee^{-(s(\tbB)-2\lambda)t/2}
     \bigl(Z_t^{(0,1)} + Z_t^{(2)} + Z_t^{(3,4)} + Z_t^{(5)}\bigr) .
 \]
First, we show
 \begin{equation}\label{Z01}
  \ee^{-(s(\tbB)-2\lambda)t/2} Z_t^{(0,1)} \as 0 \qquad \text{as \ $t \to \infty$.}
 \end{equation}
Indeed, if \ $\lambda = 0$, \ then
 \[
   \ee^{-(s(\tbB)-2\lambda)t/2} Z_t^{(0,1)}
   = \ee^{-s(\tbB)t/2}
     (\langle\bv, \bX_0\rangle + \langle\bv, \tBbeta\rangle t)
   \as 0 \qquad \text{as \ $t \to \infty$,}
 \]
 since \ $s(\tbB) \in \RR_{++}$.
\ Otherwise, if \ $\Re(\lambda) \in \bigl(-\infty, \frac{1}{2} s(\tbB)\bigr)$ \ and
 \ $\lambda \ne 0$, \ then
 \begin{align*}
  \ee^{-(s(\tbB)-2\lambda)t/2} Z_t^{(0,1)}
  &= \ee^{-(s(\tbB)-2\Re(\lambda))t/2+\ii\Im(\lambda)t}
     \biggl(\langle\bv, \bX_0\rangle
            - \frac{\langle\bv, \tBbeta\rangle}{\lambda}
              (\ee^{-\lambda t} - 1)\biggr) \\
  &= \biggl(\langle\bv, \bX_0\rangle
            + \frac{\langle\bv, \tBbeta\rangle}{\lambda}\biggr)
     \ee^{-(s(\tbB)-2\Re(\lambda))t/2+\ii\Im(\lambda)t}
     - \frac{\langle\bv, \tBbeta\rangle}{\lambda} \ee^{-s(\tbB)t/2}
   \as 0
 \end{align*}
 as \ $t \to \infty$.

For each \ $t \in \RR_+$, \ we have
 \[
   \begin{pmatrix}
    \Re\bigl(\ee^{-(s(\tbB)-2\lambda)t/2}
             \bigl(Z_t^{(2)} + Z_t^{(3,4)} + Z_t^{(5)}\bigr)\bigr) \\
    \Im\bigl(\ee^{-(s(\tbB)-2\lambda)t/2}
             \bigl(Z_t^{(2)} + Z_t^{(3,4)} + Z_t^{(5)}\bigr)\bigr)
   \end{pmatrix}
   = \bQ(t) \bM_t
 \]
 with
 \[
   \bQ(t)
   := \begin{pmatrix}
       \Re(\ee^{-(s(\tbB)-2\lambda)t/2}) & -\Im(\ee^{-(s(\tbB)-2\lambda)t/2}) \\
       \Im(\ee^{-(s(\tbB)-2\lambda)t/2}) & \Re(\ee^{-(s(\tbB)-2\lambda)t/2})
      \end{pmatrix} , \qquad t \in \RR_+ ,
 \]
 and
 \[
   \bM_t
   := \begin{pmatrix}
       \Re\bigl(Z_t^{(2)} + Z_t^{(3,4)} + Z_t^{(5)}\bigr) \\
       \Im\bigl(Z_t^{(2)} + Z_t^{(3,4)} + Z_t^{(5)}\bigr)
      \end{pmatrix} , \qquad t \in \RR_+ .
 \]
The assumption \ $\Re(\lambda) \in \bigl(-\infty, \frac{1}{2} s(\tbB)\bigr)$
 \ implies
 \[
   \bQ(t)
   = \ee^{-(s(\tbB)-2\Re(\lambda))t/2}
     \begin{pmatrix}
      \cos(\Im(\lambda)t) & -\sin(\Im(\lambda)t) \\
      \sin(\Im(\lambda)t) & \cos(\Im(\lambda)t)
     \end{pmatrix}
   \to \bzero \qquad \text{as \ $t \to \infty$.}
 \]
For each \ $t \in \RR_+$, \ we can write
 \ $\bM_t = \bM_t^{(2)} + \bM_t^{(3,4)} + \bM_t^{(5)}$ \ with
 \[
   \bM_t^{(2)}
   := \begin{pmatrix}
       \Re\bigl(Z_t^{(2)}\bigr) \\
       \Im\bigl(Z_t^{(2)}\bigr)
      \end{pmatrix} , \qquad
   \bM_t^{(3,4)}
   := \begin{pmatrix}
       \Re\bigl(Z_t^{(3,4)}\bigr) \\
       \Im\bigl(Z_t^{(3,4)}\bigr)
      \end{pmatrix} , \qquad
   \bM_t^{(5)}
   := \begin{pmatrix}
       \Re\bigl(Z_t^{(5)}\bigr) \\
       \Im\bigl(Z_t^{(5)}\bigr)
      \end{pmatrix} .
 \]
Note that under the moment condition \eqref{moment_4_2},
 \ $(\bM_t^{(2)})_{t\in\RR_+}$, \ $(\bM_t^{(3,4)})_{t\in\RR_+}$ \ and
 \ $(\bM_t^{(5)})_{t\in\RR_+}$ \ are square-integrable martingales (see, e.g.,
 Ikeda and Watanabe \cite[pages 55 and 63]{IkeWat}).
One can also observe that, by the decomposition of \ $(\ee^{-\lambda t} \langle\bv, \bX_t\rangle)_{t\in\RR_+}$ \ given at the beginning of this step,
  \ $(\ee^{-\lambda t} \langle\bv, \bX_t\rangle-\langle\bv, \tBbeta\rangle \int_0^t \ee^{-\lambda u} \, \dd u)_{t\in\RR_+}$ \
 is a martingale with respect to the filtration \ $\sigma(\bX_u : u \in[0,t])$, \ $t\in\RR_+$, \ which follows by Barczy et al.\ \cite[Lemma 2.6]{BarPalPap} as well.

The aim of the following discussion is to apply Theorem \ref{THM_Cri_Pra} for the
 2-dimensional martingale \ $(\bM_t)_{t\in\RR_+}$ \ with the scaling \ $\bQ(t)$,
 \ $t \in \RR_+$.

 \textbf{Step 2.} Now we prove that condition \eqref{ltm_cond1} of Theorem \ref{THM_Cri_Pra}
 holds for \ $(\bM_t)_{t\in\RR_+}$ \ with the scaling \ $\bQ(t)$, \ $t \in \RR_+$.
\ For each \ $t \in \RR_+$, \ by Theorem I.4.52 in Jacod and Shiryaev \cite{JSh}, we
 have
 \begin{align*}
  [\bM^{(2)}]_t
  &= \begin{pmatrix}
      [\Re\bigl(Z^{(2)}\bigr), \Re\bigl(Z^{(2)}\bigr)]_t
       & [\Re\bigl(Z^{(2)}\bigr), \Im\bigl(Z^{(2)}\bigr)]_t \\
      [\Im\bigl(Z^{(2)}\bigr), \Re\bigl(Z^{(2)}\bigr)]_t
       & [\Im\bigl(Z^{(2)}\bigr), \Im\bigl(Z^{(2)}\bigr)]_t
     \end{pmatrix} \\
  &= \begin{pmatrix}
      \langle\Re\bigl(Z^{(2)}\bigr), \Re\bigl(Z^{(2)}\bigr)\rangle_t
       & \langle\Re\bigl(Z^{(2)}\bigr), \Im\bigl(Z^{(2)}\bigr)\rangle_t \\
      \langle\Im\bigl(Z^{(2)}\bigr), \Re\bigl(Z^{(2)}\bigr)\rangle_t
       & \langle\Im\bigl(Z^{(2)}\bigr), \Im\bigl(Z^{(2)}\bigr)\rangle_t
     \end{pmatrix} \\
  &= 2 \sum_{\ell=1}^d
        c_\ell
        \int_0^t
         \begin{pmatrix}
          \Re(\ee^{-\lambda u} \langle\bv, \be_\ell\rangle) \\
          \Im(\ee^{-\lambda u} \langle\bv, \be_\ell\rangle)
         \end{pmatrix}
         \begin{pmatrix}
          \Re(\ee^{-\lambda u} \langle\bv, \be_\ell\rangle) \\
          \Im(\ee^{-\lambda u} \langle\bv, \be_\ell\rangle)
         \end{pmatrix}^\top
         X_{u,\ell}
         \, \dd u ,
 \end{align*}
 since \ $(\bM_t^{(2)})_{t\in\RR_+}$ \ is continuous, where
 \ $([\bM^{(2)}]_t)_{t\in\RR_+}$ \ and
 \ $(\langle\bM^{(2)}\rangle_t)_{t\in\RR_+}$ \ denotes the quadratic
 variation process and the predictable quadratic variation process of
 \ $(\bM^{(2)}_t)_{t\in\RR_+}$, \ respectively.
 Moreover, we have \ $\bM_t^{(3,4)} = \sum_{\ell=1}^d \tbY_t^{(\ell)}$ \ with
 \[
   \tbY_t^{(\ell)}
   := \begin{pmatrix}
       \Re\bigl(\tY_t^{(\ell)}\bigr) \\
       \Im\bigl(\tY_t^{(\ell)}\bigr)
      \end{pmatrix} , \qquad
   \tY_t^{(\ell)}
   := \int_0^t \int_{\cU_d} \int_{\cU_1}
       \ee^{-\lambda u} \langle\bv, \bz\rangle \bbone_{\{w\leq X_{u-,\ell}\}}
       \, \tN_\ell(\dd u, \dd\bz, \dd w)
 \]
 for \ $t \in \RR_+$ \ and \ $\ell \in \{1, \ldots, d\}$.
\ For each \ $t \in \RR_+$ \ and \ $\ell \in \{1, \ldots, d\}$,
 \ $(\tbY_t^{(\ell)})_{t\in\RR_+}$ \ is a square-integrable purely discontinuous
 martingale, see, e.g., Jacod and Shiryaev \cite[Definition II.1.27 and Theorem II.1.33]{JSh}).
Hence, for each \ $t \in \RR_+$ \ and \ $k, \ell \in \{1, \ldots, d\}$, \ by Lemma
 I.4.51 in Jacod and Shiryaev \cite{JSh}, we have
 \[
   [\tbY^{(k)}, \tbY^{(\ell)}]_t
   = \sum_{s\in[0,t]}
      (\tbY_s^{(k)} - \tbY_{s-}^{(k)})
      (\tbY_s^{(\ell)} - \tbY_{s-}^{(\ell)})^\top .
 \]
Further, by the proof of part (a) of Theorem II.1.33 in Jacod and Shiryaev \cite{JSh},
 for each \ $t \in \RR_+$ \ and \ $k \in \{1, \ldots, d\}$,
 \[
   [\tbY^{(k)}]_t
   = \int_0^t \int_{\cU_d} \int_{\cU_1}
      \begin{pmatrix}
        \Re(\ee^{-\lambda u} \langle\bv, \bz\rangle) \\
        \Im(\ee^{-\lambda u} \langle\bv, \bz\rangle)
       \end{pmatrix}
       \begin{pmatrix}
        \Re(\ee^{-\lambda u} \langle\bv, \bz\rangle) \\
        \Im(\ee^{-\lambda u} \langle\bv, \bz\rangle)
       \end{pmatrix}^\top
        \bbone_{\{w\leq X_{u-,k}\}}
        \, N_k(\dd u, \dd\bz, \dd w) .
 \]
The aim of the following discussion is to show that for each \ $t \in \RR_+$ \ and \ $k, \ell \in \{1, \ldots, d\}$
 \ with \ $k \ne \ell$, \ we have \ $[\tbY^{(k)}, \tbY^{(\ell)}]_t = \bzero$
 \ almost surely.
By the bilinearity of quadratic variation process, for all \ $\vare \in \RR_{++}$ \ and \ $t \in \RR_+$,
 \ we have
 \begin{align}\label{help18}
 \begin{split}
  [\tbY^{(k)}, \tbY^{(\ell)}]_t
    &= [\tbY^{(k,\vare)}, \tbY^{(\ell,\vare)}]_t
      + [\tbY^{(k)} - \tbY^{(k,\vare)} , \tbY^{(\ell)} - \tbY^{(\ell,\vare)} ]_t \\
    &\phantom{=\;}  + [\tbY^{(k,\vare)} , \tbY^{(\ell)} - \tbY^{(\ell,\vare)} ]_t
      + [\tbY^{(k)} - \tbY^{(k,\vare)} , \tbY^{(\ell,\vare)} ]_t ,
   \end{split}
 \end{align}
 where, for all \ $\vare \in \RR_{++}$, \ $k \in \{1, \ldots, d\}$ \ and \ $t \in \RR_+$,
 \[
   \tbY_t^{(k,\vare)}
   := \begin{pmatrix}
       \Re\bigl(\tY_t^{(k,\vare)}\bigr) \\
       \Im\bigl(\tY_t^{(k, \vare)}\bigr)
      \end{pmatrix} , \qquad
   \tY_t^{(k,\vare)}
   := \int_0^t \int_{\cU_d} \int_{\cU_1}
       \ee^{-\lambda u} \langle\bv, \bz\rangle
       \bbone_{\{\|\bz\|\geq\vare\}}
       \bbone_{\{w\leq X_{u-,k}\}}
       \, \tN_k(\dd u, \dd\bz, \dd w) ,
 \]
 which is well-defined and square-integrable, since, by \eqref{EX} and \eqref{moment_condition_CBI2},
 \begin{align*}
   &\int_0^t \int_{\cU_d}
     \ee^{-2\Re(\lambda)u} |\langle\bv, \bz\rangle|^2
     \bbone_{\{\|\bz\|\geq\vare\}} \EE(X_{u,k})
     \, \dd u \, \mu_k(\dd\bz) \\
   &\leq C_4 \|\bv\|^2
         \int_0^t
          \ee^{(s(\tbB)-2\Re(\lambda))u} \, \dd u
         \int_{\cU_d}
          \|\bz\|^2 \bbone_{\{\|\bz\|\geq\vare\}}
          \, \mu_k(\dd\bz)
    <\infty .
 \end{align*}
For each \ $\vare \in \RR_{++}$, \ $t \in \RR_+$ \ and \ $k, \ell \in \{1, \ldots, d\}$, \ we have
 \[
   [\tbY^{(k,\vare)}, \tbY^{(\ell,\vare)}]_t
   = \sum_{s\in[0,t]}
      (\tbY_s^{(k,\vare)} - \tbY_{s-}^{(k,\vare)})
      (\tbY_s^{(\ell,\vare)} - \tbY_{s-}^{(\ell,\vare)})^\top
   = \sum_{s\in[0,t]}
      (\bY_s^{(k,\vare)} - \bY_{s-}^{(k,\vare)})
      (\bY_s^{(\ell,\vare)} - \bY_{s-}^{(\ell,\vare)})^\top
 \]
 with
 \[
   \bY_t^{(k,\vare)}
   := \begin{pmatrix}
       \Re\bigl(Y_t^{(k,\vare)}\bigr) \\
       \Im\bigl(Y_t^{(k, \vare)}\bigr)
      \end{pmatrix} , \qquad
   Y_t^{(k,\vare)}
   := \int_0^t \int_{\cU_d} \int_{\cU_1}
       \ee^{-\lambda u} \langle\bv, \bz\rangle
       \bbone_{\{\|\bz\|\geq\vare\}}
       \bbone_{\{w\leq X_{u-,k}\}}
       \, N_k(\dd u, \dd\bz, \dd w) ,
 \]
 where the first equality follows by the proof of part (a) of Theorem II.1.33 in Jacod and Shiryaev \cite{JSh},
  and the second equality, by \eqref{EX}, part (vi) of Definition \ref{Def_admissible} and \eqref{moment_condition_CBI2}, since
 \begin{align*}
   &\int_0^t \int_{\cU_d}
     \ee^{-\Re(\lambda) u} |\langle\bv, \bz\rangle|
     \bbone_{\{\|\bz\|\geq\vare\}} \EE(X_{u,k})
     \, \dd u \,\mu_k(\dd\bz) \\
   &\leq C_4 \|\bv\|
         \int_0^t
          \ee^{(s(\tbB)-\Re(\lambda))u} \, \dd u
         \int_{\cU_d}
          \|\bz\| \bbone_{\{\|\bz\|\geq\vare\}}
          \, \mu_k(\dd\bz) \\
   &\leq \frac{C_4\|\bv\|}{\vare}
         \int_0^t
          \ee^{(s(\tbB)-\Re(\lambda))u} \, \dd u
         \int_{\cU_d} \|\bz\|^2 \,\mu_k(\dd\bz)
    < \infty ,
 \end{align*}
 and hence we have
 \begin{align*}
  \tY_t^{(k,\vare)}
  = Y_t^{(k,\vare)}
    - \int_0^t \int_{\cU_d} \int_{\cU_1}
       \ee^{-\lambda u} \langle\bv, \bz\rangle
       \bbone_{\{\|\bz\|\geq\vare\}}
       \bbone_{\{w\leq X_{u-,k}\}}
       \, \dd u \,\mu_k(\dd\bz) \, \dd w .
 \end{align*}
For each \ $\vare \in \RR_{++}$ \ and \ $k \in \{1, \ldots, d\}$, \ the jump times of
 \ $(\bY_t^{(k,\vare)})_{t\in\RR_+}$ \ is a subset of the jump times of the Poisson
 process \ $(N_k([0,t] \times \cU_d \times \cU_1))_{t\in\RR_+}$.
\ For each \ $k, \ell \in \{1, \ldots, d\}$ \ with \ $k \ne \ell$, \ the Poisson
 processes \ $(N_k([0,t] \times \cU_d \times \cU_1))_{t\in\RR_+}$ \ and
 \ $(N_\ell([0,t] \times \cU_d \times \cU_1))_{t\in\RR_+}$ \ are independent, hence
 they can jump simultaneously with probability zero, see, e.g., Revuz and Yor
 \cite[Chapter XII, Proposition 1.5]{RevYor}.
Consequently, for each \ $\vare \in \RR_{++}$, \ $t \in \RR_+$ \ and \ $k, \ell \in \{1, \ldots, d\}$
 \ with \ $k \ne \ell$, \ we have \ $[\tbY^{(k,\vare)}, \tbY^{(\ell,\vare)}]_t = \bzero$
 \ almost surely.

Moreover, for each \ $t \in \RR_+$, \ $\vare \in \RR_{++}$, \ $i, j \in \{1, 2\}$ \ and \ $k, \ell \in \{1, \ldots, d\}$ \ with \ $k \ne \ell$, \ by the Kunita--Watanabe inequality, we have
 \begin{align*}
  \bigl|\langle \be_i, [\tbY^{(k)} - \tbY^{(k,\vare)} , \tbY^{(\ell)} - \tbY^{(\ell,\vare)}]_t \be_j \rangle\bigr|
  & = \bigl|[\langle \be_i, \tbY^{(k)} - \tbY^{(k,\vare)} \rangle , \langle \be_j, \tbY^{(\ell)} - \tbY^{(\ell,\vare)}  \rangle ]_t \bigr|\\
  &\leq [\langle \be_i , \tbY^{(k)} - \tbY^{(k,\vare)} \rangle]_t^{1/2} \,
        [\langle \be_j ,  \tbY^{(k)} - \tbY^{(k,\vare)}\rangle ]_t^{1/2} , \\
  \bigl|\langle \be_i, [\tbY^{(k,\vare)} , \tbY^{(\ell)} - \tbY^{(\ell,\vare)}]_t \be_j \rangle \bigr|
  &\leq [\langle \be_i, \tbY^{(k,\vare)}\rangle ]_t^{1/2} \,
        [\langle \be_j, \tbY^{(\ell)} - \tbY^{(\ell,\vare)}  \rangle ]_t^{1/2} , \\
  \bigl|\langle \be_i, [\tbY^{(k)} - \tbY^{(k,\vare)} , \tbY^{(\ell,\vare)} ]_t \be_j \rangle\bigr|
  &\leq [\langle\be_i, \tbY^{(k)} - \tbY^{(k,\vare)} \rangle]_t^{1/2} \,
        [\langle\be_j, \tbY^{(\ell,\vare)} \rangle]_t^{1/2} .
\end{align*}
Hence it is enough to check that \ $[\langle \be_j, \tbY^{(\ell,\vare)} \rangle]_t$ \ is  stochastically bounded in \ $\vare \in \RR_{++}$ \ and
 \[
   [\langle\be_j, \tbY^{(\ell)} - \tbY^{(\ell,\vare)}\rangle]_t \mean 0 \qquad \text{as \ $\vare \downarrow 0$}
 \]
 for all \ $t \in \RR_+$, \ $j \in \{1, 2\}$ \ and \ $\ell \in \{1,\ldots, d\}$.
\  Indeed, in this case
 \begin{align*}
  &\bigl|\langle \be_i, [\tbY^{(k)} - \tbY^{(k,\vare)} , \tbY^{(\ell)} - \tbY^{(\ell,\vare)}]_t \be_j \rangle\bigr|
     \stoch 0 \qquad \text{as \ $\vare \downarrow 0$,}\\
  &\bigl|\langle \be_i, [\tbY^{(k,\vare)} , \tbY^{(\ell)} - \tbY^{(\ell,\vare)}]_t \be_j \rangle \bigr|
     \stoch 0 \qquad \text{as \ $\vare \downarrow 0$,}\\
  &\bigl|\langle \be_i, [\tbY^{(k)} - \tbY^{(k,\vare)} , \tbY^{(\ell,\vare)} ]_t \be_j \rangle\bigr|
        \stoch 0 \qquad \text{as \ $\vare \downarrow 0$,}
 \end{align*}
 and, by \eqref{help18}, for each \ $t \in \RR_+$ \ and \ $k, \ell \in \{1, \ldots, d\}$
 \ with \ $k \ne \ell$, \ we have \ $[\tbY^{(k)}, \tbY^{(\ell)}]_t = \bzero$
 \ almost surely.
By the proof of part (a) of Theorem II.1.33 in Jacod and Shiryaev \cite{JSh},
 \[
   [\tbY^{(\ell,\vare)}]_t
   = \int_0^t \int_{\cU_d} \int_{\cU_1}
      \begin{pmatrix}
        \Re(\ee^{-\lambda u} \langle\bv, \bz\rangle) \\
        \Im(\ee^{-\lambda u} \langle\bv, \bz\rangle)
       \end{pmatrix}
       \begin{pmatrix}
        \Re(\ee^{-\lambda u} \langle\bv, \bz\rangle) \\
        \Im(\ee^{-\lambda u} \langle\bv, \bz\rangle)
       \end{pmatrix}^\top
       \bbone_{\{\|\bz\|\geq\vare\}}
        \bbone_{\{w\leq X_{u-,\ell}\}}
        \, N_\ell(\dd u, \dd\bz, \dd w) ,
 \]
 and
 \begin{align*}
  &[\tbY^{(\ell)} - \tbY^{(\ell,\vare)}]_t \\
  &= \int_0^t \int_{\cU_d} \int_{\cU_1}
      \begin{pmatrix}
        \Re(\ee^{-\lambda u} \langle\bv, \bz\rangle) \\
        \Im(\ee^{-\lambda u} \langle\bv, \bz\rangle)
       \end{pmatrix}
       \begin{pmatrix}
        \Re(\ee^{-\lambda u} \langle\bv, \bz\rangle) \\
        \Im(\ee^{-\lambda u} \langle\bv, \bz\rangle)
       \end{pmatrix}^\top
       \bbone_{\{\|\bz\|<\vare\}}
        \bbone_{\{w\leq X_{u-,\ell}\}}
        \, N_\ell(\dd u, \dd\bz, \dd w) .
 \end{align*}
Consequently, using that \ $\|\bz \bz^\top\| \leq \|\bz\|^2$, \ $\bz \in \RR^2$, \ we have
 \[
   \bigl|[\langle\be_j, \tbY^{(\ell,\vare)}\rangle]_t\bigr|
   \leq \int_0^t \int_{\cU_d} \int_{\cU_1}
         |\ee^{-\lambda u} \langle \bv,\bz\rangle|^2
         \bbone_{\{w\leq X_{u-,\ell}\}}
         \, N_\ell(\dd u, \dd\bz, \dd w)
 \]
  for all \ $\vare \in \RR_{++}$ \ and \ $j\in\{1, 2\}$, \ where the right-hand side is finite almost surely, since
 \begin{align*}
   &\EE\left( \int_0^t \int_{\cU_d} \int_{\cU_1}
         |\ee^{-\lambda u} \langle \bv,\bz\rangle|^2
         \bbone_{\{w\leq X_{u-,\ell}\}}
         \, N_\ell(\dd u, \dd\bz, \dd w) \right)
     =  \int_0^t \int_{\cU_d}
        |\ee^{-\lambda u} \langle \bv,\bz\rangle|^2
        \EE(X_{u,\ell})\,\dd u\,\mu_\ell(\dd \bz) \\
   &\leq C_4 \|\bv\|^2 \int_0^t \ee^{(s(\tbB) -2\Re(\lambda))u}\,\dd u \int_{\cU_d} \|\bz\|^2 \,\mu_{\ell}(\dd\bz)
     <\infty.
 \end{align*}
Further,
 \begin{align*}
  \EE\bigl(\bigl|[\langle\be_j, \tbY^{(\ell)} - \tbY^{(\ell,\vare)}\rangle]_t\bigr|\bigr)
  &\leq \EE\left(\int_0^t \int_{\cU_d} \int_{\cU_1}
         |\ee^{-\lambda u} \langle \bv,\bz\rangle|^2
         \bbone_{\{\|\bz\|<\vare\}}
         \bbone_{\{w\leq X_{u-,\ell}\}}
         \, N_\ell(\dd u, \dd\bz, \dd w)\right) \\
   &\leq \int_0^t \int_{\cU_d} |\ee^{-\lambda u} \langle \bv,\bz\rangle|^2 \bbone_{\{\|\bz\|<\vare\}}
        \EE(X_{u,\ell})\,\dd u\,\mu_\ell(\dd \bz)\\
   & \leq C_4 \|\bv\|^2 \int_0^t \ee^{(s(\tbB) -2\Re(\lambda))u}\,\dd u \int_{\cU_d} \|\bz\|^2 \bbone_{\{\|\bz\|<\vare\}}\,\mu_{\ell}(\dd\bz)
   \to 0
 \end{align*}
 as \ $\vare \downarrow 0$.
\ Consequently, for each \ $t \in \RR_+$ \ and \ $k, \ell \in \{1, \ldots, d\}$
 \ with \ $k \ne \ell$, \ we have \ $[\tbY^{(k)}, \tbY^{(\ell)}]_t = \bzero$
 \ almost surely.

In a similar way,
 \[
   [\bM^{(5)}]_t
   = \int_0^t \int_{\cU_d}
      \begin{pmatrix}
       \Re(\ee^{-\lambda u} \langle\bv, \br\rangle) \\
       \Im(\ee^{-\lambda u} \langle\bv, \br\rangle)
      \end{pmatrix}
      \begin{pmatrix}
       \Re(\ee^{-\lambda u} \langle\bv, \br\rangle) \\
       \Im(\ee^{-\lambda u} \langle\bv, \br\rangle)
      \end{pmatrix}^\top
      M(\dd u, \dd\br) , \qquad t \in \RR_+ ,
 \]
 and \ $[\tbY^{(\ell)}, \bM^{(5)}]_t = \bzero$, \ $\ell \in \{1, \ldots, d\}$ \ almost surely.
Consequently, for each \ $t \in \RR_+$, \ we have
 \ $[\bM^{(3,4)} + \bM^{(5)}]_t = [\bM^{(3,4)}]_t + [\bM^{(5)}]_t$ \ with
 \ $[\bM^{(3,4)}]_t = \sum_{\ell=1}^d [\tbY^{(\ell)}]_t$.
\ Since \ $(\bM^{(2)}_t)_{t\in\RR_+}$ \ is a continuous martingale and
 \ $(\bM^{(3,4)}_t + \bM^{(5)}_t)_{t\in\RR_+}$ \ is a purely discontinuous
 martingale, by Corollary I.4.55 in Jacod and Shiryaev \cite{JSh}, we have
 \ $[\bM^{(2)}, \bM^{(3,4)} + \bM^{(5)}]_t = \bzero$, \ $t \in \RR_+$.
\ Consequently,
 \[
   [\bM]_t = [\bM^{(2)}]_t + [\bM^{(3,4)}]_t + [\bM^{(5)}]_t, \qquad t \in \RR_+ .
 \]
For each \ $t \in \RR_+$, \ we have
 \begin{align*}
   \bQ(t) [\bM^{(2)}]_t \bQ(t)^\top
    = 2 \sum_{\ell=1}^d
        c_\ell \int_0^t \Bf(t - \tau, \be_\ell) \ee^{-s(\tbB)\tau} X_{\tau,\ell}
        \, \dd \tau
 \end{align*}
 with
 \[
   \Bf(w, \bz)
   := \begin{pmatrix}
       \Re(\ee^{-(s(\tbB)-2\lambda)w/2} \langle\bv, \bz\rangle) \\
       \Im(\ee^{-(s(\tbB)-2\lambda)w/2} \langle\bv, \bz\rangle)
      \end{pmatrix}
      \begin{pmatrix}
       \Re(\ee^{-(s(\tbB)-2\lambda)w/2} \langle\bv, \bz\rangle) \\
       \Im(\ee^{-(s(\tbB)-2\lambda)w/2} \langle\bv, \bz\rangle)
      \end{pmatrix}^\top ,
      \qquad w \in \RR_+ , \quad \bz \in \RR^d .
 \]
First, we show
 \begin{equation}\label{bM2}
  \bQ(t) [\bM^{(2)}]_t \bQ(t)^\top
  - 2 w_{\bu,\bX_0}
    \sum_{\ell=1}^d
     c_\ell \langle\be_\ell, \tbu\rangle \int_0^t \Bf(w, \be_\ell) \, \dd w
  \as \bzero
  \qquad \text{as \ $t \to \infty$.}
 \end{equation}
For each \ $t, T \in \RR_+$, \ we have
 \[
   \bQ(t + T) [\bM^{(2)}]_{t+T} \bQ(t + T)^\top
   - 2 w_{\bu,\bX_0} \sum_{\ell=1}^d c_\ell \langle\be_\ell, \tbu\rangle
     \int_0^{t+T} \Bf(w, \be_\ell) \, \dd w
   = \bDelta_{t,T}^{(1)} + \bDelta_{t,T}^{(2)}
 \]
 with
 \begin{align*}
  \bDelta_{t,T}^{(1)}
  &:= 2 \sum_{\ell=1}^d
         c_\ell
         \int_0^T
          \Bf(t + T - \tau, \be_\ell)
          (\ee^{-s(\tbB)\tau} X_{\tau,\ell}
           - w_{\bu,\bX_0} \langle\be_\ell, \tbu\rangle)
          \, \dd \tau , \\
  \bDelta_{t,T}^{(2)}
  &:= 2 \sum_{\ell=1}^d
         c_\ell
         \int_T^{t+T}
          \Bf(t + T - \tau, \be_\ell)
          (\ee^{-s(\tbB)\tau} X_{\tau,\ell}
           - w_{\bu,\bX_0} \langle\be_\ell, \tbu\rangle)
          \, \dd \tau .
 \end{align*}
For each \ $t, T \in \RR_+$, \ we have
 \[
   \|\bDelta_{t,T}^{(1)}\|
   \leq 2 \biggl(\sup_{\tau\in[0,T]}
                  \|\ee^{-s(\tbB)\tau} \bX_\tau - w_{\bu,\bX_0} \tbu\|\biggr)
          \sum_{\ell=1}^d
           c_\ell \int_0^T \|\Bf(t + T - \tau, \be_\ell)\| \, \dd \tau ,
 \]
 where
 \ $\sup_{\tau\in[0,T]} \|\ee^{-s(\tbB)\tau} \bX_\tau - w_{\bu,\bX_0} \tbu\|
    < \infty$
 \ almost surely, since \ $(\bX_t)_{t\in\RR_+}$ \ has c\`adl\`ag sample paths (due to Theorem 4.6 in Barczy et al. \cite{BarLiPap1}).
 Then, using that \ $\|\bz \bz^\top\| \leq \|\bz\|^2$, \ $\bz \in \RR^2$, \ we have
 \begin{align}\label{help9}
  \begin{split}
  &\int_0^T \|\Bf(t + T - \tau, \be_\ell)\| \, \dd \tau
   = \int_t^{t+T} \|\Bf(w, \be_\ell)\| \, \dd w
   \leq \int_t^{t+T}
         |\ee^{-(s(\tbB)-2\lambda)w/2} \langle\bv, \be_\ell\rangle|^2 \, \dd w \\
  &\leq \|\bv\|^2 \int_t^{t+T} \ee^{-(s(\tbB)-2\Re(\lambda))w} \, \dd w
   \leq \|\bv\|^2 \int_t^\infty \ee^{-(s(\tbB)-2\Re(\lambda))w} \, \dd w \\
  &= \frac{\|\bv\|^2}{s(\tbB)-2\Re(\lambda)} \ee^{-(s(\tbB)-2\Re(\lambda))t}
   \to 0
  \end{split}
 \end{align}
 as \ $t \to \infty$.
\ Hence for each \ $T \in \RR_+$, \ we obtain
 \[
   \limsup_{t\to\infty} \|\bDelta_{t,T}^{(1)}\| = 0
 \]
 almost surely.
Moreover, for each \ $t, T \in \RR_+$, \ we have
 \[
   \|\bDelta_{t,T}^{(2)}\|
   \leq 2 \biggl(\sup_{\tau\in[T,\infty)}
                  \|\ee^{-s(\tbB)\tau} \bX_\tau - w_{\bu,\bX_0} \tbu\|\biggr)
          \sum_{\ell=1}^d
           c_\ell \int_T^{t+T} \|\Bf(t + T - \tau, \be_\ell)\| \, \dd \tau
 \]
 almost surely, where
 \begin{align}\label{help10}
  \begin{split}
   \int_T^{t+T} \|\Bf(t + T - \tau, \be_\ell)\| \, \dd \tau
   &= \int_0^t \|\Bf(w, \be_\ell)\| \, \dd w\\
   &\leq \|\bv\|^2 \int_0^\infty \ee^{-(s(\tbB)-2\Re(\lambda))w} \, \dd w
    = \frac{\|\bv\|^2}{s(\tbB)-2\Re(\lambda)} .
  \end{split}
 \end{align}
Consequently, for each \ $T \in \RR_+$, \ we obtain
 \begin{align*}
  &\limsup_{t\to\infty}
    \biggl\|\bQ(t) [\bM^{(2)}]_t \bQ(t)^\top
            - 2 w_{\bu,\bX_0} \sum_{\ell=1}^d c_\ell \langle\be_\ell, \tbu\rangle
              \int_0^t \Bf(w, \be_\ell) \, \dd w\biggr\| \\
  &= \limsup_{t\to\infty}
      \biggl\|\bQ(t + T) [\bM^{(2)}]_{t+T} \bQ(t + T)^\top
              - 2 w_{\bu,\bX_0} \sum_{\ell=1}^d c_\ell \langle\be_\ell, \tbu\rangle
                \int_0^{t+T} \Bf(w, \be_\ell) \, \dd w\biggr\| \\
  &\leq \limsup_{t\to\infty} \|\bDelta_{t,T}^{(1)}\|
        + \limsup_{t\to\infty} \|\bDelta_{t,T}^{(2)}\| \\
  &\leq \frac{2\|\bv\|^2}{s(\tbB)-2\Re(\lambda)}
        \biggl(\sup_{\tau\in[T,\infty)}
                  \|\ee^{-s(\tbB)\tau} \bX_\tau - w_{\bu,\bX_0} \tbu\|\biggr)
        \sum_{\ell=1}^d c_\ell
 \end{align*}
 almost surely.
Letting \ $T \to \infty$, \ by Theorem 3.3 in Barczy et al.\ \cite{BarPalPap}
 (which can be used, since the moment condition \eqref{moment_4_2} yields the moment condition
  \eqref{moment_condition_xlogx} with \ $\lambda = s(\tbB)$), \ we obtain \eqref{bM2}.
Moreover,
 \ $\int_0^t \Bf(w, \be_\ell) \, \dd w \to \int_0^\infty \Bf(w, \be_\ell) \, \dd w$
 \ as \ $t \to \infty$, \ since we have
 \[
   \int_0^\infty \|\Bf(w, \be_\ell)\| \, \dd w \\
   \leq \|\bv\|^2 \int_0^\infty \ee^{-(s(\tbB)-2\Re(\lambda))w} \, \dd w
   = \frac{\|\bv\|^2}{s(\tbB)-2\Re(\lambda)}
   < \infty .
 \]
Consequently,
 \begin{equation}\label{bM2+}
  \bQ(t) [\bM^{(2)}]_t \bQ(t)^\top
  \as 2 w_{\bu,\bX_0}
      \sum_{\ell=1}^d
       c_\ell \langle\be_\ell, \tbu\rangle \int_0^\infty \Bf(w, \be_\ell) \, \dd w
  \qquad \text{as \ $t \to \infty$.}
 \end{equation}
Next, by Theorem 3.3 in Barczy et al.\ \cite{BarPalPap}, we show that
 \[
   \bQ(t) [\bM^{(3,4)}]_t \bQ(t)^\top
   - w_{\bu,\bX_0}
     \sum_{\ell=1}^d
      \langle\be_\ell, \tbu\rangle
      \int_0^t \int_{\cU_d} \Bf(w,\bz) \, \dd w \, \mu_\ell(\dd\bz)
   \mean \bzero
 \]
 as \ $t \to \infty$.
\ Since
 \begin{equation}\label{QM34Q}
  \bQ(t) [\bM^{(3,4)}]_t \bQ(t)^\top
  = \sum_{\ell=1}^d
     \int_0^t \int_{\cU_d} \int_{\cU_1}
      \Bf(t - u,\bz) \ee^{-s(\tbB)u} \bbone_{\{w\leq X_{u-,\ell}\}}
       \, N_\ell(\dd u, \dd\bz, \dd w) ,
 \end{equation}
 it is enough to show that for each \ $\ell \in \{1, \ldots, d\}$ \ and
 \ $i, j \in \{1, 2\}$, \ we have
 \begin{align}\label{help11}
  \begin{split}
   &\int_0^t \int_{\cU_d} \int_{\cU_1}
     f_{i,j}(t - u, \bz) \ee^{-s(\tbB)u} \bbone_{\{w\leq X_{u-,\ell}\}}
        \, N_\ell(\dd u, \dd\bz, \dd w) \\
   &- w_{\bu,\bX_0}
      \langle\be_\ell, \tbu\rangle
      \int_0^t \int_{\cU_d} f_{i,j}(t - u, \bz) \, \dd u \, \mu_\ell(\dd\bz)
    \mean 0 \qquad \text{as \ $t \to \infty$,}
  \end{split}
 \end{align}
 where \ $\Bf(w,\bz) =: (f_{i,j}(w,\bz))_{i,j\in\{1,2\}}$, \ $w \in \RR_+$,
 \ $\bz \in \RR^d$.
\ For each \ $t \in \RR_+$, \ $\ell \in \{1, \ldots, d\}$ \ and
 \ $i, j \in \{1, 2\}$, \ we have
 \begin{equation}\label{III}
  \begin{aligned}
   &\EE\biggl(\biggl|\int_0^t \int_{\cU_d} \int_{\cU_1}
                      f_{i,j}(t - u,\bz) \ee^{-s(\tbB)u}
                      \bbone_{\{w\leq X_{u-,\ell}\}}
                      \, N_\ell(\dd u, \dd\bz, \dd w) \\
   &\phantom{\EE\biggl(\biggl|}
                     - w_{\bu,\bX_0} \langle\be_\ell, \tbu\rangle
                       \int_0^t \int_{\cU_d}
                        f_{i,j}(t - u, \bz)
                        \, \dd u \, \mu_\ell(\dd\bz)\biggr|\biggr)
    \leq I_{t,1} + I_{t,2} ,
  \end{aligned}
 \end{equation}
 where
 \begin{align*}
  I_{t,1}
  := \EE\biggl(\biggl|\int_0^t \int_{\cU_d} \int_{\cU_1}
                       f_{i,j}(t - u,\bz)
                       \ee^{-s(\tbB)u} \bbone_{\{w\leq X_{u-,\ell}\}}
                       \, \tN_\ell(\dd u, \dd\bz, \dd w)\biggr|\biggr)
 \end{align*}
 and
 \begin{align*}
  I_{t,2}
  &:= \EE\biggl(\biggl|\int_0^t \int_{\cU_d} \int_{\cU_1}
                        f_{i,j}(t - u,\bz)
                         \ee^{-s(\tbB)u} \bbone_{\{w\leq X_{u,\ell}\}}
                         \, \dd u \,\mu_\ell(\dd\bz)\, \dd w \\
   &\phantom{:= \EE\biggl(\biggl|}
                        - w_{\bu,\bX_0} \langle\be_\ell, \tbu\rangle
                          \int_0^t \int_{\cU_d}
                           f_{i,j}(t - u, \bz)
                           \, \dd u \, \mu_\ell(\dd\bz)\biggr|\biggr) \\
   &= \EE\biggl(\biggl|\int_0^t \int_{\cU_d}
                        f_{i,j}(t - u, \bz)
                        (\ee^{-s(\tbB)u}  X_{u,\ell}
                         - w_{\bu,\bX_0} \langle\be_\ell, \tbu\rangle)
                        \, \dd w \, \mu_\ell(\dd\bz)\biggr|\biggr) .
 \end{align*}
 Here, for each \ $\ell \in \{1, \ldots, d\}$ \ and \ $i, j \in \{1 ,2\}$, \ using
 Ikeda and Watanabe \cite[page 63]{IkeWat}, \eqref{EX} and that
 \ $|\Re(a)| \leq |a|$ \ and \ $|\Im(a)| \leq |a|$ \ for each \ $a \in \CC$, \ we
 have
 \begin{equation}\label{It1}
  \begin{aligned}
   I_{t,1}^2
   &\leq \EE\Biggl(\Biggl|
          \int_0^t \int_{\cU_d} \int_{\cU_1}
           f_{i,j}(t - u,\bz) \ee^{-s(\tbB)u} \bbone_{\{w\leq X_{u-,\ell}\}}
           \, \tN_\ell(\dd u, \dd\bz, \dd w)\Biggr|^2\Biggr) \\
   &= \int_0^t \int_{\cU_d}
       |f_{i,j}(t - u,\bz)|^2 \ee^{-2s(\tbB)u} \EE(X_{u,\ell})
       \, \dd u \, \mu_\ell( \dd\bz) \\
   &\leq \int_0^t \int_{\cU_d}
          |\ee^{-(s(\tbB) - 2\lambda)(t-u)/2} \langle \bv, \bz\rangle|^4
          \ee^{-2s(\tbB)u}  \EE(\|\bX_u\|)
          \, \dd u \, \mu_\ell(\dd\bz) \\
   &\leq C_4 \|\bv\|^4
         \int_{\cU_d} \|\bz\|^4 \, \mu_\ell(\dd\bz) \,
          \ee^{-2(s(\tbB) - 2\Re(\lambda))t}
         \int_0^t \ee^{(s(\tbB)-4\Re(\lambda))u} \, \dd u
    \to 0
  \end{aligned}
 \end{equation}
 as \ $t \to \infty$.
Indeed, if \ $s(\tbB) \ne 4\Re(\lambda)$, \
  using that \ $2\Re(\lambda)<s(\tbB)$,
\ we get
 \[
   I_{t,1}^2
   \leq C_4 \|\bv\|^4 \int_{\cU_d} \|\bz\|^4  \, \mu_\ell(\dd\bz)
        \frac{\ee^{-s(\tbB)t}-\ee^{-2(s(\tbB)-2\Re(\lambda))t}}
             {s(\tbB)-4\Re(\lambda)}
   \to 0
 \]
 as \ $t \to \infty$, \ since
 \ $\int_{\cU_d} \|\bz\|^4 \, \mu_\ell(\dd\bz) < \infty$.
\ Otherwise, if \ $s(\tbB) = 4 \Re(\lambda)$, \ then we obtain
 \[
   I_{t,1}^2
   \leq C_4 \|\bv\|^4
        \int_{\cU_d} \|\bz\|^4  \, \mu_\ell(\dd\bz) \, t \ee^{-s(\tbB)t}
   \to 0 \qquad \text{as \ $t \to \infty$.}
 \]
Further, for each \ $\ell\in\{1, \ldots, d\}$, \ $i, j \in\{1, 2\}$, \ and
 \ $t, T \in \RR_+$, \ we have
 \begin{equation}\label{ItT2}
  I_{t+T,2}  \leq  J_{t,T}^{(1)} + J_{t,T}^{(2)}
 \end{equation}
 with
 \begin{align*}
  J_{t,T}^{(1)}
  &:= \EE\biggl(\biggl|\int_0^T \int_{\cU_d}
       f_{i,j}(t + T - \tau, \bz)
       (\ee^{-s(\tbB)\tau} X_{\tau,\ell}
        - w_{\bu,\bX_0} \langle\be_\ell, \tbu\rangle)
       \, \dd \tau \, \mu_\ell(\dd\bz)\biggr|\biggr) , \\
  J_{t,T}^{(2)}
  &:= \EE\biggl(\biggl|\int_T^{t+T} \int_{\cU_d}
       f_{i,j}(t + T - \tau, \bz)
       (\ee^{-s(\tbB)\tau} X_{\tau,\ell}
        - w_{\bu,\bX_0} \langle\be_\ell, \tbu\rangle)
       \, \dd \tau \, \mu_\ell(\dd\bz)\biggr|\biggr) .
 \end{align*}
By Theorem 3.3 in Barczy et al.\ \cite{BarPalPap}, we have
 \ $K := \sup_{\tau\in\RR_+}
          \EE(\|\ee^{-s(\tbB)\tau} \bX_\tau - w_{\bu,\bX_0} \tbu\|)
      < \infty$,
 \ and hence, similarly as in \eqref{help9}, for any \ $T \in \RR_+$,
 \begin{align*}
  J_{t,T}^{(1)}
  &\leq K \int_0^T \int_{\cU_d}
           |f_{i,j}(t + T - \tau, \bz)| \, \dd \tau \, \mu_\ell(\dd\bz)
   = K \int_t^{t+T} \int_{\cU_d}
        |f_{i,j}(w, \bz)| \, \dd w \, \mu_\ell(\dd\bz) \\
  &\leq K \|\bv\|^2
        \int_{\cU_d} \|\bz\|^2 \, \mu_\ell(\dd\bz)
        \int_t^\infty \ee^{-(s(\tbB)-2\Re(\lambda))w} \, \dd w \\
  &= \frac{K\|\bv\|^2}{s(\tbB)-2\Re(\lambda)}
     \int_{\cU_d} \|\bz\|^2 \, \mu_\ell(\dd\bz)
     \, \ee^{-(s(\tbB) -2\Re(\lambda))t}
   \to 0
 \end{align*}
 as \ $t \to \infty$.
\ Further, similarly as in \eqref{help10}, for each \ $t, T \in \RR_+$,
 \begin{align*}
  J_{t,T}^{(2)}
  &\leq \sup_{\tau\in[T,\infty)}
         \EE(|\ee^{-s(\tbB)\tau} X_{\tau,\ell}
              - w_{\bu,\bX_0} \langle\be_\ell, \tbu\rangle|)
        \int_T^{t+T} \int_{\cU_d}
         |f_{i,j}(t + T - \tau, \bz)| \, \dd \tau \, \mu_\ell(\dd\bz) \\
  &\leq \sup_{\tau\in[T,\infty)}
         \EE(|\ee^{-s(\tbB)\tau} X_{\tau,\ell}
              - w_{\bu,\bX_0} \langle\be_\ell, \tbu\rangle|)
         \frac{\|\bv\|^2}{s(\tbB)-2\Re(\lambda)}
         \int_{\cU_d} \|\bz\|^2 \, \mu_\ell(\dd\bz) .
 \end{align*}
Consequently, for each \ $T \in \RR_+$, \ we obtain
 \begin{align*}
  \limsup_{t\to\infty} I_{t,2}
  &= \limsup_{t\to\infty} I_{t+T,2}
  \leq \limsup_{t\to\infty} J_{t,T}^{(1)}
       + \limsup_{t\to\infty} J_{t,T}^{(2)} \\
  &\leq \sup_{\tau\in[T,\infty)}
         \EE(|\ee^{-s(\tbB)\tau} X_{\tau,\ell}
              - w_{\bu,\bX_0} \langle\be_\ell, \tbu\rangle|)
         \frac{\|\bv\|^2}{s(\tbB)-2\Re(\lambda)}
         \int_{\cU_d} \|\bz\|^2 \, \mu_\ell(\dd\bz) .
 \end{align*}
Letting \ $T \to \infty$, \ by Theorem 3.3 in Barczy et al. \cite{BarPalPap}, we have
 \ $\lim_{t\to\infty} I_{t,2} = 0$, \ as desired.
All in all, \ $\lim_{t\to\infty} (I_{t,1} + I_{t,2}) = 0$, \ yielding
 \eqref{help11}.
Moreover, for each \ $\ell \in \{1, \ldots, d\}$,
 \[
   \int_0^t \int_{\cU_d} \Bf(t-u, \bz) \, \dd u \, \mu_\ell(\dd\bz)
   = \int_0^t \int_{\cU_d} \Bf(w, \bz) \, \dd w \, \mu_\ell(\dd\bz)
   \to \int_0^\infty \int_{\cU_d} \Bf(w, \bz) \, \dd w \, \mu_\ell(\dd\bz)
 \]
 as \ $t \to \infty$, \ since we have
 \begin{align*}
  \int_0^\infty \int_{\cU_d} \|\Bf(w, \bz)\| \, \dd w \, \mu_\ell(\dd\bz)
  &\leq \|\bv\|^2 \int_0^\infty \int_{\cU_d}
          \ee^{-(s(\tbB)-2\Re(\lambda))w} \|\bz\|^2 \, \dd w \, \mu_\ell(\dd\bz) \\
  &= \frac{\|\bv\|^2}{s(\tbB)-2\Re(\lambda)}
     \int_{\cU_d} \|\bz\|^2 \, \mu_\ell(\dd\bz)
   < \infty .
 \end{align*}
Consequently,
 \begin{equation}\label{bM34+}
  \bQ(t) [\bM^{(3,4)}]_t \bQ(t)^\top
  \mean w_{\bu,\bX_0}
      \sum_{\ell=1}^d
       \langle\be_\ell, \tbu\rangle
       \int_0^\infty \int_{\cU_d} \Bf(w,\bz) \, \dd w \, \mu_\ell(\dd\bz)
  \qquad \text{as \ $t \to \infty$.}
 \end{equation}
Further,
 \begin{align*}
  \bQ(t) [\bM^{(5)}]_t \bQ(t)^\top
  &= \int_0^t \int_{\cU_d} \Bf(t - u,\br) \ee^{-s(\tbB)u} \, M(\dd u, \dd\br)
  \mean \bzero
 \end{align*}
 as \ $t \to \infty$, \ since if \ $\Re(\lambda) \in (-\infty, \frac{1}{2}s(\tbB))$, \ then
 \begin{equation}\label{M5}
  \begin{aligned}
   \EE(\|\bQ(t) [\bM^{(5)}]_t \bQ(t)^\top\|)
   &\leq \EE\left(\int_0^t \int_{\cU_d}
                   \|\Bf(t - u, \br) \ee^{-s(\tbB)u}\|
                   \, M(\dd u, \dd\br) \right) \\
   &= \int_0^t \int_{\cU_d}
       \|\Bf(t - u,\br) \ee^{-s(\tbB)u} \| \, \dd u\, \nu(\dd\br) \\
   &\leq \int_0^t \int_{\cU_d}
          \bigl|\ee^{-(s(\tbB)-2\lambda)(t-u)/2} \langle\bv, \br\rangle\bigr|^2
          \ee^{-s(\tbB)u}
          \, \dd u \, \nu(\dd\br) \\
   &\leq \|\bv\|^2
         \ee^{-(s(\tbB)-2\Re(\lambda))t} \int_0^t \ee^{-2\Re(\lambda)u} \, \dd u
         \int_{\cU_d} \|\br\|^2 \, \nu(\dd\br)
    \to 0
  \end{aligned}
 \end{equation}
 as \ $t \to \infty$.
\ Indeed, if \ $\Re(\lambda) \ne 0$, \ then
 \[
   \ee^{-(s(\tbB)-2\Re(\lambda))t} \int_0^t \ee^{-2\Re(\lambda)u} \, \dd u
   = \frac{1}{2\Re(\lambda)}
     \Bigl(\ee^{-(s(\tbB)-2\Re(\lambda))t}  - \ee^{-s(\tbB)t}\Bigr)
     \to 0
 \]
 as \ $t \to \infty$, \ and if \ $\Re(\lambda) = 0$, \ then
  \ $\ee^{-(s(\tbB)-2\Re(\lambda))t} \int_0^t \ee^{-2\Re(\lambda)u} \, \dd u = t \ee^{-s(\tbB)t} \to 0$ \ as \ $t \to \infty$.
\ Consequently, by \eqref{bM2+} and \eqref{bM34+}, we get
 \begin{align}\label{help13}
  \begin{split}
   \bQ(t) [\bM]_t \bQ(t)^\top
   &\stoch 2 w_{\bu,\bX_0}
           \sum_{\ell=1}^d
            c_\ell \langle\be_\ell, \tbu\rangle
            \int_0^\infty \Bf(w, \be_\ell) \, \dd w \\
   &\phantom{\stoch\;}
           + w_{\bu,\bX_0}
             \sum_{\ell=1}^d
              \langle\be_\ell, \tbu\rangle
              \int_0^\infty \int_{\cU_d}
               \Bf(w, \bz) \, \dd w \, \mu_\ell(\dd\bz)
   = w_{\bu,\bX_0} \bSigma_\bv
  \end{split}
 \end{align}
 as \ $t \to \infty$, \ hence the condition \eqref{ltm_cond1} of Theorem \ref{THM_Cri_Pra} holds.
Indeed, for each \ $a \in \CC$, \ we have the identity
 \begin{equation}\label{help12_identity}
  \begin{aligned}
   \begin{pmatrix} \Re(a) \\ \Im(a) \end{pmatrix}
   \begin{pmatrix} \Re(a) \\ \Im(a) \end{pmatrix}^\top
   &= \begin{pmatrix}
       \Re(a)^2 & \Re(a) \Im(a) \\
       \Re(a) \Im(a) & \Im(a)^2
      \end{pmatrix} \\
   &= \frac{1}{2} |a|^2 \bI_2
      + \frac{1}{2}
        \begin{pmatrix}
         \Re(a^2) & \Im(a^2) \\
         \Im(a^2) & - \Re(a^2)
        \end{pmatrix} .
  \end{aligned}
 \end{equation}
Hence, for each \ $\ell \in \{1, \ldots, d\}$, \ applying \eqref{help12_identity}
 with \ $a=\ee^{-(s(\tbB)- 2\lambda)w/2} \langle\bv, \be_\ell\rangle$, \ we have
 \begin{align*}
  \int_0^\infty \Bf(w, \be_\ell) \, \dd w
  &= \frac{1}{2}
     \int_0^\infty
      \begin{pmatrix}
       \ee^{-(s(\tbB)-2\Re(\lambda))w} |\langle\bv, \be_\ell\rangle|^2 & 0 \\
       0 & \ee^{-(s(\tbB)-2\Re(\lambda))w} |\langle\bv, \be_\ell\rangle|^2 \\
      \end{pmatrix}
      \dd w
      \end{align*}
       \begin{align*}
  &\phantom{=\;}
     + \frac{1}{2}
       \int_0^\infty
        \begin{pmatrix}
         \Re(\ee^{-(s(\tbB)-2\lambda)w} \langle \bv, \be_\ell\rangle^2)
          & \Im(\ee^{-(s(\tbB)-2\lambda)w} \langle\bv, \be_\ell\rangle^2) \\
         \Im(\ee^{-(s(\tbB)-2\lambda)w} \langle\bv, \be_\ell\rangle^2)
          & -\Re(\ee^{-(s(\tbB)-2\lambda)w}\langle \bv,\be_\ell\rangle^2) \\
        \end{pmatrix}
        \dd w \\
  &= \frac{|\langle\bv, \be_\ell\rangle|^2}{2(s(\tbB)-2\Re(\lambda))} \bI_2 \\
  &\phantom{=\;}
     + \frac{1}{2}
       \begin{pmatrix}
        \Re(\int_0^\infty \ee^{-(s(\tbB)-2\lambda)w} \, \dd w \,
            \langle \bv ,\be_\ell\rangle^2)
         & \Im(\int_0^\infty \ee^{-(s(\tbB)-2\lambda)w} \, \dd w \,
              \langle\bv,\be_\ell\rangle^2) \\
        \Im(\int_0^\infty \ee^{-(s(\tbB)-2\lambda)w} \, \dd w \,
            \langle\bv,\be_\ell\rangle^2)
         & -\Re(\int_0^\infty \ee^{-(s(\tbB)-2\lambda)w} \, \dd w \,
                \langle \bv,\be_\ell\rangle^2)  \\
       \end{pmatrix} \\
  &= \frac{|\langle\bv, \be_\ell\rangle|^2}{2(s(\tbB)-2\Re(\lambda))} \bI_2
     + \frac{1}{2}
       \begin{pmatrix}
        \Re\left(\frac{\langle\bv,\be_\ell\rangle^2}{s(\tbB)-2\lambda}\right)
         & \Im\left(\frac{\langle\bv,\be_\ell\rangle^2}
                         {s(\tbB)-2\lambda}\right) \\[2mm]
        \Im\left(\frac{\langle\bv,\be_\ell\rangle^2}{s(\tbB)-2\lambda}\right)
         & -\Re\left(\frac{\langle\bv,\be_\ell\rangle^2}{s(\tbB)-2\lambda} \right)
       \end{pmatrix} ,
 \end{align*}
 and similarly
 \begin{align*}
  \int_0^\infty \int_{\cU_d} \Bf(w,\bz) \, \dd w \, \mu_\ell(\dd\bz) =
 \end{align*}
 \begin{align*}
  &= \frac{\int_{\cU_d}|\langle\bv,\bz \rangle|^2 \, \mu_\ell(\dd\bz)}
          {2(s(\tbB)-2\Re(\lambda))}\bI_2
     + \frac{1}{2}
       \begin{pmatrix}
        \Re\left(\frac{\int_{\cU_d}\langle\bv,\bz\rangle^2\,\mu_\ell(\dd\bz)}
                      {s(\tbB)-2\lambda}\right)
         & \Im\left(\frac{\int_{\cU_d}\langle\bv,\bz\rangle^2\,\mu_\ell(\dd\bz)}
                         {s(\tbB)-2\lambda}\right) \\[3mm]
        \Im\left(\frac{\int_{\cU_d}\langle\bv,\bz\rangle^2\,\mu_\ell(\dd\bz)}
                      {s(\tbB)-2\lambda}\right)
         & -\Re\left(\frac{\int_{\cU_d}\langle\bv,\bz\rangle^2\,\mu_\ell(\dd\bz)}
                          {s(\tbB)-2\lambda}\right)
       \end{pmatrix} ,
 \end{align*}
 yielding \eqref{help13}.
Note that \ $ \bSigma_\bv$ \ is non-negative definite irrespective of
 \ $\tBbeta \ne \bzero$ \ or \ $\tBbeta = \bzero$, \ since \ $\bc \in \RR_+^d$,
 \ $\tbu \in \RR_{++}^d$, \ and \ $\Bf(w,\bz)$ \ is non-negative definite for any
 \ $w \in \RR_+$ \ and \ $\bz \in \RR^d$.

\textbf{Step 3.}  Now we turn to prove that condition \eqref{ltm_cond2} of Theorem \ref{THM_Cri_Pra}
 holds for \ $(\bM_t)_{t\in\RR_+}$ \ with the scaling \ $\bQ(t)$, \ $t \in \RR_+$, \ namely,
 \[
  \EE\biggl(\sup_{u\in[0,t]} \|\bQ(t) (\bM_u - \bM_{u-})\|\biggr)
  \to 0 \qquad \text{as \ $t \to \infty$.}
 \]
Since \ $(\bM^{(2)}_t)_{t\in\RR_+}$ \ has continuous sample paths, we have for each
 \ $t \in \RR_+$,
 \begin{equation}\label{sup}
  \begin{aligned}
   &\sup_{u\in[0,t]} \|\bQ(t) (\bM_u - \bM_{u-})\|
    = \sup_{u\in[0,t]}
       \|\bQ(t) (\bM^{(3,4)}_u - \bM^{(3,4)}_{u-})
         + \bQ(t) (\bM^{(5)}_u - \bM^{(5)}_{u-})\|\\
   &\leq \|\bQ(t)\|
         \sum_{\ell=1}^d \sup_{u\in[0,t]} \|\tbY^{(\ell)}_u - \tbY^{(\ell)}_{u-}\|
         + \|\bQ(t)\| \sup_{u\in[0,t]} \|\bM^{(5)}_u - \bM^{(5)}_{u-}\|
  \end{aligned}
 \end{equation}
 almost surely.
Since \ $\bQ(t) \bQ(t)^\top = \ee^{-(s(\tbB)-2\Re(\lambda))t} \bI_2$,
 \ $t \in \RR_+$, \ we have \ $\|\bQ(t)\| = \ee^{-(s(\tbB)-2\Re(\lambda))t/2}$,
 \ $t \in \RR_+$.
\ Hence it is enough to show that
 \begin{align}\label{help_maximal_jump_1}
  \ee^{-(s(\tbB)-2\Re(\lambda))t/2}
  \EE\biggl(\sup_{u\in[0,t]} \|\tbY^{(\ell)}_u - \tbY^{(\ell)}_{u-}\|\biggr)
  \to 0 \qquad \text{as \ $t \to \infty$}
 \end{align}
 for every \ $\ell \in \{1, \ldots, d\}$, \ and
 \begin{align}\label{help_maximal_jump_2}
  \ee^{-(s(\tbB)-2\Re(\lambda))t/2}
  \EE\biggl(\sup_{u\in[0,t]} \|\bM^{(5)}_u - \bM^{(5)}_{u-}\|\biggr)
  \to 0 \qquad \text{as \ $t \to \infty$.}
 \end{align}
First, we prove \eqref{help_maximal_jump_1} for each \ $\ell \in \{1, \ldots, d\}$.
\  By Cauchy-Schwarz's inequality, for each \ $\vare \in \RR_{++}$, \ $t \in \RR_+$ \ and
  \ $\ell \in \{1, \ldots, d\}$, \ we have
 \begin{align}\label{CS1}
  \begin{split}
   &\ee^{-(s(\tbB)-2\Re(\lambda))t/2}
    \EE\biggl(\sup_{u\in[0,t]}
               \|\tbY^{(\ell)}_u
                 - \tbY^{(\ell)}_{u-}\|\biggr) \\
   &\leq \ee^{-(s(\tbB)-2\Re(\lambda))t/2}
        \EE\biggl(\sup_{u\in[0,t]}
                   \|\tbY^{(\ell,\vare)}_u
                     - \tbY^{(\ell,\vare)}_{u-}\|
           \biggr) \\
   &\phantom{\leq\;}
         + \ee^{-(s(\tbB)-2\Re(\lambda))t/2}
           \EE\biggl(\sup_{u\in[0,t]}
         \|(\tbY^{(\ell)}_u
            - \tbY^{(\ell,\vare)}_{u})
           - (\tbY^{(\ell)}_{u-}
              - \tbY^{(\ell,\vare)}_{u-})\|\biggr) \\
   &\leq \ee^{-(s(\tbB)-2\Re(\lambda))t/2}
         \Biggl(\EE\biggl(\sup_{u\in[0,t]}
           \|\tbY^{(\ell,\varepsilon)}_u
             - \tbY^{(\ell,\vare)}_{u-}\|^4\biggr)
             \Biggr)^{1/4} \\
   &\phantom{\leq\;}
         + \ee^{-(s(\tbB)-2\Re(\lambda))t/2}
         \Biggl(\EE\biggl(\sup_{u\in[0,t]}
          \|(\tbY^{(\ell)}_u
             - \tbY^{(\ell,\vare)}_{u} )
            - (\tbY^{(\ell)}_{u-}
               - \tbY^{(\ell,\vare)}_{u-})\|^2\biggr) \Biggr)^{1/2} .
 \end{split}
 \end{align}
Here, by \eqref{EX}, for each \ $\vare \in \RR_{++}$, \ $t \in \RR_+$ \ and \ $\ell \in \{1, \ldots, d\}$, \  we have
 \begin{equation}\label{CS1+}
  \begin{aligned}
   &\EE\biggl(\sup_{u\in[0,t]}
     \|\tbY^{(\ell,\vare)}_u
       - \tbY^{(\ell,\vare)}_{u-}\|^4\biggr)
    = \EE\biggl(\sup_{u\in[0,t]} |\tY^{(\ell,\vare)}_u - \tY^{(\ell,\vare)}_{u-}|^4\biggr)
   = \EE\biggl(\sup_{u\in[0,t]} |Y^{(\ell,\vare)}_u - Y^{(\ell,\vare)}_{u-}|^4\biggr)\\
  & \leq \EE\Biggl(\sum_{u\in[0,t]} |Y^{(\ell,\vare)}_u - Y^{(\ell,\vare)}_{u-}|^4\Biggr) \\
  &= \EE\biggl(\int_0^t \int_{\cU_d} \int_{\cU_1}
                |\ee^{-\lambda u}
                \langle\bv, \bz\rangle
                \bbone_{\{\|\bz\|\geq\vare\}}
                \bbone_{\{w\leq X_{u-,\ell}\}}|^4
                \, N_\ell(\dd u, \dd\bz, \dd w)\biggr) \\
  &= \int_0^t \int_{\cU_d}
      \ee^{-4\Re(\lambda)u} |\langle\bv, \bz\rangle|^4 \bbone_{\{\|\bz\|\geq\vare\}} \EE(X_{u,\ell})
      \, \dd u \, \mu_\ell(\dd\bz) \\
  &\leq C_4 \|\bv\|^4
        \int_0^t
         \ee^{(s(\tbB)-4\Re(\lambda))u} \, \dd u
        \int_{\cU_d}
         \|\bz\|^4 \bbone_{\{\|\bz\|\geq\vare\}}
         \, \mu_\ell(\dd\bz) .
  \end{aligned}
 \end{equation}
Hence, by \eqref{It1}  and \ $2\Re(\lambda)<s(\tbB)$,\   for each \
$\vare \in \RR_{++}$ \ and \ $\ell \in \{1, \ldots, d\}$, \ we get
 \begin{align}\label{help19}
   \ee^{-(s(\tbB)-2\Re(\lambda))t/2}
   \Biggl(\EE\biggl(\sup_{u\in[0,t]}
     \|\tbY^{(\ell,\vare)}_u
       - \tbY^{(\ell,\vare)}_{u-}\|^4\biggr)\Biggr)^{1/4}
   \to 0 \qquad \text{as \ $t \to \infty$.}
 \end{align}
Further, since
 \[
   \tbY_t^{(\ell)} - \tbY_t^{(\ell,\vare)}
   = \int_0^t \int_{\cU_d} \int_{\cU_1}
      \ee^{-\lambda u} \langle\bv, \bz\rangle
      \bbone_{\{\|\bz\|<\vare\}}
      \bbone_{\{w\leq X_{u-,\ell}\}}
      \, \tN_\ell(\dd u, \dd\bz, \dd w),
   \qquad t \in \RR_+ ,
 \]
 by the proof of part (a) of Theorem II.1.33 in Jacod and Shiryaev \cite{JSh}, we get
 \begin{align}\label{help17}
 \begin{split}
  &\EE\biggl(\sup_{u\in[0,t]}
      \|(\tbY^{(\ell)}_u - \tbY^{(\ell,\vare)}_{u})
        - (\tbY^{(\ell)}_{u-}
           - \tbY^{(\ell,\vare)}_{u-})\|^2\biggr) \\
  &\leq \EE\biggl(\sum_{u\in[0,t]}
         \|(\tbY^{(\ell)}_u - \tbY^{(\ell,\vare)}_{u})
           - (\tbY^{(\ell)}_{u-}
              - \tbY^{(\ell,\vare)}_{u-})\|^2\biggr) \\
  &= \EE\Biggl(\int_0^t \int_{\cU_d} \int_{\cU_1}
       |\ee^{-\lambda u}  \langle\bv, \bz\rangle|^2
       \bbone_{\{\|\bz\|<\vare\}}
       \bbone_{\{w\leq X_{u,\ell}\}}
       \, N_\ell(\dd u, \dd\bz, \dd w)\Biggr) \\
  &\leq C_4 \|\bv\|^2
        \int_0^t \int_{\cU_d}
         \ee^{(s(\tbB)-2\Re(\lambda)) u} \|\bz\|^2
         \bbone_{\{\|\bz\|<\vare\}}
         \, \mu_{\ell}(\dd \bz) \\
  &\leq C_4 \|\bv\|^2
        \frac{\ee^{(s(\tbB)-2\Re(\lambda))t}}
             {s(\tbB)-2\Re(\lambda)}
        \int_{\cU_d}
         \|\bz\|^2 \bbone_{\{\|\bz\|<\vare\}}
         \, \mu_{\ell}(\dd \bz) .
  \end{split}
 \end{align}
Hence, by \eqref{CS1}  and \eqref{help19}, for all \ $\vare \in \RR_{++}$ \ and \ $\ell \in \{1, \ldots, d\}$, \ we have
 \begin{align*}
  \limsup_{t\to\infty}
   \ee^{-(s(\tbB)-2\Re(\lambda))t/2}
  &\EE\biggl(\sup_{u\in[0,t]}
             \|\tbY^{(\ell)}_u
               - \tbY^{(\ell)}_{u-}\|\biggr) \\
  &\leq \Biggl(\frac{C_4\|\bv\|^2}
                    {s(\tbB)-2\Re(\lambda)}
               \int_{\cU_d}
                \|\bz\|^2 \bbone_{\{\|\bz\|<\vare\}}
                \, \mu_{\ell}(\dd \bz) \Biggr)^{1/2} ,
 \end{align*}
 which tends to 0 as \ $\vare \downarrow 0$ \ due to \eqref{moment_condition_CBI2}.
Hence we conclude \eqref{help_maximal_jump_1} for each \ $\ell \in \{1, \ldots, d\}$.

Next, we prove \eqref{help_maximal_jump_2}.
By Cauchy-Schwarz's inequality, for each \ $t \in \RR_+$, \ we have
 \begin{equation}\label{CS2}
  \EE\biggl(\sup_{u\in[0,t]}
             \|\bM^{(5)}_u - \bM^{(5)}_{u-}\|\biggr)
  \leq \Biggl(\EE\biggl(\sup_{u\in[0,t]}
                   \|\bM^{(5)}_u - \bM^{(5)}_{u-}\|^2\biggr)\!\Biggr)^{1/2}\!
   = \Biggl(\EE\biggl(\sup_{u\in[0,t]} |Z^{(5)}_u - Z^{(5)}_{u-}|^2\biggr)\Biggr)^{1/2} ,
 \end{equation}
 hence it is enough to prove that
 \[
   \ee^{-(s(\tbB)-2\Re(\lambda))t}
   \EE\biggl(\sup_{u\in[0,t]} |Z^{(5)}_u - Z^{(5)}_{u-}|^2\biggr)
   \to 0 \qquad \text{as \ $t \to \infty$.}
 \]
Since \ $\int_{\cU_d} \|\br\| \, \nu(\dd \br) < \infty$, \ for each \ $t \in \RR_+$, \ we have
 \ $Z^{(5)}_t = Z^*_t - \int_0^t \int_{\cU_d}
      \ee^{-\lambda u} \langle\bv, \br\rangle
      \, \dd u \, \nu(\dd\br)$
 with
 \ $Z^*_t := \int_0^t \int_{\cU_d}
              \ee^{-\lambda u} \langle\bv, \br\rangle \, M(\dd u, \dd\br)$,
 \ hence
 \begin{equation}\label{Z5}
  \begin{aligned}
  &\EE\biggl(\sup_{u\in[0,t]} |Z^{(5)}_u - Z^{(5)}_{u-}|^2\biggr)
   = \EE\biggl(\sup_{u\in[0,t]} |Z^*_u - Z^*_{u-}|^2\biggr)
   \leq \EE\Biggl(\sum_{u\in[0,t]} |Z^*_u - Z^*_{u-}|^2\Biggr) \\
  &= \EE\biggl(\int_0^t \int_{\cU_d}
                |\ee^{-\lambda u} \langle\bv, \br\rangle|^2
                \, M(\dd u, \dd\br)\biggr)
   = \int_0^t \int_{\cU_d}
      \ee^{-2\Re(\lambda)u} |\langle\bv, \br\rangle|^2
      \, \dd u \, \nu(\dd\br) \\
  &\leq \|\bv\|^2 \int_0^t \ee^{-2\Re(\lambda)u} \, \dd u
        \int_{\cU_d} \|\br\|^2 \, \nu(\dd\br)
  \end{aligned}
 \end{equation}
 hence, by \eqref{M5}, we conclude \eqref{help_maximal_jump_2}.
Consequently, by Theorem \ref{THM_Cri_Pra}, we obtain
 \[
   \bQ(t) \bM_t \distr (w_{\bu,\bX_0} \bSigma_\bv)^{1/2} \bN \qquad
   \text{as \ $t \to \infty$,}
 \]
 where \ $\bN$ \ is a \ $2$-dimensional random vector with
 \ $\bN \distre \cN_2(\bzero, \bI_2)$ \ independent of
 \ $w_{\bu,\bX_0} \bSigma_\bv$.
\ Clearly,
 \ $(w_{\bu,\bX_0} \bSigma_\bv)^{1/2} \bN
    = \sqrt{w_{\bu,\bX_0}} \, \bSigma_\bv^{1/2} \bN
    \distre \sqrt{w_{\bu,\bX_0}} \bZ_\bv$.
\ By the decomposition
 \[
   \ee^{-s(\tbB)t/2}
   \begin{pmatrix}
    \Re(\langle \bv, \bX_t \rangle) \\
    \Im(\langle \bv, \bX_t \rangle)
   \end{pmatrix}
   = \begin{pmatrix}
      \Re(\ee^{-(s(\tbB)-2\lambda)t/2} Z_t^{(0,1)}) \\
      \Im(\ee^{-(s(\tbB)-2\lambda)t/2} Z_t^{(0,1)})
     \end{pmatrix}
     + \bQ(t) \bM_t , \qquad t \in \RR_+ ,
 \]
 the convergence \eqref{Z01} and Slutsky's lemma  (see, e.g., van der Vaart \cite[Lemma 2.8]{vanderVaart}), we obtain \eqref{convvweak1}.
\proofend

\noindent
\textbf{Proof of part (ii) of Theorem \ref{convCBIweak1}.}
We use a similar approach as in the proof of part (iii) of Theorem \ref{convCBIweak1}.
We divide the proof into three main steps.

\textbf{Step 1.}
We use the same representation of \ $\ee^{-\lambda t} \langle\bv, \bX_t\rangle$,
 \ $t \in \RR_+$ \ as in the proof of part (iii) of Theorem \ref{convCBIweak1}.
We have
 \begin{equation}\label{Z01m}
  t^{-1/2} \ee^{-(s(\tbB)-2\lambda)t/2} Z_t^{(0,1)} \as 0 \qquad
  \text{as \ $t \to \infty$,}
 \end{equation}
 since \ $\Re(\lambda) = \frac{1}{2} s(\tbB) > 0$ \ implies \ $\lambda \ne 0$,
 \ hence
 \begin{align*}
  t^{-1/2} \ee^{-(s(\tbB)-2\lambda)t/2} Z_t^{(0,1)}
  &= t^{-1/2} \ee^{\ii\Im(\lambda)t}
     \biggl(\langle\bv, \bX_0\rangle
            - \frac{\langle\bv, \tBbeta\rangle}{\lambda}
              (\ee^{-\lambda t} - 1)\biggr) \\
  &= t^{-1/2} \ee^{\ii\Im(\lambda)t}
     \biggl(\langle\bv, \bX_0\rangle
            + \frac{\langle\bv, \tBbeta\rangle}{\lambda}\biggr)
     - t^{-1/2} \ee^{-s(\tbB)t/2} \frac{\langle\bv, \tBbeta\rangle}{\lambda}
   \as 0
 \end{align*}
 as \ $t \to \infty$.

For each \ $t \in \RR_+$, \ with the notations of the proof of part (iii) of Theorem
 \ref{convCBIweak1}, we have
 \[
   \begin{pmatrix}
    \Re\bigl(t^{-1/2}\ee^{-(s(\tbB)-2\lambda)t/2}
             \bigl(Z_t^{(2)} + Z_t^{(3,4)} + Z_t^{(5)}\bigr)\bigr) \\
    \Im\bigl(t^{-1/2}\ee^{-(s(\tbB)-2\lambda)t/2}
             \bigl(Z_t^{(2)} + Z_t^{(3,4)} + Z_t^{(5)}\bigr)\bigr)
   \end{pmatrix}
   = t^{-1/2} \bQ(t) \bM_t ,
 \]
 where now
 \[
   \bQ(t)
   = \begin{pmatrix}
      \Re(\ee^{\ii\Im(\lambda)t}) & -\Im(\ee^{\ii\Im(\lambda)t}) \\
      \Im(\ee^{\ii\Im(\lambda)t}) & \Re(\ee^{\ii\Im(\lambda)t})
     \end{pmatrix}
   = \begin{pmatrix}
      \cos(\Im(\lambda)t) & -\sin(\Im(\lambda)) \\
      \sin(\Im(\lambda)t) & \cos(\Im(\lambda)t)
     \end{pmatrix} , \qquad t \in \RR_+ .
 \]
We are again going to apply Theorem \ref{THM_Cri_Pra} for the 2-dimensional
 martingale \ $(\bM_t)_{t\in\RR_+}$ \ now with the scaling \ $t^{-1/2} \bQ(t)$,
 \ $t \in \RR_+$.
\ We clearly have \ $t^{-1/2} \bQ(t) \to \bzero$ \ as \ $t \to \infty$.

\textbf{Step 2.}  Now we prove that condition \eqref{ltm_cond1} of Theorem \ref{THM_Cri_Pra}
holds.  For each \ $t \in \RR_+$, \ with the notations of the proof of part (iii) of Theorem
 \ref{convCBIweak1}, we have
 \[
   t^{-1} \bQ(t) [\bM^{(2)}]_t \bQ(t)^\top
   = \frac{2}{t}
     \sum_{\ell=1}^d
      c_\ell \int_0^t \Bf(t - \tau, \be_\ell) \ee^{-s(\tbB)\tau} X_{\tau,\ell}
      \, \dd \tau ,
 \]
 where now
 \[
   \Bf(w, \bz)
   = \begin{pmatrix}
      \Re(\ee^{\ii\Im(\lambda)w} \langle\bv, \bz\rangle) \\
      \Im(\ee^{\ii\Im(\lambda)w} \langle\bv, \bz\rangle)
     \end{pmatrix}
     \begin{pmatrix}
      \Re(\ee^{\ii\Im(\lambda)w} \langle\bv, \bz\rangle) \\
      \Im(\ee^{\ii\Im(\lambda)w} \langle\bv, \bz\rangle)
     \end{pmatrix}^\top ,
   \qquad w \in \RR_+ , \quad \bz \in \RR^d .
 \]
First, we show
 \begin{equation}\label{bM2m}
  t^{-1} \bQ(t) [\bM^{(2)}]_t \bQ(t)^\top
  - \frac{2w_{\bu,\bX_0}}{t}
    \sum_{\ell=1}^d c_\ell \langle\be_\ell, \tbu\rangle \int_0^t \Bf(w, \be_\ell) \, \dd w
  \as \bzero \qquad \text{as \ $t \to \infty$.}
 \end{equation}
For each \ $t, T \in \RR_{++}$, \ we have
 \[
   (t+T)^{-1} \bQ(t + T) [\bM^{(2)}]_{t+T} \bQ(t + T)^\top
   - \frac{2w_{\bu,\bX_0}}{t+T}
     \sum_{\ell=1}^d c_\ell \langle\be_\ell, \tbu\rangle \int_0^{t+T} \Bf(w, \be_\ell) \, \dd w
   = \bDelta_{t,T}^{(1)} + \bDelta_{t,T}^{(2)}
 \]
 with
 \begin{align*}
  \bDelta_{t,T}^{(1)}
  &:= \frac{2}{t+T}
      \sum_{\ell=1}^d
        c_\ell
        \int_0^T
         \Bf(t + T - \tau, \be_\ell)
         (\ee^{-s(\tbB)\tau} X_{\tau,\ell}
          - w_{\bu,\bX_0} \langle\be_\ell, \tbu\rangle)
         \, \dd \tau , \\
  \bDelta_{t,T}^{(2)}
  &:= \frac{2}{t+T}
      \sum_{\ell=1}^d
       c_\ell
       \int_T^{t+T}
        \Bf(t + T - \tau, \be_\ell)
        (\ee^{-s(\tbB)\tau} X_{\tau,\ell}
         - w_{\bu,\bX_0} \langle\be_\ell, \tbu\rangle)
        \, \dd \tau .
 \end{align*}
For each \ $t, T \in \RR_+$, \ we have
 \[
   \|\bDelta_{t,T}^{(1)}\|
   \leq \frac{2}{t+T}
        \biggl(\sup_{\tau\in[0,T]}
                \|\ee^{-s(\tbB)\tau} \bX_\tau - w_{\bu,\bX_0} \tbu\|\biggr)
        \sum_{\ell=1}^d
          c_\ell \int_0^T \|\Bf(t + T - \tau, \be_\ell)\| \, \dd \tau ,
 \]
 where
 \ $\sup_{\tau\in[0,T]} \|\ee^{-s(\tbB)\tau} \bX_\tau - w_{\bu,\bX_0} \tbu\|
    < \infty$
 \ almost surely since \ $(\bX_t)_{t\in\RR_+}$ \ has c\`adl\`ag sample paths, and
 using that \ $\|\bz \bz^\top\| \leq \|\bz\|^2$, \ $\bz \in \RR^2$, \ we have
 \[
   \int_0^T \|\Bf(t + T - \tau, \be_\ell)\| \, \dd \tau
   = \int_t^{t+T} \|\Bf(w, \be_\ell)\| \, \dd w
   \leq \int_t^{t+T}
         |\ee^{\ii\Im(\lambda)w} \langle\bv, \be_\ell\rangle|^2 \, \dd w
   \leq \|\bv\|^2 T .
 \]
Hence for each \ $T \in \RR_+$, \ we obtain
 \[
   \limsup_{t\to\infty} \|\bDelta_{t,T}^{(1)}\| = 0
 \]
 almost surely.
Moreover, for each \ $t, T \in \RR_+$, \ we have
 \[
   \|\bDelta_{t,T}^{(2)}\|
   \leq \frac{2}{t+T}
        \biggl(\sup_{\tau\in[T,\infty)}
                \|\ee^{-s(\tbB)\tau} \bX_\tau - w_{\bu,\bX_0} \tbu\|\biggr)
        \sum_{\ell=1}^d
         c_\ell \int_T^{t+T} \|\Bf(t + T - \tau, \be_\ell)\| \, \dd \tau
 \]
 almost surely, where
 \[
   \int_T^{t+T} \|\Bf(t + T - \tau, \be_\ell)\| \, \dd \tau
   = \int_0^t \|\Bf(w, \be_\ell)\| \, \dd w\\
   \leq \int_0^t |\ee^{\ii\Im(\lambda)w} \langle\bv, \be_\ell\rangle|^2 \, \dd w
   \leq \|\bv\|^2 t .
 \]
 Consequently, for each \ $T \in \RR_+$, \ we obtain
 \begin{align*}
  &\limsup_{t\to\infty}
    \biggl\|t^{-1} \bQ(t) [\bM^{(2)}]_t \bQ(t)^\top
            - \frac{2w_{\bu,\bX_0}}{t}
              \sum_{\ell=1}^d c_\ell \langle\be_\ell, \tbu\rangle
              \int_0^t \Bf(w, \be_\ell) \, \dd w\biggr\| \\
  &= \limsup_{t\to\infty}
      \biggl\|(t+T)^{-1} \bQ(t + T) [\bM^{(2)}]_{t+T} \bQ(t + T)^\top
              - \frac{2w_{\bu,\bX_0}}{t+T}
                \sum_{\ell=1}^d c_\ell \langle\be_\ell, \tbu\rangle
                \int_0^{t+T} \Bf(w, \be_\ell) \, \dd w\biggr\| \\
  &\leq \limsup_{t\to\infty} \|\bDelta_{t,T}^{(1)}\|
        + \limsup_{t\to\infty} \|\bDelta_{t,T}^{(2)}\| \\
  &\leq 2 \|\bv\|^2
        \biggl(\sup_{\tau\in[T,\infty)}
                  \|\ee^{-s(\tbB)\tau} \bX_\tau - w_{\bu,\bX_0} \tbu\|\biggr)
        \sum_{\ell=1}^d c_\ell
 \end{align*}
 almost surely.
Letting \ $T \to \infty$, \ by Theorem 3.3 in Barczy et al.\ \cite{BarPalPap}, we
 obtain \eqref{bM2m}.
The aim of the following discussion is to show
 \begin{equation}\label{int_mean}
  \frac{1}{t} \int_0^t \Bf(w, \be_\ell) \, \dd w
  \to \frac{1}{2} |\langle\bv, \be_\ell\rangle|^2 \bI_2
      + \frac{1}{2}
        \begin{pmatrix}
         \Re(\langle\bv, \be_\ell\rangle^2) & \Im(\langle\bv, \be_\ell\rangle^2) \\
         \Im(\langle\bv, \be_\ell\rangle^2) & -\Re(\langle\bv, \be_\ell\rangle^2)
        \end{pmatrix}
        \bbone_{\{\Im(\lambda)=0\}}
 \end{equation}
 as \ $t \to \infty$.
Applying \eqref{help12_identity} for
 \ $a = \ee^{\ii\Im(\lambda)w} \langle\bv, \be_\ell\rangle$, \ we obtain
 \[
   \Bf(w, \be_\ell)
   = \frac{1}{2} |\langle\bv, \be_\ell\rangle|^2 \bI_2
      + \frac{1}{2}
        \begin{pmatrix}
         \Re((\ee^{\ii\Im(\lambda)w} \langle\bv, \be_\ell\rangle)^2)
          & \Im((\ee^{\ii\Im(\lambda)w} \langle\bv, \be_\ell\rangle)^2) \\
         \Im((\ee^{\ii\Im(\lambda)w} \langle\bv, \be_\ell\rangle)^2)
          & - \Re((\ee^{\ii\Im(\lambda)w} \langle\bv, \be_\ell\rangle)^2)
        \end{pmatrix} .
 \]
Thus, if \ $\Im(\lambda) = 0$, \ then we have
 \[
   \frac{1}{t} \int_0^t \Bf(w, \be_\ell) \, \dd w
   = \frac{1}{2} |\langle\bv, \be_\ell\rangle|^2 \bI_2
     + \frac{1}{2}
       \begin{pmatrix}
        \Re(\langle\bv, \be_\ell\rangle)^2)
         & \Im(\langle\bv, \be_\ell\rangle)^2) \\
        \Im(\langle\bv, \be_\ell\rangle)^2)
         & - \Re(\langle\bv, \be_\ell\rangle)^2)
       \end{pmatrix}
 \]
 for all \ $t \in \RR_+$.
\ If \ $\Im(\lambda) \ne 0$, \ then we have
 \begin{align*}
  \frac{1}{t}
  \int_0^t \Re((\ee^{\ii\Im(\lambda)w} \langle\bv, \be_\ell\rangle)^2) \, \dd w
  &= \frac{1}{t}
     \Re\biggl(\langle\bv, \be_\ell\rangle^2
               \int_0^t \ee^{2\ii\Im(\lambda)w} \, \dd w\biggr) \\
  &= \frac{1}{t}
     \Re\biggl(\frac{\langle\bv, \be_\ell\rangle^2}{2\ii\Im(\lambda)}
               (\ee^{2\ii\Im(\lambda)t} - 1)\biggr)
   \leq \frac{\|\bv\|^2}{|\Im(\lambda)|t}
   \to 0
 \end{align*}
 as \ $t \to \infty$, \ and, in a similar way,
 \ $\frac{1}{t}
    \int_0^t \Im((\ee^{\ii\Im(\lambda)w} \langle\bv, \be_\ell\rangle)^2) \, \dd t
    \to 0$
 \ as \ $t \to \infty$.
\ Hence
 \ $\frac{1}{t} \int_0^t \Bf(w, \be_\ell) \, \dd w
    \to \frac{1}{2} |\langle\bv, \be_\ell\rangle|^2 \bI_2$
 \ as \ $t \to \infty$, \ and we conclude \eqref{int_mean}.

Next, using Theorem 3.3 in Barczy et al.\ \cite{BarPalPap}, we show that
 \begin{equation}\label{bM34m}
   t^{-1} \bQ(t) [\bM^{(3,4)}]_t \bQ(t)^\top
   - t^{-1} w_{\bu,\bX_0}
     \sum_{\ell=1}^d
      \langle\be_\ell, \tbu\rangle
      \int_0^t \int_{\cU_d} \Bf(w,\bz) \, \dd w \, \mu_\ell(\dd\bz)
   \mean \bzero
 \end{equation}
 as \ $t \to \infty$.
\ By the help of \eqref{QM34Q}, it is enough to show that for each
 \ $\ell \in \{1, \ldots, d\}$ \ and \ $i, j \in \{1, 2\}$, \ we have
 \begin{align}\label{help11m}
  \begin{split}
   &t^{-1} \int_0^t \int_{\cU_d} \int_{\cU_1}
     f_{i,j}(t - u, \bz) \ee^{-s(\tbB)u} \bbone_{\{w\leq X_{u-,\ell}\}}
        \, N_\ell(\dd u, \dd\bz, \dd w) \\
   &- t^{-1} w_{\bu,\bX_0}
      \langle\be_\ell, \tbu\rangle
      \int_0^t \int_{\cU_d} f_{i,j}(t - u, \bz) \, \dd u \, \mu_\ell(\dd\bz)
    \mean 0 \qquad \text{as \ $t \to \infty$.}
  \end{split}
 \end{align}
For each \ $t \in \RR_+$, \ $\ell \in \{1, \ldots, d\}$ \ and \ $i, j \in \{1, 2\}$,
 \ we use again the estimation \eqref{III}.
For each \ $\ell \in \{1, \ldots, d\}$ \ and \ $i, j \in \{1, 2\}$, \ as in
 \eqref{It1}, we have
  \begin{align*}
   ((t + T)^{-1} I_{t,1})^2
   &\leq (t + T)^{-2}
         \int_0^t \int_{\cU_d}
          |\ee^{\ii\Im(\lambda)(t-u)} \langle \bv, \bz\rangle|^4
          \ee^{-2s(\tbB)u}  \EE(\|\bX_u\|)
          \, \dd u \, \mu_\ell(\dd\bz) \\
   &\leq C_4 \|\bv\|^4 (t + T)^{-2}
         \int_{\cU_d} \|\bz\|^4 \, \mu_\ell(\dd\bz)
         \int_0^t \ee^{-s(\tbB)u} \, \dd u
    \to 0
  \end{align*}
 as \ $t \to \infty$, \ since
 \ $\int_0^t \ee^{-s(\tbB)u} \, \dd u \leq \int_0^\infty \ee^{-s(\tbB)u} \, \dd u
    = \frac{1}{s(\tbB)}$
 for every \ $t \in \RR_+$, \ and
 \ $\int_{\cU_d} \|\bz\|^4 \, \mu_\ell(\dd\bz) < \infty$.
For each \ $t \in \RR_+$ \ and \ $T \in \RR_+$, \ we use again the decomposition
 \eqref{ItT2}.
Similarly as in \eqref{help9}, for any \ $T \in \RR_+$,
 \begin{align*}
  (t + T)^{-1}  J_{t,T}^{(1)}
  &\leq \frac{K}{t+T}
        \int_0^T \int_{\cU_d}
         |f_{i,j}(t + T - \tau, \bz)| \, \dd\tau \, \mu_\ell(\dd\bz) \\
  &\leq \frac{K}{t+T}
        \int_0^T \int_{\cU_d}
         |\langle\bv, \bz\rangle|^2 \, \dd \tau \, \mu_\ell(\dd\bz)
   \leq \frac{K\|\bv\|^2T}{t+T} \int_{\cU_d} \|\bz\|^2 \, \mu_\ell(\dd\bz)
   \to 0
 \end{align*}
 as \ $t \to \infty$.
\ Further, similarly as in \eqref{help10}, for each \ $t, T \in \RR_+$,
 \begin{align*}
  \frac{J_{t,T}^{(2)}}{t+T}
  &\leq \frac{1}{t+T}
        \sup_{\tau\in[T,\infty)}
         \EE(|\ee^{-s(\tbB)\tau} X_{\tau,\ell}
              - w_{\bu,\bX_0} \langle\be_\ell, \tbu\rangle|)
        \int_T^{t+T} \int_{\cU_d}
         |f_{i,j}(t + T - \tau, \bz)| \, \dd \tau \, \mu_\ell(\dd\bz) \\
  &\leq \frac{\|\bv\|^2t}{t+T}
        \sup_{\tau\in[T,\infty)}
         \EE(|\ee^{-s(\tbB)\tau} X_{\tau,\ell}
              - w_{\bu,\bX_0} \langle\be_\ell, \tbu\rangle|)
        \int_{\cU_d} \|\bz\|^2 \, \mu_\ell(\dd\bz).
 \end{align*}
Consequently, for each \ $T \in \RR_+$, \ we obtain
 \begin{align*}
  \limsup_{t\to\infty} t^{-1} I_{t,2}
  &= \limsup_{t\to\infty} (t + T)^{-1} I_{t+T,2}
  \leq \limsup_{t\to\infty} (t + T)^{-1} J_{t,T}^{(1)}
       + \limsup_{t\to\infty} (t + T)^{-1} J_{t,T}^{(2)} \\
  &\leq \|\bv\|^2
        \sup_{\tau\in[T,\infty)}
         \EE(|\ee^{-s(\tbB)\tau} X_{\tau,\ell}
              - w_{\bu,\bX_0} \langle\be_\ell, \tbu\rangle|)
        \int_{\cU_d} \|\bz\|^2 \, \mu_\ell(\dd\bz) .
 \end{align*}
Letting \ $T \to \infty$, \ by Theorem 3.3 in Barczy et al.\ \cite{BarPalPap}, we have
 \ $\lim_{t\to\infty} t^{-1} I_{t,2} = 0$, \ as desired.
All in all, \ $\lim_{t\to\infty} t^{-1} (I_{t,1} + I_{t,2}) = 0$, \ yielding
 \eqref{help11m}.
As in case of \eqref{int_mean}, one can derive
 \begin{equation}\label{int_mean34}
  \begin{aligned}
   &\frac{1}{t} \int_0^t \int_{\cU_d} \Bf(w, \bz) \, \dd w \, \mu_\ell(\bz)
    \to \frac{1}{2}
        \int_{\cU_d} |\langle\bv, \bz\rangle|^2 \, \mu_\ell(\dd\bz) \bI_2 \\
   &\hspace*{35mm}
        + \frac{1}{2}
          \begin{pmatrix}
           \Re\bigl(\int_{\cU_d} \langle\bv, \bz\rangle^2 \, \mu_\ell(\dd\bz)\bigr)
            & \Im\bigl(\int_{\cU_d}
                        \langle\bv, \bz\rangle^2 \, \mu_\ell(\dd\bz)\bigr) \\
           \Im\bigl(\int_{\cU_d} \langle\bv, \bz\rangle^2 \, \mu_\ell(\dd\bz)\bigr)
            & -\Re\bigl(\int_{\cU_d}
                         \langle\bv, \bz\rangle^2 \, \mu_\ell(\dd\bz)\bigr)
          \end{pmatrix}
          \bbone_{\{\Im(\lambda)=0\}}
  \end{aligned}
 \end{equation}
 as \ $t \to \infty$.
\ Indeed, we can apply \eqref{help12_identity} for
 \ $a = \ee^{\ii\Im(\lambda)w} \langle\bv, \bz\rangle$.
\ In case of \ $\Im(\lambda) = 0$, \ we obtain
 \begin{align*}
  \frac{1}{t} \int_0^t \int_{\cU_d} \Bf(w, \bz) \, \dd w \, \mu_\ell(\dd\bz)
  &= \frac{1}{2}
     \int_{\cU_d} |\langle\bv, \bz\rangle|^2 \, \mu_\ell(\dd\bz) \bI_2 \\
  &\quad
     + \frac{1}{2}
       \begin{pmatrix}
        \Re\bigl(\int_{\cU_d} \langle\bv, \bz\rangle^2 \, \mu_\ell(\dd\bz)\bigr)
         & \Im\bigl(\int_{\cU_d}
                     \langle\bv, \bz\rangle^2 \, \mu_\ell(\dd\bz)\bigr) \\
        \Im\bigl(\int_{\cU_d} \langle\bv, \bz\rangle^2 \, \mu_\ell(\dd\bz)\bigr)
         & - \Re\bigl(\int_{\cU_d}
                       \langle\bv, \bz\rangle^2 \, \mu_\ell(\dd\bz)\bigr)
       \end{pmatrix}
 \end{align*}
 for all \ $t \in \RR_+$.
\ If \ $\Im(\lambda) \ne 0$, \ then we have
 \begin{align*}
  &\frac{1}{t}
   \int_0^t \int_{\cU_d}
    \Re((\ee^{\ii\Im(\lambda)w} \langle\bv, \bz\rangle)^2)
    \, \dd w \, \mu_\ell(\dd\bz)
   = \frac{1}{t}
     \Re\biggl(\int_{\cU_d} \langle\bv, \bz\rangle^2 \, \mu_\ell(\dd\bz)
               \int_0^t \ee^{2\ii\Im(\lambda)w} \, \dd w\biggr) \\
  &= \frac{1}{t}
     \Re\biggl(\int_{\cU_d} \langle\bv, \bz\rangle^2 \, \mu_\ell(\dd\bz)
               \frac{\ee^{2\ii\Im(\lambda)t} - 1}{2\ii\Im(\lambda)}\biggr)
   \leq \frac{\|\bv\|^2}{|\Im(\lambda)|t} \int_{\cU_d} \|\bz\|^2 \, \mu_\ell(\dd\bz)
   \to 0
 \end{align*}
 as \ $t \to \infty$, \ and, in a similar way,
 \ $\frac{1}{t}
    \int_0^t \int_{\cU_d} \Im((\ee^{\ii\Im(\lambda)w} \langle\bv, \bz\rangle)^2)
     \, \dd w \, \mu_\ell(\dd\bz)
    \to 0$
 \ as \ $t \to \infty$.
\ Hence
 \ $\frac{1}{t} \int_0^t \int_{\cU_d} \Bf(w, \bz) \, \dd w \, \mu_\ell(\dd\bz)
    \to \frac{1}{2}
        \int_{\cU_d} |\langle\bv, \bz\rangle|^2 \, \mu_\ell(\dd\bz) \bI_2$
 \ as \ $t \to \infty$, \ and we conclude \eqref{int_mean34}.

Further,
 \[
   t^{-1} \bQ(t) [\bM^{(5)}]_t \bQ(t)^\top
   = t^{-1} \int_0^t \int_{\cU_d} \Bf(t - u,\br) \ee^{-s(\tbB)u} \, M(\dd u, \dd\br)
   \mean \bzero
 \]
 as \ $t \to \infty$, \ since \ $\Re(\lambda) = \frac{1}{2}s(\tbB) > 0$
 \ implies \ $\Re(\lambda) \ne 0$, \ and hence
 \begin{equation}\label{M5m}
  \begin{aligned}
   t^{-1} \EE(\|\bQ(t) [\bM^{(5)}]_t \bQ(t)^\top\|)
   &\leq t^{-1}
         \EE\left(\int_0^t \int_{\cU_d}
                   \|\Bf(t - u, \br) \ee^{-s(\tbB)u}\|
                   \, M(\dd u, \dd\br) \right) \\
   &= t^{-1}
      \int_0^t \int_{\cU_d}
       \|\Bf(t - u,\br) \ee^{-s(\tbB)u} \| \, \dd u\, \nu(\dd\br) \\
   &\leq t^{-1}
         \int_0^t \int_{\cU_d}
          \bigl|\langle\bv, \br\rangle\bigr|^2 \ee^{-s(\tbB)u}
          \, \dd u \, \nu(\dd\br) \\
   &\leq \|\bv\|^2 t^{-1}
         \int_0^t \ee^{-s(\tbB)u} \, \dd u
         \int_{\cU_d} \|\br\|^2 \, \nu(\dd\br)
    \to 0
  \end{aligned}
 \end{equation}
 as \ $t \to \infty$.
\ Consequently, by \eqref{bM2m}, \eqref{int_mean}, \eqref{bM34m}, \eqref{int_mean34}
 and \eqref{M5m}, we get
 \[
   t^{-1} \bQ(t) [\bM]_t \bQ(t)^\top \stoch w_{\bu,\bX_0} \bSigma_\bv
   \qquad \text{as \ $t \to \infty$.}
 \]

\textbf{Step 3.}  Now we turn to prove that condition \eqref{ltm_cond2} of Theorem \ref{THM_Cri_Pra}
 holds, namely,
 \[
  \EE\biggl(\sup_{u\in[0,t]} t^{-1/2} \|\bQ(t) (\bM_u - \bM_{u-})\|\biggr)
  \to 0 \qquad \text{as \ $t \to \infty$.}
 \]
By \eqref{sup} and \ $\|\bQ(t)\| = 1$, \ $t \in \RR_+$, \ it is enough to show that
 \begin{align}\label{help_maximal_jump_1m}
  t^{-1/2}
  \EE\biggl(\sup_{u\in[0,t]} \|\tbY^{(\ell)}_u - \tbY^{(\ell)}_{u-}\|\biggr)
  \to 0 \qquad \text{as \ $t \to \infty$}
 \end{align}
 for every \ $\ell \in \{1, \ldots, d\}$, \ and
 \begin{align}\label{help_maximal_jump_2m}
  t^{-1/2} \EE\biggl(\sup_{u\in[0,t]} \|\bM^{(5)}_u - \bM^{(5)}_{u-}\|\biggr)
  \to 0 \qquad \text{as \ $t \to \infty$.}
 \end{align}
First, we prove \eqref{help_maximal_jump_1m} for each \ $\ell \in \{1, \ldots, d\}$.
\ By \eqref{CS1}, it is enough to prove that for all \ $\vare \in \RR_{++}$, \
 \[
   t^{-2}
   \EE\biggl(\sup_{u\in[0,t]} |\tbY^{(\ell,\vare)}_u - \tbY^{(\ell,\varepsilon)}_{u-}|^4\biggr)
   \to 0 \qquad \text{as \ $t \to \infty$}
 \]
 and
 \[
   \limsup_{\vare\downarrow0} \limsup_{t\to\infty}
    t^{-1}
    \EE\biggl(\sup_{u\in[0,t]}
    \|(\tbY^{(\ell)}_u
       - \tbY^{(\ell,\vare)}_{u})
      - (\tbY^{(\ell)}_{u-}
         - \tbY^{(\ell,\vare)}_{u-} )\|^2\biggr)
   = 0 .
 \]
By \eqref{CS1+}, for all \ $\vare \in \RR_{++}$, \ we get
 \begin{align*}
  t^{-2}
  \EE\biggl(\sup_{u\in[0,t]}
   |\tbY^{(\ell,\vare)}_u - \tbY^{(\ell,\vare)}_{u-}|^4\biggr)
  &= t^{-2}
     \EE\biggl(\sup_{u\in[0,t]}
      |\tY^{(\ell,\vare)}_u - \tY^{(\ell,\vare)}_{u-}|^4\biggr) \\
  &\leq C_4 \|\bv\|^4 t^{-2}
        \int_0^t \ee^{-s(\tbB)u} \, \dd u
        \int_{\cU_d}
         \|\bz\|^4 \bbone_{\{\|\bz\|\geq\vare\}}
         \, \mu_\ell(\dd\bz)
   \to 0
 \end{align*}
 as \ $t \to \infty$.
\ Further, by \eqref{help17}, for all \ $t \in \RR_{++}$,
 \begin{align*}
  t^{-1}
  \EE\biggl(\sup_{u\in[0,t]}
    \|(\tbY^{(\ell)}_u - \tbY^{(\ell,\vare)}_{u})
      - (\tbY^{(\ell)}_{u-} - \tbY^{(\ell,\vare)}_{u-} )\|^2\biggr)
  \leq C_4 \|\bv\|^2
       \int_{\cU_d}
        \|\bz\|^2 \bbone_{\{\|\bz\|<\vare\}}
        \, \mu_\ell(\dd\bz)
  \to 0
 \end{align*}
 as \ $\vare \downarrow 0$ \ due to \eqref{moment_condition_CBI2}.
Hence we conclude \eqref{help_maximal_jump_1m} for each \ $\ell \in \{1, \ldots, d\}$.

Next, we prove \eqref{help_maximal_jump_2m}.
By \eqref{CS2}, it is enough to prove that
 \[
   t^{-1} \EE\biggl(\sup_{u\in[0,t]} |Z^{(5)}_u - Z^{(5)}_{u-}|^2\biggr)
   \to 0 \qquad \text{as \ $t \to \infty$.}
 \]
By \eqref{Z5}, we get
 \[
   t^{-1} \EE\biggl(\sup_{u\in[0,t]} |Z^{(5)}_u - Z^{(5)}_{u-}|^2\biggr)
   \leq t^{-1} \|\bv\|^2 \int_0^t \ee^{-s(\tbB)u} \, \dd u
        \int_{\cU_d} \|\br\|^2 \, \nu(\dd\br)
   \to 0
 \]
 as \ $t \to \infty$, \ hence we conclude \eqref{help_maximal_jump_2m}.
Consequently, by Theorem \ref{THM_Cri_Pra}, we obtain
 \[
   t^{-1/2} \bQ(t) \bM_t \distr (w_{\bu,\bX_0} \bSigma_\bv)^{1/2} \bN \qquad
   \text{as \ $t \to \infty$,}
 \]
 where \ $\bN$ \ is a \ $2$-dimensional random vector with
 \ $\bN \distre \cN_2(\bzero, \bI_2)$ \ independent of
 \ $w_{\bu,\bX_0} \bSigma_\bv$.
\ Clearly,
 \ $(w_{\bu,\bX_0} \bSigma_\bv)^{1/2} \bN
    = \sqrt{w_{\bu,\bX_0}} \, \bSigma_\bv^{1/2} \bN
    \distre \sqrt{w_{\bu,\bX_0}} \bZ_\bv$.
\ By the decomposition
 \[
   t^{-1/2} \ee^{-s(\tbB)t/2}
   \begin{pmatrix}
    \Re(\langle \bv, \bX_t \rangle) \\
    \Im(\langle \bv, \bX_t \rangle)
   \end{pmatrix}
   = t^{-1/2}
     \begin{pmatrix}
      \Re(\ee^{\ii\Im(\lambda)t} Z_t^{(0,1)}) \\
      \Im(\ee^{\ii\Im(\lambda)t} Z_t^{(0,1)})
     \end{pmatrix}
     + t^{-1/2} \bQ(t) \bM_t , \qquad t \in \RR_+ ,
 \]
 the convergence \eqref{Z01m} and Slutsky's lemma, we obtain \eqref{convvweak2}.
\proofend

\noindent
\textbf{Proof of Theorem \ref{convCBIasL1+}.}
First, suppose that the conditions (i) and (ii) hold.
In the special case of \ $\bX_0 \ase \bzero$, \ applying Lemma
 \ref{decomposition_CBI} with \ $T = 1$, \ we have
 \ $\bX_{t+1} \distre \bX_t^{(1)} + \bX_t^{(2,1)}$ \ for each
 \ $t \in \RR_+$, \ where
 \ $(\bX_s^{(1)})_{s\in\RR_+}$ \ and
 \ $(\bX_s^{(2,1)})_{s\in\RR_+}$ \ are independent multi-type
 CBI processes with \ $\bX_0^{(1)} \ase \bzero$,
 \ $\bX_0^{(2,1)} \distre \bX_1$, \ and with parameters
 \ $(d, \bc, \Bbeta, \bB, \nu, \bmu)$ \ and \ $(d, \bc, \bzero, \bB, 0, \bmu)$,
 \ respectively.
Without loss of generality, we may and do suppose that \ $(\bX_s)_{s\in\RR_+}$,
 \ $(\bX_s^{(1)})_{s\in\RR_+}$ \ and \ $(\bX_s^{(2,1)})_{s\in\RR_+}$ \ are
 independent.
Then, for each \ $t \in \RR_+$, \ we have
 \ $\ee^{-\lambda(t+1)} \langle\bv, \bX_{t+1}\rangle
    \distre \ee^{-\lambda} (\ee^{-\lambda t} \langle\bv, \bX_t^{(1)}\rangle)
            + \ee^{-\lambda} (\ee^{-\lambda t} \langle\bv, \bX_t^{(2,1)}\rangle)$.
\ By \eqref{convwv}, we obtain
 \ $w_{\bv,\bzero}
    \distre \ee^{-\lambda} w_{\bv,\bzero}^{(1)}
            + \ee^{-\lambda} w_{\bv,\bX_0^{(2,1)}}^{(2,1)}$,
 \ where \ $w_{\bv,\bzero}^{(1)}$ \ and \ $w_{\bv,\bX_0^{(2,1)}}^{(2,1)}$ \ denote the almost sure limit of \ $\ee^{-\lambda t} \langle\bv, \bX_t^{(1)}\rangle$
 \ and \ $\ee^{-\lambda t} \langle\bv, \bX_t^{(2,1)}\rangle$ \ as
 \ $t \to \infty$, \ respectively.
Since, for each \ $t \in \RR_+$, \ we have \ $\bX_t^{(1)} \distre \bX_t$, \ we
 conclude \ $w_{\bv,\bzero}^{(1)} \distre w_{\bv,\bzero}$.
\ The independence of \ $(\bX_s)_{s\in\RR_+}$ \ and
 \ $(\bX_s^{(2,1)})_{s\in\RR_+}$ \ implies the independence of \ $w_{\bv,\bzero}$ \ and \ $w_{\bv,\bX_0^{(2,1)}}^{(2,1)}$, \ hence
 \ $w_{\bv,\bzero}
    \distre \ee^{-\lambda} w_{\bv,\bzero}
            + \ee^{-\lambda} w_{\bv,\bX_0^{(2,1)}}^{(2,1)}$.
\ Taking the real and imaginary parts, we get
 \begin{align*}
  \begin{pmatrix} \Re(w_{\bv,\bzero}) \\ \Im(w_{\bv,\bzero}) \end{pmatrix}
  &\distre \begin{pmatrix}
      \Re(\ee^{-\lambda} w_{\bv,\bzero}) \\
      \Im(\ee^{-\lambda} w_{\bv,\bzero})
     \end{pmatrix}
     + \begin{pmatrix}
        \Re(\ee^{-\lambda} w_{\bv,\bX_0^{(2,1)}}^{(2,1)}) \\[2mm]
        \Im(\ee^{-\lambda} w_{\bv,\bX_0^{(2,1)}}^{(2,1)})
       \end{pmatrix} \\
  &= \begin{pmatrix}
      \Re(\ee^{-\lambda}) & - \Im(\ee^{-\lambda}) \\
      \Im(\ee^{-\lambda}) & \Re(\ee^{-\lambda})
     \end{pmatrix}\!
     \begin{pmatrix} \Re(w_{\bv,\bzero}) \\ \Im(w_{\bv,\bzero}) \end{pmatrix}
     + \begin{pmatrix}
        \Re(\ee^{-\lambda}) & - \Im(\ee^{-\lambda}) \\
        \Im(\ee^{-\lambda}) & \Re(\ee^{-\lambda})
       \end{pmatrix}\!
       \begin{pmatrix}
        \Re(w_{\bv,\bX_0^{(2,1)}}^{(2,1)}) \\[2mm]
        \Im(w_{\bv,\bX_0^{(2,1)}}^{(2,1)})
       \end{pmatrix} \\
  &=: \bA \begin{pmatrix} \Re(w_{\bv,\bzero}) \\ \Im(w_{\bv,\bzero}) \end{pmatrix}
      + \bA \bC ,
 \end{align*}
 which is a 2-dimensional stochastic fixed point equation.
We are going to apply Corollary \ref{Cor_sfpe}.
We have
 \ $\det(\bA) = (\Re(\ee^{-\lambda}))^2 + (\Im(\ee^{-\lambda}))^2
    = |\ee^{-\lambda}|^2
    = \ee^{-2\Re(\lambda)} \ne 0$.
\ The eigenvalues of the matrix \ $\bA$ \ are \ $\ee^{-\lambda}$ \ and
 \ $\ee^{-\overline{\lambda}}$, \ hence the spectral radius of \ $\bA$ \ is
 \ $r(\bA) = \ee^{-\Re(\lambda)} \in (0, 1)$.
\ Next we check that \ $\bA \bC$ \ is not deterministic.
Suppose that, on the contrary, \ $\bA \bC$ \ is deterministic.
Then \ $w_{\bv,\bX_0^{(2,1)}}^{(2,1)}$ \ is deterministic, since \ $\bA$ \ is
 invertible.
By Lemma 2.6 in Barczy et al.\ \cite{BarPalPap}, the process
 \ $(\ee^{-s\tbB} \bX_s^{(2,1)})_{s\in\RR_+}$ \ is a
 $d$-dimensional martingale with respect to the filtration
 \ $\cF^{\bX^{(2,1)}}_s := \sigma(\bX_u^{(2,1)} : u \in [0, s])$,
 \ $s \in \RR_+$, \ hence
 \ $(\ee^{-\lambda s} \langle\bv, \bX_s^{(2,1)}\rangle)_{s\in\RR_+}$ \ is a
 complex martingale with respect to the same filtration.
By \eqref{convwv}, we have
 \ $\ee^{-\lambda s} \langle\bv, \bX_s^{(2,1)} \rangle
    \to w_{\bv,\bX_0^{(2,1)}}^{(2,1)}$
 \ as \ $s \to \infty$ \ in \ $L_1$ \ and almost surely, hence
 \ $\langle \bv, \bX_0^{(2,1)}\rangle
    = \EE(w_{\bv,\bX_0^{(2,1)}}^{(2,1)} \mid \cF_0^{\bX^{(2,1)}})
    = w_{\bv,\bX_0^{(2,1)}}^{(2,1)}$
 \ almost surely, see, e.g., Karatzas and Shreve
 \cite[Chapter I, Problem 3.20]{KarShr}.
Thus \ $\langle \bv, \bX_0^{(2,1)}\rangle$ \ is deterministic as well.
Then \ $\langle \bv, \bX_1\rangle$ \ is also deterministic since
 \ $\bX_1 \distre \bX_0^{(2,1)}$.
\ However, applying Lemma \ref{deterministic_CBI} for the process
 \ $(\bX_s)_{s\in\RR_+}$, \ we obtain that \ $\langle \bv, \bX_1\rangle$ \ is not deterministic, since the condition (i) of this theorem implies that the process \ $(\bX_s)_{s\in\RR_+}$ \ is non-trivial, and the condition (ii) of this theorem yields that the condition (ii)/(b) of Lemma \ref{deterministic_CBI} does not hold.
Thus we get a contradiction, and we conclude that \ $\bA \bC$ \ is not
 deterministic.
Moreover, we have \ $\EE(\|\bC\|) = \EE(|w_{\bv,\bX_0^{(2,1)}}^{(2,1)}|) < \infty$,
 \ see \eqref{convwv}.
Applying Corollary \ref{Cor_sfpe}, we conclude that the distribution of
 \ $w_{\bv,\bzero}$ \ does not have atoms.
In particular, we obtain \ $\PP(w_{\bv,\bzero} = 0) = 0$.

If the conditions (i) and (ii) hold, but \ $\bX_0 \ase \bzero$ \ does not necessarily
 holds, then we apply Lemma \ref{decomposition_CBI} with \ $T = 0$, \ and
  we obtain that \ $\bX_t \distre \bX_t^{(1)} + \bX_t^{(2,0)}$ \ for each \ $t \in \RR_+$, \ where \ $(\bX_s^{(1)})_{s\in\RR_+}$ \ and
 \ $(\bX_s^{(2,0)})_{s\in\RR_+}$ \ are
 independent multi-type CBI processes with \ $\bX_0^{(1)} \ase \bzero$,
 \ $\bX_0^{(2,0)} \distre \bX_0$, \ and with parameters
 \ $(d, \bc, \Bbeta, \bB, \nu, \bmu)$ \ and \ $(d, \bc, \bzero, \bB, 0, \bmu)$,
 \ respectively.
Then, for each \ $t \in \RR_+$, \ we have
 \ $\ee^{-\lambda t} \langle\bv, \bX_t\rangle
    \distre \ee^{-\lambda t} \langle\bv, \bX_t^{(1)}\rangle
            + \ee^{-\lambda t} \langle\bv, \bX_t^{(2,0)}\rangle$.
\ By \eqref{convwv}, we obtain
 \ $w_{\bv,\bX_0}
    \distre w_{\bv,\bzero}^{(1)} + w_{\bv,\bX_0^{(2,0)}}^{(2,0)}$,
 \ where \ $w_{\bv,\bzero}^{(1)}$ \ and \ $w_{\bv,\bX_0^{(2,0)}}^{(2,0)}$ \ denotes the
 almost sure limit of \ $\ee^{-\lambda s} \langle\bv, \bX_s^{(1)}\rangle$ \ and of
 \ $\ee^{-\lambda s} \langle\bv, \bX_s^{(2,0)}\rangle$ \ as \ $s \to \infty$,
 \ respectively.
The independence of \ $(\bX_s^{(1)})_{s\in\RR_+}$ \ and
 \ $(\bX_s^{(2,0)})_{s\in\RR_+}$ \ implies the independence of
 \ $w_{\bv,\bzero}^{(1)}$ \ and \ $w_{\bv,\bX_0^{(2,0)}}^{(2,0)}$.
\ We have already shown that \ $w_{\bv,\bzero}^{(1)} \distre w_{\bv,\bzero}$ \ does
 not have atoms, yielding that \ $w_{\bv,\bX_0}$ \ does not have atoms, since for
 each \ $z \in \CC$, \ we have
 \[
   \PP(w_{\bv,\bX_0} = z)
   = \PP\Bigl(w_{\bv,\bzero}^{(1)} = z - w_{\bv,\bX_0^{(2,0)}}^{(2,0)}\Bigr)
   = \EE\Bigl(\PP\Bigl(w_{\bv,\bzero}^{(1)} = z - w_{\bv,\bX_0^{(2,0)}}^{(2,0)}
                       \,\Big|\, w_{\bv,\bX_0^{(2,0)}}^{(2,0)}\Bigr)\Bigr)
   = \EE(0)
   = 0 .
 \]
In particular, we obtain \ $\PP(w_{\bv,\bX_0} = 0) = 0$.

If the condition \textup{(ii)} does not hold, then, as in part
 (ii) $\Longrightarrow$ (iii) of the proof of Lemma \ref{deterministic_CBI}, we
 obtain that in the representation \eqref{reprZ} of
 \ $\ee^{-\lambda t} \langle\bv, \bX_t\rangle$, \ the terms \ $Z_t^{(2)}$, \ $Z_t^{(3,4)}$, \ and \ $Z_t^{(5)}$
 \ are $0$ almost surely, so
 \ $\ee^{-\lambda t} \langle\bv, \bX_t\rangle
    = \langle\bv, \bX_0\rangle
      + \langle\bv, \tBbeta\rangle \int_0^t \ee^{- \lambda u} \, \dd u$
 \ for all \ $t \in \RR_+$ \ almost surely, and hence, taking the limit \ $t\to\infty$, \ we have
 \ $w_{\bv,\bX_0} = \langle\bv, \bX_0\rangle + \lambda^{-1} \langle\bv, \tBbeta\rangle$
 \ almost surely.

If \ $\lambda = s(\tbB)$, \ $\bv = \bu$ \ and the conditions \textup{(i)} and
 \textup{(ii)} hold, then we have already derived \ $\PP(w_{\bu,\bX_0} = 0) = 0$.

If \ $\lambda = s(\tbB)$, \ $\bv = \bu$ \ and the condition \textup{(i)} holds but
 the condition \textup{(ii)} does not hold, then we have already derived
 \ $\PP(w_{\bu,\bX_0} = 0)
    = \PP(\langle\bu, \bX_0\rangle + s(\tbB)^{-1} \langle\bu, \tBbeta\rangle = 0)$,
 \ and this probability is 0, since \ $\bu \in \RR^d_{++}$,
 \ $\PP(\bX_0 \in \RR^d_+) = 1$, \ $s(\tbB) \in \RR_{++}$ \ and
 \ $\tBbeta \in \RR^d_+\setminus \{\bzero\}$ \ yielding that \ $\langle\bu, \tBbeta\rangle > 0$.

If \ $\lambda = s(\tbB)$, \ $\bv = \bu$ \ and the conditions \textup{(i)} and
 \textup{(ii)} do not hold, then we have already derived
 \ $\PP(w_{\bu,\bX_0} = 0)
    = \PP(\langle\bu, \bX_0\rangle + s(\tbB)^{-1} \langle\bu, \tBbeta\rangle = 0)
    = \PP(\langle\bu, \bX_0\rangle = 0)$,
 \ and this equals \ $\PP(\bX_0 = \bzero)$, \ since \ $\bu \in \RR^d_{++}$ \ and
 \ $\PP(\bX_0 \in \RR^d_+) = 1$.
\proofend

\noindent
\textbf{ Proof of Lemma \ref{trivial}.}
Note that \ $w_{\bu,\bX_0} \ase 0$ \ if and only if \ $\EE(w_{\bu,\bX_0}) = 0$.
\ By \eqref{convwv}, we have \ $\ee^{-s(\tbB)t} \langle\bu, \bX_t\rangle \mean w_{\bu,\bX_0}$ \ as \ $t \to \infty$.
\ By \eqref{EXcond}, we obtain \ $\EE(\bX_t) = \ee^{t\tbB} \EE(\bX_0) + \int_0^t \ee^{u\tbB} \tBbeta \, \dd u$, \ $t \in \RR_+$,
 \ hence
 \[
   \EE(\ee^{-s(\tbB)t} \langle\bu, \bX_t\rangle)
    = \langle\bu, \EE(\bX_0) \rangle + \langle\bu,
      \tBbeta\rangle \int_0^t \ee^{-s(\tbB)(t-u)} \, \dd u
    \to \langle\bu, \EE(\bX_0) \rangle
        +  \frac{\langle\bu, \tBbeta\rangle}{s(\tbB)}
    = \EE(w_{\bu,\bX_0})
 \]
 as \ $t \to \infty$, \ where \ $\bu \in \RR_{++}^d$ \ and \ $s(\tbB) > 0$,
 \ thus \ $\EE(w_{\bu,\bX_0}) = 0$ \ if and only if \ $\bX_0 \ase \bzero$ \ and \ $\tBbeta = \bzero$.
\proofend

\noindent
\textbf{Proof of part  (i) of Theorem \ref{convCBIweakr}.}
By Theorem 3.3 in Barczy et al.\ \cite{BarPalPap}, we have
 \ $\ee^{-s(\tbB)t} \bX_t \as w_{\bu,\bX_0} \tbu$ \ as \ $t \to \infty$, \ hence
 \ $\bbone_{\{\bX_t\ne\bzero\}} = \bbone_{\{\ee^{-s(\tbB)t}\bX_t\ne\bzero\}} \to 1$
 \ as \ $t \to \infty$ \ on the event \ $\{w_{\bu,\bX_0} > 0\}$, \ since
 \ $\tbu \in \RR_{++}^d$.
\ By \eqref{convwv}, we have
 \ $\ee^{-s(\tbB)t} \langle\bu, \bX_t\rangle \as w_{\bu,\bX_0}$ \ and
 \ $\ee^{-\lambda t} \langle \bv, \bX_t\rangle \as w_{\bv,\bX_0}$ \ as
 \ $t \to \infty$.
\ Using that \ $\langle\bu, \bX_t\rangle \ne 0$ \ if and only if
 \ $\bX_t \ne \bzero$, \ we have
 \begin{align*}
  &\bbone_{\{\bX_t\ne\bzero\}}
    \frac{1}{\langle \bu, \bX_t \rangle^{\Re(\lambda)/s(\tbB)}}
    \begin{pmatrix}
     \cos(\Im(\lambda)t) & \sin(\Im(\lambda)t) \\
     - \sin(\Im(\lambda)t) & \cos(\Im(\lambda)t)
    \end{pmatrix}
    \begin{pmatrix}
     \Re(\langle \bv, \bX_t \rangle) \\
     \Im(\langle \bv, \bX_t \rangle)
    \end{pmatrix} \\
   &= \frac{\bbone_{\{\bX_t\ne\bzero\}}}
           {(\ee^{-s(\tbB)t}\langle\bu,\bX_t\rangle)^{\Re(\lambda)/s(\tbB)}
            \ee^{\Re(\lambda)t}}
      \begin{pmatrix}
       \cos(\Im(\lambda)t) & \sin(\Im(\lambda)t) \\
       - \sin(\Im(\lambda)t) & \cos(\Im(\lambda)t)
      \end{pmatrix}
      \begin{pmatrix}
       \Re(\ee^{\lambda t} \ee^{-\lambda t} \langle\bv, \bX_t\rangle) \\
       \Im(\ee^{\lambda t} \ee^{-\lambda t} \langle\bv,\bX_t\rangle)
      \end{pmatrix} \\
   &= \frac{\bbone_{\{\bX_t\ne\bzero\}}}
           {(\ee^{-s(\tbB)t}\langle\bu,\bX_t\rangle)^{\Re(\lambda)/s(\tbB)}}
      \begin{pmatrix}
       \Re(\ee^{-\lambda t} \langle\bv, \bX_t\rangle) \\
       \Im(\ee^{-\lambda t} \langle\bv,\bX_t\rangle)
      \end{pmatrix}
    \to \frac{1}{w_{\bu,\bX_0}^{\Re(\lambda)/s(\tbB)}}
         \begin{pmatrix}
          \Re(w_{\bv,\bX_0}) \\
          \Im(w_{\bv,\bX_0})
         \end{pmatrix}
 \end{align*}
 as \ $t \to \infty$ \ on the event \ $\{w_{\bu,\bX_0} > 0\}$, \ as desired.
\proofend

\noindent
\textbf{Proof of part  (iii) of Theorem \ref{convCBIweakr}.}
First, note that the moment condition \eqref{moment_4_2} yields the moment condition \eqref{moment_condition_xlogx}
 with \ $\lambda=s(\tbB)$, \ so, by  Lemma \ref{trivial},
 \ $w_{\bu,\bX_0} \ase 0$ \ if and only if \ $(\bX_t)_{t\in\RR_+}$ \ is trivial.
For each \ $t \in \RR_+$, \ we have the decomposition
 \[
   \bbone_{\{\bX_t\ne\bzero\}} \frac{1}{\sqrt{\langle \bu, \bX_t \rangle}}
   \begin{pmatrix}
    \Re(\langle \bv, \bX_t \rangle) \\
    \Im(\langle \bv, \bX_t \rangle)
   \end{pmatrix}
   = \bbone_{\{\bX_t\ne\bzero\}}
     \frac{\sqrt{w_{\bu,\bX_0}}}
          {\sqrt{\ee^{-s(\tbB)t} \langle \bu, \bX_t \rangle}}
     \frac{\ee^{-s(\tbB)t/2}}{\sqrt{w_{\bu,\bX_0}}}
     \begin{pmatrix}
      \Re(\langle \bv, \bX_t \rangle) \\
      \Im(\langle \bv, \bX_t \rangle)
     \end{pmatrix}
 \]
 on the event \ $\{w_{\bu,\bX_0} > 0\}$.
\ As we have seen in the proof of part  (i) of Theorem \ref{convCBIweakr}, we have
 \ $\bbone_{\{\bX_t\ne\bzero\}} \to 1$ \ as \ $t \to \infty$ \ on the event \ $\{w_{\bu,\bX_0} > 0\}$, \
 and \ $\ee^{-s(\tbB)t} \langle\bu, \bX_t\rangle \as w_{\bu,\bX_0}$ \ as \ $t\to\infty$.
\ In case of \ $\bSigma_\bv = \bzero$, \ \eqref{convvweak1} yields
 \[
   \ee^{-s(\tbB)t/2}
   \begin{pmatrix}
    \Re(\langle \bv, \bX_t \rangle) \\
    \Im(\langle \bv, \bX_t \rangle)
   \end{pmatrix}
   \stoch \bzero
   \qquad \text{as \ $t \to \infty$,}
 \]
 hence, using the above decomposition, by Slutsky's lemma, we obtain
 \[
   \bbone_{\{\bX_t\ne\bzero\}} \frac{1}{\sqrt{\langle \bu, \bX_t \rangle}}
   \begin{pmatrix}
    \Re(\langle \bv, \bX_t \rangle) \\
    \Im(\langle \bv, \bX_t \rangle)
   \end{pmatrix}
   \distrw \bzero
   \qquad \text{as \ $t \to \infty$.}
 \]
 In case of \ $\bSigma_\bv \ne \bzero$, \ as in the proof of \eqref{convvweak2}, we
 may apply Theorem \ref{THM_Cri_Pra} to obtain
 \[
   \biggl(\ee^{-s(\tbB)t/2}
          \begin{pmatrix}
           \Re(\langle \bv, \bX_t \rangle) \\
           \Im(\langle \bv, \bX_t \rangle)
          \end{pmatrix}, \sqrt{w_{\bu,\bX_0}}\biggr)
   \distr \big((w_{\bu,\bX_0} \bSigma_\bv)^{1/2} \bN, \sqrt{w_{\bu,\bX_0}}\big)
   \qquad \text{as \ $t \to \infty$,}
 \]
 where \ $\bN$ \ is a \ $2$-dimensional random vector with
 \ $\bN \distre \cN_2(\bzero, \bI_2)$ \ independent of
 \ $w_{\bu,\bX_0} \bSigma_\bv$, \ and hence independent of \ $w_{\bu,\bX_0}$
 \ because \ $\bSigma_\bv \ne \bzero$ \ and \ $\bSigma_\bv$ \ is deterministic.
Applying the continuous mapping theorem, we get
 \[
   \frac{\ee^{-s(\tbB)t/2}}{\sqrt{w_{\bu,\bX_0}}}
   \begin{pmatrix}
    \Re(\langle \bv, \bX_t \rangle) \\
    \Im(\langle \bv, \bX_t \rangle)
   \end{pmatrix}
   \distrw \bSigma_\bv^{1/2} \bN \qquad \text{as \ $t \to \infty$.}
 \]
Hence, using again the above decomposition, by Slutsky's lemma and \eqref{convwv},
 \begin{align*}
  \bbone_{\{\bX_t\ne\bzero\}} \frac{1}{\sqrt{\langle \bu, \bX_t \rangle}}
  \begin{pmatrix}
   \Re(\langle \bv, \bX_t \rangle) \\
   \Im(\langle \bv, \bX_t \rangle)
  \end{pmatrix}
  &\distrw \bSigma_\bv^{1/2} \bN  \qquad \text{as \ $t \to \infty$,}
 \end{align*}
 where \ $\bSigma_\bv^{1/2} \bN \distre \cN_2(\bzero, \bSigma_\bv)$, \ as desired.
\proofend

\noindent
\textbf{Proof of part  (ii) of Theorem \ref{convCBIweakr}.}
First, note that the moment condition \eqref{moment_4_2} yields the moment condition \eqref{moment_condition_xlogx}
 with \ $\lambda=s(\tbB)$, \ so, by Lemma \ref{trivial},
 \ $w_{\bu,\bX_0} \ase 0$ \ if and only if \ $(\bX_t)_{t\in\RR_+}$ \ is trivial.
For each \ $t \in \RR_+$, \ we have the
 decomposition
 \begin{align*}
  &\bbone_{\{\langle\bu,\bX_t\rangle>1\}}
   \frac{1}{\sqrt{\langle\bu,\bX_t\rangle\log(\langle\bu,\bX_t\rangle)}}
   \begin{pmatrix}
    \Re(\langle \bv, \bX_t \rangle) \\
    \Im(\langle \bv, \bX_t \rangle)
   \end{pmatrix} \\
  &= \bbone_{\{\langle\bu,\bX_t\rangle>1\}}
     \frac{\sqrt{w_{\bu,\bX_0}}}
          {\sqrt{\ee^{-s(\tbB)t}\langle\bu,\bX_t \rangle
                 t^{-1}\log(\langle \bu, \bX_t \rangle)}}
     \frac{t^{-1/2}\ee^{-s(\tbB)t/2}}{\sqrt{w_{\bu,\bX_0}}}
     \begin{pmatrix}
      \Re(\langle \bv, \bX_t \rangle) \\
      \Im(\langle \bv, \bX_t \rangle)
     \end{pmatrix}
 \end{align*}
 on the event \ $\{w_{\bu,\bX_0} > 0\}$.
\ By Theorem 3.3 in Barczy et al.\ \cite{BarPalPap}, we have
 \ $\ee^{-s(\tbB)t} \bX_t \as w_{\bu,\bX_0} \tbu$ \ as \ $t \to \infty$, \ hence
 \ $\bbone_{\{\langle\bu,\bX_t\rangle>1\}}
    = \bbone_{\{\ee^{-s(\tbB)t}\langle\bu,\bX_t\rangle-\ee^{-s(\tbB)t}>0\}} \to 1$
 \ as \ $t \to \infty$ \ on the  event \ $\{w_{\bu,\bX_0} > 0\}$, \ since
 \ $\ee^{-s(\tbB)t} \langle\bu, \bX_t\rangle - \ee^{-s(\tbB)t} \as w_{\bu,\bX_0}$.
\ In case of \ $\bSigma_\bv = \bzero$, \ \eqref{convvweak2} yields
 \[
   t^{-1/2} \ee^{-s(\tbB)t/2}
   \begin{pmatrix}
    \Re(\langle \bv, \bX_t \rangle) \\
    \Im(\langle \bv, \bX_t \rangle)
   \end{pmatrix}
   \stoch \bzero
   \qquad \text{as \ $t \to \infty$,}
 \]
 hence, using the above decomposition, by Slutsky's lemma, we obtain
 \[
   \bbone_{\{\langle\bu,\bX_t\rangle>1\}}
   \frac{1}{\sqrt{\langle\bu,\bX_t\rangle\log(\langle\bu,\bX_t\rangle)}}
   \begin{pmatrix}
    \Re(\langle \bv, \bX_t \rangle) \\
    \Im(\langle \bv, \bX_t \rangle)
   \end{pmatrix}
   \distrw \bzero
   \qquad \text{as \ $t \to \infty$,}
 \]
  since, by \eqref{convwv}, we have
 \ $\ee^{-s(\tbB)t} \langle \bu, \bX_t \rangle \as w_{\bu,\bX_0}$ \ as
 \ $t \to \infty$, \ which also implies
 \[
   t^{-1} \log(\langle \bu, \bX_t \rangle)
   = t^{-1} \log(\ee^{s(\tbB)t})
     + t^{-1} \log(\ee^{-s(\tbB)t} \langle \bu, \bX_t \rangle)
   \as s(\tbB) \in \RR_{++} \qquad \text{as \ $t \to \infty$.}
 \]
In case of \ $\bSigma_\bv \ne \bzero$, \ as in the proof of \eqref{convvweak2}, we
 may apply Theorem \ref{THM_Cri_Pra} to obtain
 \[
   \biggl(t^{-1/2} \ee^{-s(\tbB)t/2}
          \begin{pmatrix}
           \Re(\langle \bv, \bX_t \rangle) \\
           \Im(\langle \bv, \bX_t \rangle)
          \end{pmatrix}, \sqrt{w_{\bu,\bX_0}}\biggr)
   \distr \big((w_{\bu,\bX_0} \bSigma_\bv)^{1/2} \bN, \sqrt{w_{\bu,\bX_0}}\big)
   \qquad \text{as \ $t \to \infty$,}
 \]
 where \ $\bN$ \ is a \ $2$-dimensional random vector with
 \ $\bN \distre \cN_2(\bzero, \bI_2)$ \ independent of
 \ $w_{\bu,\bX_0} \bSigma_\bv$, \ and hence independent of \ $w_{\bu,\bX_0}$
 \ because \ $\bSigma_\bv \ne \bzero$ \ and \ $\bSigma_\bv$ \ is deterministic.
Applying the continuous mapping theorem, we get
 \[
   \frac{t^{-1/2}\ee^{-s(\tbB)t/2}}{\sqrt{w_{\bu,\bX_0}}}
   \begin{pmatrix}
    \Re(\langle \bv, \bX_t \rangle) \\
    \Im(\langle \bv, \bX_t \rangle)
   \end{pmatrix}
   \distrw \bSigma_\bv^{1/2} \bN \qquad \text{as \ $t \to \infty$.}
 \]
Hence, using again the above decomposition, by Slutsky's lemma and \eqref{convwv},
 \begin{align*}
  \bbone_{\{\langle\bu,\bX_t\rangle>1\}}
  \frac{1}{\sqrt{\langle\bu,\bX_t\rangle\log(\langle\bu,\bX_t\rangle)}}
  \begin{pmatrix}
   \Re(\langle \bv, \bX_t \rangle) \\
   \Im(\langle \bv, \bX_t \rangle)
  \end{pmatrix}
  &\distrw \frac{1}{s(\tbB)^{1/2}} \bSigma_\bv^{1/2} \bN  \qquad
   \text{as \ $t \to \infty$,}
 \end{align*}
 where \ $\frac{1}{s(\tbB)^{1/2}} \bSigma_\bv^{1/2} \bN \distre \cN_2(\bzero, \frac{1}{s(\tbB)} \bSigma_\bv)$, \ as desired.
\proofend

\noindent{\bf Proof of Proposition \ref{Cor_relative_frequency}.}
Theorem 3.3 in Barczy et al.\ \cite{BarPalPap} yields that
 \ $\ee^{-s(\tbB)t} \langle\be_i, \bX_t\rangle
    \as w_{\bu,\bX_0} \langle\be_i, \tbu\rangle$
 \ and
 \ $\ee^{-s(\tbB)t} \langle\be_j, \bX_t\rangle
    \as w_{\bu,\bX_0} \langle\be_j, \tbu\rangle$
 \ as \ $t \to \infty$.
\ Consequently, since \ $\tbu \in \RR_{++}$, \ we have
 \ $\bbone_{\{\langle\be_j,\bX_t\rangle\ne0\}}
    = \bbone_{\{\ee^{-s(\tbB)t}\langle\be_j,\bX_t\rangle\ne 0\}} \to 1$
 \ as \ $t \to \infty$ \ on the event \ $\{w_{\bu,\bX_0} > 0\}$, \ and hence
 \begin{align*}
  \bbone_{\{\langle\be_j,\bX_t\rangle\ne0\}}
  \frac{\langle\be_i,\bX_t\rangle}{\langle\be_j,\bX_t\rangle}
  = \bbone_{\{\langle\be_j,\bX_t\rangle\ne0\}}
    \frac{\ee^{-s(\tbB)t}\langle\be_i,\bX_t\rangle}
         {\ee^{-s(\tbB)t}\langle\be_j,\bX_t\rangle}
  \to \frac{w_{\bu,\bX_0}\langle\be_i,\tbu\rangle}
           {w_{\bu,\bX_0}\langle\be_j,\tbu\rangle}
  = \frac{\langle\be_i,\tbu\rangle}{\langle\be_j,\tbu\rangle}
  \qquad \text{as \ $t \to \infty$}
 \end{align*}
 on the event \ $\{w_{\bu,\bX_0} > 0\}$, \ thus we obtain the first convergence.
In a similar way, \ $\bbone_{\{\bX_t\ne\bzero\}} \to 1$ \ as \ $t \to \infty$ \ on the event \ $\{w_{\bu,\bX_0} > 0\}$, \ thus
 \begin{align*}
  \bbone_{\{\bX_t\ne\bzero\}}
  \frac{\langle\be_i,\bX_t\rangle}{\sum_{k=1}^d\langle\be_k,\bX_t\rangle}
  = \bbone_{\{\bX_t\ne\bzero\}}
    \frac{\ee^{-s(\tbB)t}\langle\be_i,\bX_t\rangle}
         {\sum_{k=1}^d\ee^{-s(\tbB)t}\langle\be_k,\bX_t\rangle}
  \to \frac{w_{\bu,\bX_0}\langle\be_i,\tbu\rangle}
           {\sum_{k=1}^dw_{\bu,\bX_0}\langle\be_k,\tbu\rangle}
  = \langle\be_i, \tbu\rangle
 \end{align*}
 as \ $t \to \infty$ \ on the event \ $\{w_{\bu,\bX_0} > 0\}$ \ since the sum of the coordinates of \ $\tbu$ \ is \ 1,
 \ hence we obtain the second convergence.
\proofend

\appendix

\vspace*{5mm}

\noindent{\bf\Large Appendix}

\section{A decomposition of multi-type CBI processes}
\label{deco_CBI}

The following useful decomposition of a multi-type CBI process as an independent sum of a
 CBI process starting from $\bzero$ and a CB process has been derived in Barczy et al.\ \cite[Lemma A.1]{BarPalPap}.

\begin{Lem}\label{decomposition_CBI}
If \ $(\bX_s)_{s\in\RR_+}$ \ is a multi-type CBI process with parameters
 \ $(d, \bc, \Bbeta, \bB, \nu, \bmu)$, \ then for each \ $t, T \in \RR_+$, \ we
 have \ $\bX_{t+T} \distre \bX_t^{(1)} + \bX_t^{(2,T)}$, \ where
 \ $(\bX_s^{(1)})_{s\in\RR_+}$ \ and
 \ $(\bX_s^{(2,T)})_{s\in\RR_+}$ \ are independent multi-type CBI processes with \ $\PP(\bX_0^{(1)} = \bzero) = 1$,
 \ $\bX_0^{(2,T)} \distre \bX_T$, \ and with parameters
 \ $(d, \bc, \Bbeta, \bB, \nu, \bmu)$ \ and \ $(d, \bc, \bzero, \bB, 0, \bmu)$,
 \ respectively.
\end{Lem}

\section{On deterministic projections of multi-type CBI processes}
\label{det_CBI}

\begin{Lem}\label{deterministic_CBI}
Let \ $(\bX_t)_{t\in\RR_+}$ \ be an irreducible multi-type CBI process with parameters
 \ $(d, \bc, \Bbeta, \bB, \nu, \bmu)$ \ such that \ $\EE(\|\bX_0\|) < \infty$ \ and
 the moment conditions \eqref{moment_condition_m_new} and \eqref{moment_condition_CBI2} hold.
Let \ $\lambda \in \sigma(\tbB)$, \ and let \ $\bv \in \CC^d$ \ be a left
 eigenvector of \ $\tbB$ \ corresponding to the eigenvalue \ $\lambda$.
\ Then the following three assertions are equivalent:
 \renewcommand{\labelenumi}{{\rm(\roman{enumi})}}
 \begin{enumerate}
  \item
   There exists \ $t \in \RR_{++}$ \ such that \ $\langle\bv, \bX_t\rangle$ \ is
    deterministic.
  \item
   One of the following two conditions holds:
    \renewcommand{\labelenumii}{{\rm(\alph{enumii})}}
    \begin{enumerate}
     \item
      $(\bX_t)_{t\in\RR_+}$ \ is a trivial process (see Definition
       \ref{Def_nontrivial}).
     \item
      $\langle\bv, \bX_0\rangle$ \ is deterministic,
       \ $\langle\bv, \be_\ell\rangle c_\ell = 0$ \ and
       \ $\mu_\ell(\{\bz \in \cU_d : \langle\bv, \bz\rangle \ne 0\}) = 0$ \ for
       every \ $\ell \in \{1, \ldots, d\}$, \ and
       \ $\nu(\{\br \in \cU_d : \langle\bv, \br\rangle \ne 0\}) = 0$.
    \end{enumerate}
  \item
   For each \ $t \in \RR_+$, \ $\langle\bv, \bX_t\rangle$ \ is deterministic.
 \end{enumerate}
If \ $(\langle\bv, \bX_t\rangle)_{t\in\RR_+}$ \ is deterministic, then
 \ $\langle\bv, \bX_t\rangle
    \ase \ee^{\lambda t} \langle\bv, \EE(\bX_0)\rangle
         + \langle\bv, \tBbeta\rangle \int_0^t \ee^{\lambda u} \, \dd u$,
 \ $t \in \RR_+$.
\end{Lem}

\noindent
\textbf{Proof.}
\ (i) $\Longrightarrow$ (ii).
We have the representation
 \begin{equation}\label{reprZ}
  \ee^{-\lambda t} \langle\bv, \bX_t\rangle
  = Z_t^{(0)} + Z_t^{(1)} + Z_t^{(2)} + Z_t^{(3,4)} + Z_t^{(5)}
 \end{equation}
 with
 \begin{align*}
  Z_t^{(0)}
  &:= \langle\bv, \bX_0\rangle , \\
  Z_t^{(1)}
  &:= \langle\bv, \tBbeta\rangle \int_0^t \ee^{-\lambda u} \, \dd u , \\
  Z_t^{(2)}
  &:= \sum_{\ell=1}^d
       \langle\bv, \be_\ell\rangle
       \int_0^t \ee^{-\lambda u} \sqrt{2 c_\ell X_{u,\ell}} \, \dd W_{u,\ell} , \\
  Z_t^{(3,4)}
  &:= \sum_{\ell=1}^d
       \int_0^t \int_{\cU_d} \int_{\cU_1}
        \ee^{-\lambda u} \langle\bv, \bz\rangle
        \bbone_{\{w\leq X_{u-,\ell}\}}
        \, \tN_\ell(\dd u, \dd\bz, \dd w) , \\
  Z_t^{(5)}
  &:= \int_0^t \int_{\cU_d}
       \ee^{-\lambda u} \langle\bv, \br\rangle \, \tM(\dd u, \dd\br) ,
 \end{align*}
 see Barczy et al.\ \cite[Lemma 4.1]{BarLiPap3} or Barczy et al.\ \cite[Lemma 2.7]{BarPalPap}.
Note that under the moment condition \eqref{moment_condition_CBI2}, \ $(Z_t^{(2)})_{t\in\RR_+}$, \ $(Z_t^{(3,4)})_{t\in\RR_+}$ \
 and \ $(Z_t^{(5)})_{t\in\RR_+}$ \ are square-integrable martingales with initial values \ $0$, \ hence
 \ $\EE(Z_t^{(2)}) = \EE(Z_t^{(3,4)}) = \EE(Z_t^{(5)}) = 0$.
\ Since \ $\ee^{-\lambda t} \langle\bv, \bX_t\rangle$ \ and \ $Z_t^{(1)}$ \ are
 deterministic, we obtain
 \ $\ee^{-\lambda t} \langle\bv, \bX_t\rangle
    = \EE(\ee^{-\lambda t} \langle\bv, \bX_t\rangle)
    = \EE(\langle\bv, \bX_0\rangle) + Z_t^{(1)}$.
\ Hence, by the representation \eqref{reprZ}, we get
 \ $0 = \ee^{-\lambda t} \langle\bv, \bX_t\rangle
        - \EE(\ee^{-\lambda t} \langle\bv, \bX_t\rangle)
      = \langle\bv, \bX_0\rangle - \EE(\langle\bv, \bX_0\rangle)
        + \sum_{j=2}^5 Z_t^{(j)}$
 \ almost surely.
Consequently,
 \[
   \EE(|\langle\bv, \bX_0\rangle - \EE(\langle\bv, \bX_0\rangle)
         + Z_t^{(2)} + Z_t^{(3,4)} + Z_t^{(5)} |^2)
   = 0.
 \]
By the independence of \ $\bX_0$,
 \ $(W_{u,1})_{u\geq0}$, \ldots, \ $(W_{u,d})_{u\geq0}$, \ $N_1$, \ldots, \ $N_d$
 \ and \ $M$, \ the random variables
 \ $\langle\bv, \bX_0\rangle - \EE(\langle\bv, \bX_0\rangle)$, \ $Z_t^{(2)}$, $Z_t^{(3,4)}$, \ and
 \ $Z_t^{(5)}$ \ are conditionally independent with respect to \ $(\bX_u)_{u\in[0,t]}$, \ thus
 \begin{align*}
  0 &= \EE\biggl(\biggl|\langle\bv, \bX_0\rangle - \EE(\langle\bv, \bX_0\rangle)
                        + Z_t^{(2)} + Z_t^{(3,4)} + Z_t^{(5)} \biggr|^2\biggr) \\
    &= \EE\biggl(\EE\biggl(\biggl|\langle\bv, \bX_0\rangle
                                  - \EE(\langle\bv, \bX_0\rangle)
                                  + Z_t^{(2)} + Z_t^{(3,4)} + Z_t^{(5)} \biggr|^2
                                  \,\bigg|\, (\bX_u)_{u\in[0,t]}\biggr)\biggr) \\
    &= \EE(\EE(|\langle\bv, \bX_0\rangle - \EE(\langle\bv, \bX_0\rangle)|^2
               \mid (\bX_u)_{u\in[0,t]}))
       + \EE(\EE(|Z_t^{(2)}|^2 \mid (\bX_u)_{u\in[0,t]})) \\
    &\phantom{=\;}
       + \EE(\EE(|Z_t^{(3,4)}|^2 \mid (\bX_u)_{u\in[0,t]}))
       + \EE(\EE(|Z_t^{(5)}|^2 \mid (\bX_u)_{u\in[0,t]})) \\
    &= \EE(|\langle\bv, \bX_0\rangle - \EE(\langle\bv, \bX_0\rangle)|^2)
       + \EE(|Z_t^{(2)}|^2) + \EE(|Z_t^{(3,4)}|^2) + \EE(|Z_t^{(5)}|^2),
 \end{align*}
  where we also used that \ $(Z_s^{(2)})_{s\in[0,t]}$, \ $(Z_s^{(3,4)})_{s\in[0,t]}$ \ and \ $(Z_s^{(5)})_{s\in[0,t]}$ \ are square-integrable martingales
  with initial values \ $0$ \ conditionally on \ $(\bX_u)_{u\in[0,t]}$.
\ Consequently,
 \ $\EE(|\langle\bv, \bX_0\rangle - \EE(\langle\bv, \bX_0\rangle)|^2 ) = 0$ \ and
 \ $\EE(|Z_t^{(2)}|^2) = \EE(|Z_t^{(3,4)}|^2) = \EE(|Z_t^{(5)}|^2) = 0$.
\ One can easily derive
 \[
   \EE(|Z_t^{(2)}|^2)
   = 2 \sum_{\ell=1}^d
        |\langle\bv, \be_\ell\rangle|^2 c_\ell
        \int_0^t \ee^{-2\Re(\lambda)u} \EE(X_{u,\ell}) \, \dd u ,
 \]
 hence we conclude
 \[
   |\langle\bv, \be_\ell\rangle|^2 c_\ell
   \int_0^t \ee^{-2\Re(\lambda)u} \EE(X_{u,\ell}) \, \dd u = 0 , \qquad
   \ell \in \{1, \ldots, d\} .
 \]
Consequently, for each \ $\ell \in \{1, \ldots, d\}$, \ we have
 \ $|\langle\bv, \be_\ell\rangle|^2 c_\ell = 0$ \ or
 \ $\int_0^t \ee^{-2\Re(\lambda)u} \EE(X_{u,\ell}) \, \dd u = 0$.
\ In the first case we obtain \ $\langle\bv, \be_\ell\rangle c_\ell = 0$, \ which
 is in (ii)/(b).
In the second case, using Lemma \ref{expectation_CBI} and
 \ $\ee^{-2\Re(\lambda)u} \in \RR_{++}$ \ for all \ $u \in \RR_+$, \ we conclude
 (ii)/(a).

Since \ $\EE(|Z_t^{(3,4)}|^2) = 0$, \ we have
 \[
   \EE(|Z_t^{(3,4)}|^2)
   = \sum_{\ell=1}^d
      \int_0^t \int_{\cU_d}
       \ee^{-2\Re(\lambda)u} |\langle\bv, \bz\rangle|^2 \EE(X_{u,\ell})
       \, \dd u \, \mu_\ell(\dd\bz)
   = 0 .
 \]
Using \ $\ee^{-2\Re(\lambda)u} \in \RR_{++}$ \ for all \ $u \in \RR_+$, \ we conclude
 \begin{equation*}
  \sum_{\ell=1}^d
   \int_0^t \int_{\cU_d}
    \bbone_{\{\langle\bv, \bz\rangle\ne0\}} \EE(X_{u,\ell})
    \, \dd u \, \mu_\ell(\dd\bz)
  = 0 .
 \end{equation*}
Then, using the non-negativity of the
 integrands, we obtain
 \[
   \int_0^t \int_{\cU_d}
    \bbone_{\{\langle\bv, \bz\rangle\ne0\}} \EE(X_{u,\ell})
    \, \dd u \, \mu_\ell(\dd\bz)
   = 0 , \qquad \ell \in \{1, \ldots, d\} .
 \]
By Lemma \ref{expectation_CBI}, for each \ $\ell \in \{1, \ldots, d\}$, \ we have
 either (ii)/(a), or \ $\EE(X_{u,\ell}) = \be_\ell^\top \EE(\bX_u) \in \RR_{++}$
 \ for all \ $u \in \RR_{++}$.
\ In the second case, we conclude
 \[
   \int_0^t \int_{\cU_d}
    \bbone_{\{\langle\bv, \bz\rangle\ne0\}}
    \, \dd u \, \mu_\ell(\dd\bz)
   = t \mu_\ell(\{\bz \in \cU_d : \langle\bv, \bz\rangle \ne 0\})
   = 0 ,
 \]
 and hence \ $\mu_\ell(\{\bz \in \cU_d : \langle\bv, \bz\rangle \ne 0\}) = 0$,
 \ which is in (ii)/(b).

Since \ $\EE(|Z_t^{(5)}|^2) = 0$, \ we have \ $Z_t^{(5)} = 0$ \ almost surely.
Hence the random variable
 \[
   \int_0^t \int_{\cU_d}
    \ee^{-\lambda u} \langle\bv, \br\rangle \, M(\dd u, \dd\br)
 \]
 is deterministic, since
 \ $\int_0^t \int_{\cU_d}
     \ee^{-\lambda u} \langle\bv, \br\rangle \, \dd u \, \nu(\dd\br)$
 \ is deterministic.
We have \ $Z_s^{(5)} = \EE(Z_t^{(5)} \mid \cF_s^{Z^{(5)}}) =  \EE( 0 \mid \cF_s^{Z^{(5)}}) = 0$ \
 for all \ $s \in [0,t]$ \ almost surely, where \ $\cF_s^{Z^{(5)}} :=\sigma(Z^{(5)}_u : u\in[0,s])$,
 \ since \ $(Z_s^{(5)})_{s\in\RR_+}$ \ is a martingale.
Thus \ $\PP(A_t^{(M)}) = 1$, \ where \ $A_t^{(M)}$ \ is the event such that the
 Poisson random measure \ $M$ \ has no point in the set \ $H_t$, \ where
 \[
   H_t
   := \{(u, \br) \in (0, t] \times \cU_d
        : \ee^{-\lambda u} \langle\bv, \br\rangle \ne 0\}
    = \{(u, \br) \in (0, t] \times \cU_d
        : \bbone_{\{\langle\bv, \br\rangle\ne0\}} \ne 0\} ,
 \]
 since \ $\ee^{-\lambda u} \ne 0$ \ for all \ $u \in \RR_+$.
\ The number of the points of \ $M$ \ in the set \ $H_t$ \ has a Poisson
 distribution with parameter
 \[
   \lambda_t
   := \int_0^t \int_{\cU_d}
       \bbone_{\{\langle\bv, \br\rangle\ne0\}} \, \dd u \, \nu(\dd\br) .
 \]
We have \ $1 = \PP(A_t^{(M)}) = \ee^{-\lambda_t}$, \ yielding
 \[
   \lambda_t
   = \int_0^t \int_{\cU_d}
      \bbone_{\{\langle\bv, \br\rangle\ne0\}} \, \dd u \, \nu(\dd\br)
   = t \nu(\{\br \in \cU_d : \langle\bv, \br\rangle \ne 0\})
   = 0 ,
 \]
 and hence \ $\nu(\{\br \in \cU_d : \langle\bv, \br\rangle \ne 0\}) = 0$, \ which
 is in (ii)/(b).

(ii) $\Longrightarrow$ (iii).
If (ii)/(a) holds, then \ $\langle\bv, \bX_t\rangle \ase 0$ \ for all
 \ $t \in \RR_+$.
\ If (ii)/(b) holds, then we use again the representation \eqref{reprZ} of
 \ $\langle\bv, \bX_t\rangle$.
\ We have
 \ $\langle\bv, \bX_0\rangle = \EE(\langle\bv, \bX_0\rangle)
    = \langle\bv, \EE(\bX_0)\rangle$,
 \ since \ $\langle\bv, \bX_0\rangle$ \ is deterministic.
For each \ $t \in \RR_+$, \ we have
 \[
   Z_t^{(2)}
   = \sqrt{2} \sum_{\ell=1}^d
        \langle\bv, \be_\ell\rangle \sqrt{c_\ell}
        \int_0^t \ee^{-\lambda u} \sqrt{X_{u,\ell}} \, \dd W_{u,\ell}
   = 0 ,
 \]
 since \ $\langle\bv, \be_\ell\rangle c_\ell = 0$ \ for every
 \ $\ell \in \{1, \ldots, d\}$.

Further, for each \ $t \in \RR_+$ \ and \ $n \in \NN$, \ using the notation
 \begin{align*}
   f(u, \bz, w)
   := \sum_{\ell=1}^d
       \ee^{-\lambda u} \langle\bv, \bz\rangle \bbone_{\{w\leq X_{u-,\ell}\}}
   = \ee^{-\lambda u} \langle\bv, \bz\rangle \sum_{\ell=1}^d  \bbone_{\{w\leq X_{u-,\ell}\}}
 \end{align*}
 for \ $u \in (0, t]$, \ $\bz \in \cU_d$, \ and \ $w \in \cU_1$, \ we have
 \begin{align*}
  &\int_0^t \int_{\cU_d} \int_{\cU_1}
    \bbone_{\{|f(u,\bz,w)|<n\}} \bbone_{\{\|\bz\|>1/n\}} \bbone_{\{w<n\}}
    f(u, \bz, w)
    \, \tN_\ell(\dd u, \dd\bz, \dd w) \\
  &= \int_0^t \int_{\cU_d} \int_{\cU_1}
      \bbone_{\{|f(u,\bz,w)|<n\}} \bbone_{\{\|\bz\|>1/n\}} \bbone_{\{w<n\}}
      f(u, \bz, w)
      \, N_\ell(\dd u, \dd\bz, \dd w) \\
  &\quad
     - \int_0^t \int_{\cU_d} \int_{\cU_1}
        \bbone_{\{|f(u,\bz,w)|<n\}} \bbone_{\{\|\bz\|>1/n\}} \bbone_{\{w<n\}}
        f(u, \bz, w)
        \, \dd u \, \mu_\ell(\dd\bz) \, \dd w
   = 0
 \end{align*}
 almost surely, since
 \ $\int_{\cU_d} \bbone_{\{\|\bz\|>1/n\}} \, \mu_\ell(\dd\bz)
        \leq n^2 \int_{\cU_d} \|\bz\|^2 \, \mu_\ell(\dd\bz) < \infty$
 \ due to part (vi) of Definition \ref {Def_admissible}, \eqref{moment_condition_CBI2}
  and
 \[
   (\cL_1\otimes \mu_\ell\otimes \cL_d)(\{ (u,\bz,w) \in (0,t]\times\cU_d\times \cU_1 : f(u, \bz, w) \ne 0\}) = 0
 \]
 for each \ $\ell \in \{1, \ldots, d\}$, \ where \ $\cL_1$ \ and \ $\cL_d$ \ denote the Lebesgue measure on \ $\RR$ \
 and on \ $\RR^d$, \ respectively.
Letting \ $n \to \infty$, \ by Ikeda and Watanabe \cite[page 63]{IkeWat}, we
 conclude
 \[
   Z_t^{(3,4)}
   = \sum_{\ell=1}^d
      \int_0^t \int_{\cU_d} \int_{\cU_1}
       \ee^{-\lambda u} \langle\bv, \bz\rangle
       \bbone_{\{w\leq X_{u-,\ell}\}}
       \, \tN_\ell(\dd u, \dd\bz, \dd w)
   = 0
 \]
 almost surely.

Finally, for each \ $t \in \RR_+$, \ we have
 \begin{align*}
  Z^{(5)}_t
  = \int_0^t \int_{\cU_d}
     \ee^{-\lambda u} \langle \bv,\br\rangle \, M(\dd u, \dd\br)
    - \int_0^t \int_{\cU_d}
       \ee^{-\lambda u} \langle \bv,\br\rangle \, \dd u \, \nu(\dd\br)
  = 0
 \end{align*}
 almost surely, since \ $\int_{\cU_d} \|\br\| \, \nu(\dd\br) < \infty$ \ (due to
 Definition \ref{Def_admissible} and \eqref{moment_condition_m_new}) and
 \ $\nu(\{\br \in \cU_d : \langle\bv, \br\rangle \ne 0\}) = 0$.

(iii) $\Longrightarrow$ (i) is trivial.

If \ $(\langle\bv, \bX_t\rangle)_{t\in\RR_+}$ \ is deterministic, then, by
 \eqref{EXcond}, for each \ $t \in \RR_+$, \ we have
 \ $\langle\bv, \bX_t\rangle
    = \EE(\langle\bv, \bX_t\rangle)
    = \langle\bv, \EE(\bX_t)\rangle
    = \ee^{\lambda t} \langle\bv, \EE(\bX_0)\rangle
      + \langle\bv, \tBbeta\rangle \int_0^t \ee^{\lambda u} \, \dd u$
 \ almost surely.
\proofend

\section{A stochastic fixed point equation}
\label{section_stoch_fixed_point_eq}

Under some mild conditions, the solution of a stochastic fixed point equation
 is atomless, see, e.g., Buraczewski et al. \cite[Proposition 4.3.2]{BurDamMik}.

\begin{Thm}\label{sfpe}
Let \ $(\bA, \bC)$ \ be a random element in \ $\RR^{d \times d} \times \RR^d$,
 \ where \ $d \in \NN$.
\ Assume that
 \renewcommand{\labelenumi}{{\rm(\roman{enumi})}}
 \begin{enumerate}
  \item
   $\bA$ \ is invertible almost surely,
  \item
   $\PP(\bA \bx + \bC = \bx) < 1$ \ for every \ $\bx \in \RR^d$,
  \item
   the \ $d$-dimensional fixed point equation \ $\bX \distre \bA \bX + \bC$,
    \ where \ $(\bA, \bC)$ \ and \ $\bX$ \ are independent, has a solution \ $\bX$,
    \ which is unique in distribution.
 \end{enumerate}
Then the distribution of \ $\bX$ \ does not have atoms and is of pure type, i.e.,
 it is either absolutely continuous or singular with respect to Lebesgue measure in
 \ $\RR^d$.
\end{Thm}

\begin{Cor}\label{Cor_sfpe}
Let \ $\bA \in \RR^{d \times d}$ \ with \ $\det(\bA) \ne 0$ \ and \ $r(\bA) < 1$.
\ Let \ $\bC$ \ be a $d$-dimensional non-deterministic random vector with
 \ $\EE(\|\bC\|) < \infty$.
\ Then the \ $d$-dimensional fixed point equation \ $\bX \distre \bA \bX + \bC$,
 \ where \ $\bX$ \ is independent of \ $\bC$, \ has a solution \ $\bX$ \ which is
 unique in distribution, the distribution of \ $\bX$ \ does not have atoms and is
 of pure type, i.e., it is either absolutely continuous or singular with respect to
 Lebesgue measure in \ $\RR^d$.
\end{Cor}

\noindent
\textbf{Proof.}
The first condition of Theorem \ref{sfpe} is trivially satisfied, since
 \ $\det(\bA) \ne 0$.
\ Since \ $\bC$ \ is not deterministic and for each \ $\bx \in \RR^d$, \ we have
 \ $\PP(\bA \bx + \bC = \bx) = \PP(\bC = (\bI_d - \bA) \bx)$, \ the second
 condition of Theorem \ref{sfpe} is also satisfied.
In order to check the third condition of Lemma \ref{sfpe}, first we suppose that
 \ $\bX$ \ is a solution of the stochastic fixed point equation
 \ $\bX \distre \bA \bX + \bC$, \ where \ $\bX$ \ is a $d$-dimensional random
 vector independent of \ $(\bA, \bC)$, \ equivalently, independent of \ $\bC$
 \ (since \ $\bA$ \ is deterministic and invertible).
Then, iterating this equation, for each \ $n \in \NN$, \ we obtain
 \ $\bX \distre \bA^n \bX + \sum_{k=0}^{n-1} \bA^k \bC_k$, \ where \ $\bC_k$,
 \ $k \in \ZZ_+$, \ are independent copies of \ $\bC$.
\ Since \ $r(\bA) < 1$, \ we have \ $\bA^n \to \bzero$ \ as \ $n \to \infty$,
 \ see, e.g., Horn and Johnson \cite[Theorem 5.6.12]{HorJoh}.
Moreover, \ $\sum_{k=0}^{n-1} \bA^k \bC_k \mean \sum_{k=0}^\infty \bA^k \bC_k$ \ as
 \ $n \to \infty$, \ since \ $\sum_{k=n}^\infty \bA^k \bC_k \mean \bzero$ \ as
 \ $n \to \infty$.
\ Indeed, by the Gelfand formula, we have
 \ $r(\bA) = \lim_{k\to\infty} \|\bA^k\|^{1/k}$, \ see, e.g., Horn and Johnson
 \cite[Corollary 5.6.14]{HorJoh}, hence there exists \ $k_0 \in \NN$ \ such that
 \ $\|\bA^k\|^{1/k} \leq (r(\bA)+1)/2 < 1$ \ for every \ $k \in \NN$ \ with
 \ $k \geq k_0$.
\ Thus, for each \ $n \in \NN$ \ with \ $n \geq k_0$, \ we have
 \[
   \EE\biggl(\bigg\|\sum_{k=n}^\infty \bA^k \bC_k\biggr\|\biggr)
   \leq \sum_{k=n}^\infty \|\bA^k\| \EE(\|\bC_k\|)
   \leq \EE(\|\bC\|) \sum_{k=n}^\infty \Bigl(\frac{r(\bA)+1}{2}\Bigr)^k
   \to 0
 \]
 as \ $n \to \infty$, \ hence we obtain
 \ $\sum_{k=n}^\infty \bA^k \bC_k \mean \bzero$ \ as \ $n \to \infty$, \ and
 hence \ $\sum_{k=0}^{n-1} \bA^k \bC_k \mean \sum_{k=0}^\infty \bA^k \bC_k$ \ as
 \ $n \to \infty$.
\ Consequently, if \ $\bX$ \ is a solution of \ $\bX \distre \bA \bX + \bC$,
 \ then, necessarily, \ $\bX \distre \sum_{k=0}^\infty \bA^k \bC_k$.
\ The $d$-dimensional random variable \ $\sum_{k=0}^\infty \bA^k \bC_k$ \ is a
 solution of \ $\bX \distre \bA \bX + \bC$, \ since
 \ $\sum_{k=0}^\infty \bA^k \bC_k
    = \bA \sum_{k=0}^\infty \bA^k \bC_{k+1} + \bA \bC_0$,
 \ where
 \ $\sum_{k=0}^\infty \bA^k \bC_{k+1} \distre \sum_{k=0}^\infty \bA^k \bC_k$ \ and
 \ $\sum_{k=0}^\infty \bA^k \bC_{k+1}$ \ is independent of \ $A \bC_0$ \ (equivalently,
  of \ $(\bA,\bA \bC_0)$), \ hence the third condition of Lemma \ref{sfpe} is also satisfied.
\proofend

\section{On the second moment of projections of multi-type CBI processes}
\label{second_moment_CBI}

An explicit formula for the second absolute moment of the projection of a multi-type CBI
 process on the left eigenvectors of its branching mean matrix has been presented
 together with its asymptotic behavior in the supercritical and irreducible case in
 Barczy et al.\ \cite[Proposition B.1]{BarPalPap}.

\begin{Pro}\label{second_moment_asymptotics_CBI}
If \ $(\bX_t)_{t\in\RR_+}$ \ is a multi-type CBI process with parameters
 \ $(d, \bc, \Bbeta, \bB, \nu, \bmu)$ \ such that \ $\EE(\|\bX_0\|^2) < \infty$
 \ and the moment condition \eqref{moment_condition_CBI2} holds,
  then for each left eigenvector \ $\bv \in \CC^d$ \ of \ $\tbB$
 \ corresponding to an arbitrary eigenvalue \ $\lambda \in \sigma(\tbB)$, \ we have
 \[
   \EE(|\langle\bv, \bX_t\rangle|^2)
   = E_{\bv,\lambda}(t) + \sum_{\ell=1}^d C_{\bv,\ell} I_{\lambda,\ell}(t)
     + I_\lambda(t) \int_{\cU_d} |\langle\bv, \br\rangle|^2 \, \nu(\dd\br) ,
   \qquad t \in \RR_+ ,
 \]
 where \ $C_{\bv,\ell}$, \ $\ell \in \{1, \ldots, d\}$, \ are defined in Theorem
 \ref{convCBIweak1}, and
 \begin{align*}
  E_{\bv,\lambda}(t)
  &:= \EE\biggl(\biggl|\ee^{\lambda t} \langle\bv, \bX_0\rangle
                       + \langle\bv, \tBbeta\rangle
                         \int_0^t \ee^{\lambda(t-u)} \, \dd u\biggr|^2\biggr) ,\\
  I_{\lambda,\ell}(t)
  &:= \int_0^t \ee^{2\Re(\lambda)(t-u)} \EE(X_{u,\ell}) \, \dd u , \qquad
   \ell \in \{1, \ldots, d\} , \\
  I_\lambda(t)
  &:= \int_0^t \ee^{2\Re(\lambda)(t-u)} \, \dd u .
 \end{align*}
If, in addition, \ $(\bX_t)_{t\in\RR_+}$ \ is supercritical and irreducible, then
 we have
 \[
   \lim_{t\to\infty} h(t) \EE(|\langle\bv, \bX_t\rangle|^2) = M_\bv^{(2)} ,
 \]
 where
 \begin{gather*}
  h(t) := \begin{cases}
           \ee^{-s(\tbB)t}
            & \text{if \ $\Re(\lambda)
                          \in \bigl(-\infty, \frac{1}{2} s(\tbB)\bigr)$,} \\
           t^{-1} \ee^{-s(\tbB)t}
            & \text{if \ $\Re(\lambda) = \frac{1}{2} s(\tbB)$,} \\
           \ee^{-2\Re(\lambda)t}
            & \text{if \ $\Re(\lambda)
                          \in \bigl(\frac{1}{2} s(\tbB), s(\tbB)\bigr]$,}
          \end{cases}
 \end{gather*}
 and
 \begin{gather*}
  M_\bv^{(2)}
  := \begin{cases}
      \frac{1}{s(\tbB)-2\Re(\lambda)}
      \bigl(\langle\bu, \EE(\bX_0)\rangle
            + \frac{\langle\bu, \tBbeta\rangle}{s(\tbB)}\bigr)
      \sum_{\ell=1}^d C_{\bv,\ell} \langle\be_\ell, \tbu\rangle
       & \text{if \ $\Re(\lambda)
                     \in \bigl(-\infty, \frac{1}{2} s(\tbB)\bigr)$,} \\[2mm]
      \bigl(\langle\bu, \EE(\bX_0)\rangle
            + \frac{\langle\bu, \tBbeta\rangle}{s(\tbB)}\bigr)
      \sum_{\ell=1}^d C_{\bv,\ell} \langle\be_\ell, \tbu\rangle
       & \text{if \ $\Re(\lambda) = \frac{1}{2} s(\tbB)$,} \\[2mm]
      \EE\bigl(\bigl|\langle\bv, \bX_0\rangle
                     + \frac{\langle\bv, \tBbeta\rangle}{\lambda}\bigr|^2\bigr)
      + \frac{1}{2\Re(\lambda)}
        \int_{\cU_d} |\langle\bv, \br\rangle|^2 \, \nu(\dd\br) \\
      + \sum\limits_{\ell=1}^d
         C_{\bv,\ell} \be_\ell^\top (2 \Re(\lambda) \bI_d - \tbB)^{-1}
         \bigl(\EE(\bX_0) + \frac{\tBbeta}{2\Re(\lambda)}\bigr)
       & \text{if \ $\Re(\lambda) \in \bigl(\frac{1}{2} s(\tbB), s(\tbB)\bigr]$.}
     \end{cases}
 \end{gather*}
\end{Pro}

Based on Proposition \ref{second_moment_asymptotics_CBI}, we derive the asymptotic
 behavior of the variance matrix of the real and imaginary parts of the projection
 of a multi-type CBI process on certain left eigenvectors of its branching mean
 matrix \ $\ee^\tbB$.

\begin{Pro}\label{variance_asymptotics_CBI}
If \ $(\bX_t)_{t\in\RR_+}$ \ is a supercritical and irreducible multi-type CBI
 process with parameters \ $(d, \bc, \Bbeta, \bB, \nu, \bmu)$ \ such that
 \ $\EE(\|\bX_0\|^2) < \infty$ \ and the moment condition
 \eqref{moment_condition_CBI2} holds, then for each left eigenvector
 \ $\bv \in \CC^d$ \ of \ $\tbB$ \ corresponding to an arbitrary eigenvalue
 \ $\lambda \in \sigma(\tbB)$ \ with
 \ $\Re(\lambda) \in \bigl(-\infty, \frac{1}{2} s(\tbB)\bigr]$ \ we have
 \[
   \lim_{t\to\infty}
    h(t)
    \EE\left(\begin{pmatrix}
              \Re(\langle\bv, \bX_t\rangle) \\
              \Im(\langle\bv, \bX_t\rangle)
             \end{pmatrix}
             \begin{pmatrix}
              \Re(\langle\bv, \bX_t\rangle) \\
              \Im(\langle\bv, \bX_t\rangle)
             \end{pmatrix}^\top\right)
   = \biggl(\langle\bu, \EE(\bX_0)\rangle
            + \frac{\langle\bu,\tBbeta\rangle}{s(\tbB)}\biggr)
     \bSigma_\bv ,
 \]
 where the scaling factor \ $h : \RR_{++} \to \RR_{++}$ \ and the matrix
 \ $\bSigma_\bv$ \ are defined in Proposition \ref{second_moment_asymptotics_CBI}
 and in Theorem \ref{convCBIweak1}, respectively.
\end{Pro}

\noindent
\textbf{Proof.}
For each \ $t \in \RR_+$, \ using the identity \eqref{help12_identity} for
 \ $a = \langle\bv, \bX_t\rangle \in \CC$, \ and then taking expectation, we obtain
 \begin{equation}\label{ReIm}
  \begin{aligned}
   &\EE\left(\begin{pmatrix}
              \Re(\langle\bv, \bX_t\rangle) \\
              \Im(\langle\bv, \bX_t\rangle)
             \end{pmatrix}
             \begin{pmatrix}
              \Re(\langle\bv, \bX_t\rangle) \\
              \Im(\langle\bv, \bX_t\rangle)
             \end{pmatrix}^\top\right) \\
   &= \frac{1}{2} \EE(|\langle\bv, \bX_t\rangle|^2) \bI_2
      + \frac{1}{2}
        \begin{pmatrix}
         \Re(\EE(\langle\bv, \bX_t\rangle^2))
          & \Im(\EE(\langle\bv, \bX_t\rangle^2)) \\
         \Im(\EE(\langle\bv, \bX_t\rangle^2))
          & - \Re(\EE(\langle\bv, \bX_t\rangle^2))
        \end{pmatrix} .
  \end{aligned}
 \end{equation}
The asymptotic behavior of \ $\EE(|\langle\bv, \bX_t\rangle|^2)$ \ as
 \ $t \to \infty$ \ is described in Proposition
 \ref{second_moment_asymptotics_CBI}.
The aim of the following discussion is to describe the asymptotic behavior of
 \ $\EE((\langle\bv, \bX_t\rangle)^2)$ \ as \ $t \to \infty$.
\ For each \ $t \in \RR_+$, \ we use the representation  of \ $\ee^{-\lambda t} \langle \bv,\bX_t\rangle$ \ given
 at the beginning of the proof of part (iii) of Theorem \ref{convCBIweak1}.
The independence of \ $\bX_0$, \ $(W_{u,1})_{u\in\RR_+}$, \ldots,
 \ $(W_{u,d})_{u\in\RR_+}$, \ $N_1$, \ldots, \ $N_d$ \ and \ $M$ \ implies the
 conditional independence of the random variables \ $Z_t^{(0,1)}$, \ $Z_t^{(2)}$, \ $Z_t^{(3,4)}$ \ and
  \ $Z_t^{(5)}$ \ with respect to \ $(\bX_u)_{u\in[0,t]}$ \ for every \ $t \in \RR_+$.
\ Moreover, the conditional expectations of \ $Z_t^{(2)}$, \ $Z_t^{(3,4)}$ \ and
  \ $Z_t^{(5)}$ \ with respect to \ $(\bX_u)_{u\in[0,t]}$ \ are \ $0$, \ since the processes
 \ $(Z_t^{(2)})_{t\in\RR_+}$, \ $(Z_t^{(3,4)})_{t\in\RR_+}$ \ and \ $(Z_t^{(5)})_{t\in\RR_+}$
 \ are martingales with initial values \ $0$.
\ Consequently, for all \ $t\in\RR_+$, \ we get
 \begin{align*}
   \EE\bigl((\ee^{-\lambda t} \langle\bv, \bX_t\rangle)^2
            \,\big|\, (\bX_u)_{u\in[0,t]}\bigr)
   &=  \EE\bigl(\bigl(Z_t^{(0,1)}\bigr)^2 \,\big|\, (\bX_u)_{u\in[0,t]}\bigr)
      +  \EE\bigl(\bigl(Z_t^{(2)}\bigr)^2 \,\big|\, (\bX_u)_{u\in[0,t]}\bigr) \\
   &\phantom{=} + \EE\bigl(\bigl(Z_t^{(3,4)}\bigr)^2 \,\big|\, (\bX_u)_{u\in[0,t]}\bigr)
      +  \EE\bigl(\bigl(Z_t^{(5)}\bigr)^2 \,\big|\, (\bX_u)_{u\in[0,t]}\bigr)
 \end{align*}
 almost surely.
We have
 \begin{align*}
  \EE\bigl(\bigl(Z_t^{(0,1)}\bigr)^2 \,\big|\, (\bX_u)_{u\in[0,t]}\bigr)
  &= \biggl(\langle\bv, \bX_0\rangle
            + \langle\bv, \tBbeta\rangle
              \int_0^t \ee^{-\lambda u} \, \dd u\biggr)^2 , \\
  \EE\bigl(\bigl(Z_t^{(2)}\bigr)^2 \,\big|\, (\bX_u)_{u\in[0,t]})
  &= 2 \sum_{\ell=1}^d
        \langle\bv, \be_\ell\rangle^2 c_\ell
        \int_0^t \ee^{-2\lambda u} X_{u,\ell} \, \dd u , \\
  \EE\bigl(\bigl(Z_t^{(3,4)}\bigr)^2 \,\big|\, (\bX_u)_{u\in[0,t]}\bigr)
  &= \sum_{\ell=1}^d
      \int_0^t \ee^{-2\lambda u} X_{u,\ell} \, \dd u
      \int_{\cU_d} \langle\bv, \bz\rangle^2 \, \mu_\ell(\dd\bz) , \\
  \EE\bigl(\bigl(Z_t^{(5)}\bigr)^2 \,\big|\, (\bX_u)_{u\in[0,t]}\bigr)
  &= \int_0^t \ee^{-2\lambda u} \, \dd u
     \int_{\cU_d} \langle\bv, \br\rangle^2 \, \nu(\dd\br)
 \end{align*}
 almost surely.
Taking the expectation and multiplying by \ $\ee^{2\lambda t}$, \ $t \in \RR_+$,
 \ we obtain
 \[
   \EE(\langle\bv, \bX_t\rangle^2)
   = \tE_{\bv,\lambda}(t) + \sum_{\ell=1}^d \tC_{\bv,\ell} \tI_{\lambda,\ell}(t)
     + \tI_\lambda(t) \int_{\cU_d} \langle\bv, \br\rangle^2 \, \nu(\dd\br)
 \]
 with
 \begin{align*}
  \tE_{\bv,\lambda}(t)
  &:= \EE\biggl(\biggl(\ee^{\lambda t} \langle\bv, \bX_0\rangle
                       + \langle\bv, \tBbeta\rangle
                         \int_0^t \ee^{\lambda(t-u)} \, \dd u\biggr)^2\biggr) , \\
  \tI_{\lambda,\ell}(t)
  &:= \int_0^t \ee^{2\lambda(t-u)} \EE(X_{u,\ell}) \, \dd u , \qquad
   \ell \in \{1, \ldots, d\} , \\
  \tI_\lambda(t)
  &:= \int_0^t \ee^{2\lambda(t-u)} \, \dd u ,
 \end{align*}
 and \ $\tC_{\bv,\ell}$, \ $\ell \in \{1, \ldots, d\}$ \ defined in Theorem
 \ref{convCBIweak1}.
For each \ $t \in \RR_+$, \ we have
 \[
   |\tE_{\bv,\lambda}(t)|
   \leq \EE\biggl(\biggl|\ee^{\lambda t} \langle\bv, \bX_0\rangle
                         + \langle\bv, \tBbeta\rangle
                           \int_0^t \ee^{\lambda(t-u)} \, \dd u\biggr|^2\biggr)
   = E_{\bv,\lambda}(t) .
 \]
If \ $\Re(\lambda) \in \bigl(-\infty, \frac{1}{2} s(\tbB)\bigr]$, \ then
 \ $h(t) E_{\bv,\lambda}(t) \to 0$ \ as \ $t \to \infty$, \ see the proof of
 Proposition B.1 in Barczy et al.\
 \cite{BarPalPap}, hence
 \begin{align}\label{help1}
   h(t) \tE_{\bv,\lambda}(t) \to 0 \qquad \text{as \ $t \to \infty$.}
 \end{align}
Moreover, for each \ $t \in \RR_+$ \ and \ $\ell \in \{1, \ldots, d\}$, \ by
 formula \eqref{EXcond}, we get
 \[
   \tI_{\lambda,\ell}(t)
   = \be_\ell^\top \tbA_{\lambda,1}(t) \EE(\bX_0)
     + \be_\ell^\top \tbA_{\lambda,2}(t) \tBbeta
 \]
 with
 \[
   \tbA_{\lambda,1}(t)
   := \int_0^t \ee^{2\lambda(t-u)} \ee^{u\tbB} \, \dd u , \qquad
   \tbA_{\lambda,2}(t)
   := \int_0^t
       \ee^{2\lambda(t-u)} \biggl(\int_0^u \ee^{w\tbB} \, \dd w\biggr) \dd u .
 \]
We have
 \[
   \tbA_{\lambda,1}(t)
   = \ee^{2\lambda t} \tbA_{\lambda,1,1}(t)
     + \ee^{2\lambda t} \tbA_{\lambda,1,2}(t), \qquad t \in \RR_+ ,
 \]
 with
 \begin{align*}
  \tbA_{\lambda,1,1}(t)
  := \int_0^t \ee^{(s(\tbB)-2\lambda)u} \, \tbu \bu^\top \dd u , \qquad
  \tbA_{\lambda,1,2}(t)
  := \int_0^t
      \ee^{(s(\tbB)-2\lambda)u} (\ee^{-s(\tbB)u} \ee^{u\tbB} - \tbu \bu^\top)
      \, \dd u .
 \end{align*}
If \ $\Re(\lambda) \in \bigl(-\infty, \frac{1}{2} s(\tbB)\bigr)$, \ then we have
 \[
   \ee^{-(s(\tbB)-2\lambda)t} \tbA_{\lambda,1,1}(t)
   = \ee^{-(s(\tbB)-2\lambda)t}
     \frac{\ee^{(s(\tbB)-2\lambda)t}-1}{s(\tbB)-2\lambda} \tbu \bu^\top
   = \frac{1-\ee^{-(s(\tbB)-2\lambda)t}}{s(\tbB)-2\lambda} \tbu \bu^\top
   \to \frac{\tbu \bu^\top}{s(\tbB)-2\lambda}
 \]
 as \ $t \to \infty$, \ and, by \eqref{C},
 \begin{align*}
  |\ee^{-(s(\tbB)-2\lambda)t} \tbA_{\lambda,1,2}(t)|
  &\leq C_1 \ee^{-(s(\tbB)-2\Re(\lambda))t}
        \int_0^t \ee^{(s(\tbB)-2\Re(\lambda))u} \ee^{-C_2 u} \, \dd u \\
  &\leq C_1 \ee^{-(s(\tbB)-2\Re(\lambda))t}
        \int_0^t \ee^{(s(\tbB)-2\Re(\lambda)-\tC_2)u} \, \dd u \\
  &\leq C_1 \ee^{-(s(\tbB)-2\Re(\lambda))t}
        \int_0^\infty \ee^{(s(\tbB)-2\Re(\lambda)-\tC_2)u} \, \dd u \\
  &= \frac{C_1}{s(\tbB)-2\Re(\lambda)-\tC_2} \ee^{-(s(\tbB)-2\Re(\lambda))t}
   \to 0
 \end{align*}
 as \ $t \to \infty$, \ where
 \ $\tC_2 \in (0, C_2 \land (s(\tbB) - 2\Re(\lambda)))$.
\ Hence, if \ $\Re(\lambda) \in \bigl(-\infty, \frac{1}{2} s(\tbB)\bigr)$, \ then
 \begin{align}\label{help2}
  h(t) \tbA_{\lambda,1}(t)
  = \ee^{-s(\tbB)t} \tbA_{\lambda,1}(t)
  \to \frac{\tbu \bu^\top}{s(\tbB)-2\lambda}
  \qquad \text{as \ $t \to \infty$.}
 \end{align}
If \ $\lambda = 0$, \ then, by Fubini's theorem, we obtain
 \[
   h(t) \tbA_{\lambda,2}(t)
   = \ee^{-s(\tbB)t} \tbA_{\lambda,2}(t)
   = \ee^{-s(\tbB)t} \int_0^t (t - w) \ee^{w\tbB} \, \dd w
   \to \frac{\tbu \bu^\top}{s(\tbB)^2}  \qquad
   \text{as \ $t \to \infty$,}
 \]
 see the proof of Proposition B.1 in Barczy et al.\
 \cite{BarPalPap}.
Hence, if \ $\lambda = 0$, \ then
 \begin{align}\label{help3}
  \begin{split}
   h(t) \tI_{\lambda,\ell}(t)
   =\ee^{-s(\tbB)t} \tI_{\lambda,\ell}(t)
   &\to \frac{1}{s(\tbB)} \be_\ell^\top \tbu \bu^\top
       \biggl(\EE(\bX_0) + \frac{\tBbeta}{s(\tbB)}\biggr)\\
   &= \frac{\langle\be_\ell, \tbu\rangle}{s(\tbB)}
     \biggl(\langle\bu, \EE(\bX_0)\rangle
            + \frac{\langle\bu, \tBbeta\rangle}{s(\tbB)}\biggr)
  \end{split}
 \end{align}
 as \ $t \to \infty$.
\ If \ $\lambda \in \sigma(\tbB) \setminus \{0\}$ \ with
 \ $\Re(\lambda) \in (-\infty, \frac{1}{2} s(\tbB))$, \ then, by Fubini's theorem,
 we obtain
 \begin{align*}
  \ee^{-s(\tbB)t} \tbA_{\lambda,2}(t)
  &= \ee^{-(s(\tbB)-2\lambda)t}
     \int_0^t
      \ee^{-2\lambda u} \biggl(\int_0^u \ee^{w\tbB} \, \dd w\biggr) \dd u \\
  &= \ee^{-(s(\tbB)-2\lambda)t}
     \int_0^t
      \biggl(\int_w^t \ee^{-2\lambda u} \, \dd u\biggr) \ee^{w\tbB} \, \dd w \\
  &= \frac{1}{2\lambda} \ee^{-(s(\tbB)-2\lambda)t}
     \int_0^t
      \bigl(\ee^{-2\lambda w} - \ee^{-2\lambda t}\bigr) \ee^{w\tbB} \, \dd w \\
  &= \frac{1}{2\lambda}
     \biggl(\ee^{-(s(\tbB)-2\lambda)t}
            \int_0^t \ee^{-2\lambda w} \ee^{w\tbB} \, \dd w
            - \ee^{-s(\tbB)t} \int_0^t \ee^{w\tbB} \, \dd w\biggr) \\
  &\to \frac{1}{2\lambda}
       \biggl(\frac{\tbu \bu^\top}{s(\tbB)-2\lambda}
              - \frac{\tbu \bu^\top}{s(\tbB)}\biggr)
   = \frac{\tbu \bu^\top}{(s(\tbB)-2\lambda)s(\tbB)}
 \end{align*}
 as \ $t \to \infty$, \ since
 \[
   \ee^{-(s(\tbB)-2\lambda)t} \int_0^t \ee^{-2\lambda w} \ee^{w\tbB} \, \dd w
   = \ee^{-s(\tbB)t} \tbA_{\lambda,1}(t)
   \to \frac{\tbu \bu^\top}{s(\tbB)-2\lambda} , \qquad
   \ee^{-s(\tbB)t} \int_0^t \ee^{w\tbB} \, \dd w
   \to \frac{\tbu \bu^\top}{s(\tbB)}
 \]
 as \ $t \to \infty$, \ by \eqref{help2} and the proof of Proposition B.1 in Barczy
 et al.\ \cite{BarPap}.
Hence, if \ $\lambda \in \sigma(\tbB) \setminus \{0\}$ \ with
 \ $\Re(\lambda) \in (-\infty, \frac{1}{2} s(\tbB))$, \ then
 \begin{align}\label{help4}
 \begin{split}
  h(t)\tI_{\lambda,\ell}(t)
  = \ee^{-s(\tbB)t} \tI_{\lambda,\ell}(t)
  &\to \frac{1}{s(\tbB)-2\lambda} \be_\ell^\top \tbu \bu^\top
       \biggl(\EE(\bX_0) + \frac{\tBbeta}{s(\tbB)}\biggr) \\
  &= \frac{\langle\be_\ell, \tbu\rangle}{s(\tbB)-2\lambda}
     \biggl(\langle\bu, \EE(\bX_0)\rangle
            + \frac{\langle\bu, \tBbeta\rangle}{s(\tbB)}\biggr) \qquad
  \text{as \ $t \to \infty$.}
 \end{split}
 \end{align}
If \ $\Re(\lambda) = \frac{1}{2} s(\tbB)$ \ and \ $\Im(\lambda) = 0$, \ then we
 have
 \[
   t^{-1} \ee^{-(s(\tbB)-2\lambda)t}  \tbA_{\lambda,1,1}(t)
   = t^{-1} \int_0^t \tbu \bu^\top \dd u
   = \tbu \bu^\top , \qquad t \in \RR_+
 \]
 and, by \eqref{C},
 \begin{align*}
  |t^{-1} \ee^{-(s(\tbB)-2\lambda)t} \tbA_{\lambda,1,2}(t)|
  &= |t^{-1} \tbA_{\lambda,1,2}(t)|
   \leq C_1 t^{-1} \int_0^t \ee^{(s(\tbB)-2\Re(\lambda))u} \ee^{-C_2 u} \, \dd u \\
  &\leq C_1 t^{-1} \int_0^\infty \ee^{-C_2 u} \, \dd u
   = \frac{C_1}{C_2} t^{-1}
   \to 0 \qquad \text{as \ $t \to \infty$.}
 \end{align*}
Hence, if \ $\Re(\lambda) = \frac{1}{2} s(\tbB)$ \ and \ $\Im(\lambda) = 0$, \ then
 \[
   h(t) \tbA_{\lambda,1}(t)
   = t^{-1} \ee^{-s(\tbB)t} \tbA_{\lambda,1}(t)
   \to \tbu \bu^\top \qquad \text{as \ $t \to \infty$.}
 \]
If \ $\Re(\lambda) = \frac{1}{2} s(\tbB)$ \ and \ $\Im(\lambda) = 0$, \ then, by
 Fubini's theorem, we obtain
 \begin{align*}
  t^{-1} \ee^{-s(\tbB)t} \tbA_{\lambda,2}(t)
  &= t^{-1}
     \int_0^t \ee^{-s(\tbB)u} \biggl(\int_0^u \ee^{w\tbB} \, \dd w\biggr) \dd u
   = t^{-1}
     \int_0^t
      \biggl(\int_w^t \ee^{-s(\tbB)u} \, \dd u\biggr) \ee^{w\tbB} \, \dd w \\
  &= \frac{1}{s(\tbB)} t^{-1}
     \int_0^t \bigl(\ee^{-s(\tbB)w} - \ee^{-s(\tbB)t}\bigr) \ee^{w\tbB} \, \dd w \\
  &= \frac{1}{s(\tbB)} t^{-1}
     \biggl(\int_0^t \ee^{-s(\tbB)w} \ee^{w\tbB} \, \dd w
            - \ee^{-s(\tbB)t} \int_0^t \ee^{w\tbB} \, \dd w\biggr)
   \to \frac{\tbu \bu^\top}{s(\tbB)}
 \end{align*}
 as \ $t \to \infty$, \ since
 \[
   t^{-1} \int_0^t \ee^{-s(\tbB)w} \ee^{w\tbB} \, \dd w
   \to \tbu \bu^\top , \qquad
   \ee^{-s(\tbB)t} \int_0^t \ee^{w\tbB} \, \dd w
   \to \frac{\tbu \bu^\top}{s(\tbB)} \qquad
   \text{as \ $t\to\infty$,}
 \]
 see part (v) of Lemma A.2 and the proof of Proposition B.1 in
 Barczy et al.\ \cite{BarPap}.
Consequently, if \ $\Re(\lambda) = \frac{1}{2} s(\tbB)$ \ and \ $\Im(\lambda) = 0$,
 \ then
 \begin{align}\label{help5}
   t^{-1} \ee^{-s(\tbB)t} \tI_{\lambda,\ell}(t)
   \to \be_\ell^\top \tbu \bu^\top
       \biggl(\EE(\bX_0) + \frac{\tBbeta}{s(\tbB)}\biggr)
   = \langle\be_\ell, \tbu\rangle
     \biggl(\langle\bu, \EE(\bX_0)\rangle
            + \frac{\langle\bu, \tBbeta\rangle}{s(\tbB)}\biggr)
 \end{align}
 as \ $t \to \infty$.

If \ $\Re(\lambda) = \frac{1}{2} s(\tbB)$ \ and \ $\Im(\lambda) \ne 0$, \ then we
 have
 \begin{align*}
   t^{-1} \ee^{-(s(\tbB)-2\lambda)t} \tbA_{\lambda,1,1}(t)
  & = t^{-1} \ee^{-(s(\tbB)-2\lambda)t}
     \frac{\ee^{(s(\tbB)-2\lambda)t}-1}{s(\tbB)-2\lambda} \tbu \bu^\top \\
  & = \frac{1}{(s(\tbB)-2\lambda)t} (1 - \ee^{2\ii\Im(\lambda)t})
     \tbu \bu^\top
   \to \bzero
 \end{align*}
 as \ $t \to \infty$ \ and
 \begin{align*}
  |t^{-1} \ee^{-(s(\tbB)-2\lambda)t} \tbA_{\lambda,1,2}(t)|
  &\leq C_1 t^{-1} \ee^{-(s(\tbB)-2\Re(\lambda))t}
        \int_0^t \ee^{(s(\tbB)-2\Re(\lambda))u} \ee^{-C_2 u} \, \dd u \\
  &\leq C_1 t^{-1} \int_0^\infty \ee^{-C_2 u} \, \dd u
   = \frac{C_1}{C_2} t^{-1}
   \to 0 \qquad \text{as \ $t \to \infty$.}
 \end{align*}
Hence, if \ $\Re(\lambda) = \frac{1}{2} s(\tbB)$ \ and \ $\Im(\lambda) \ne 0$, \ then
 \[
   h(t) \tbA_{\lambda,1}(t)
   = t^{-1} \ee^{-s(\tbB)t} \tbA_{\lambda,1}(t)
   \to \bzero \qquad \text{as \ $t \to \infty$.}
 \]
If \ $\Re(\lambda) = \frac{1}{2} s(\tbB)$ \ and \ $\Im(\lambda) \ne 0$, \ then, by
 Fubini's theorem, as above, we obtain
 \begin{align*}
  t^{-1} \ee^{-s(\tbB)t} \tbA_{\lambda,2}(t)
  &= t^{-1} \ee^{-(s(\tbB)-2\lambda)t}
     \int_0^t
      \biggl(\int_w^t \ee^{-2\lambda u} \, \dd u\biggr) \ee^{w\tbB} \, \dd w \\
  &= \frac{1}{2\lambda t} \ee^{-(s(\tbB)-2\lambda)t}
     \int_0^t (\ee^{-2\lambda w} - \ee^{-2\lambda t}) \ee^{w\tbB} \, \dd w \\
  &= \frac{1}{2\lambda t}
     \biggl(\ee^{-(s(\tbB)-2\lambda)t}
            \int_0^t \ee^{-2\lambda w} \ee^{w\tbB} \, \dd w
            - \ee^{-s(\tbB)t} \int_0^t \ee^{w\tbB} \, \dd w\biggr)
   \to \bzero
 \end{align*}
 as \ $t \to \infty$.
\ Indeed,
 \ $\ee^{-s(\tbB)t} \int_0^t \ee^{w\tbB} \, \dd w \to \frac{\tbu\bu^\top}{s(\tbB)}$
 \ as \ $t \to \infty$, \ and using that \ $\Re(\lambda) = \frac{1}{2} s(\tbB)$
 \ and \ $\Im(\lambda) \ne 0$, \ for all \ $t \in \RR_+$ \ we have
 \begin{align*}
  & t^{-1} \ee^{-(s(\tbB)-2\lambda)t}
   \int_0^t \ee^{-2\lambda w} \ee^{w\tbB} \, \dd w
   = t^{-1} \ee^{2\ii\Im(\lambda)t}
     \int_0^t \ee^{-2\ii\Im(\lambda)w} \ee^{-s(\tbB)w} \ee^{w\tbB} \, \dd w \\
  &= t^{-1} \ee^{2\ii\Im(\lambda)t}
     \biggl(\int_0^t \ee^{-2\ii\Im(\lambda)w} \, \dd w\biggr) \tbu \bu^\top
     + t^{-1} \ee^{2\ii\Im(\lambda)t}
       \int_0^t
        \ee^{-2\ii\Im(\lambda)w} (\ee^{-s(\tbB)w} \ee^{w\tbB} - \tbu \bu^\top)
        \, \dd w ,
 \end{align*}
 where
 \begin{align*}
  \biggl|t^{-1} \ee^{2\ii\Im(\lambda)t}
         \int_0^t \ee^{-2\ii\Im(\lambda)w} \, \dd w \biggr|
  = t^{-1}
    \biggl|\frac{\ee^{-2\ii\Im(\lambda)t}-1}{-2\ii\Im(\lambda)}\biggr|
  \leq \frac{1}{t|\Im(\lambda)|} \to 0
 \end{align*}
 as \ $t \to \infty$, \ and, by \eqref{C},
 \begin{align*}
  &\biggl\|t^{-1} \ee^{2\ii\Im(\lambda)t}
           \int_0^t
            \ee^{-2\ii\Im(\lambda)w}
            (\ee^{-s(\tbB) w} \ee^{w\tbB} - \tbu \bu^\top)
            \, \dd w\biggr\|
   \leq t^{-1} \int_0^t \|\ee^{-s(\tbB)w} \ee^{w\tbB} - \tbu \bu^\top\| \, \dd w \\
  &\leq t^{-1} C_1 \int_0^t \ee^{-C_2 w} \, \dd w
   \leq t^{-1} C_1 \int_0^\infty \ee^{-C_2 w} \, \dd w
   = \frac{C_1}{C_2 t}
   \to 0 \qquad \text{as \ $t \to \infty$.}
 \end{align*}
Consequently, if \ $\Re(\lambda) = \frac{1}{2} s(\tbB)$ \ and
 \ $\Im(\lambda) \ne 0$, \ then
 \begin{align}\label{help6}
  h(t) \tI_{\lambda,\ell}(t)
  = t^{-1} \ee^{-s(\tbB)t} \tI_{\lambda,\ell}(t)
  \to 0 \qquad \text{as \ $t \to \infty$.}
 \end{align}
By the help of \eqref{help3}, \eqref{help4}, \eqref{help5} and
 \eqref{help6}, we have
 \begin{align}\label{help7}
  \lim_{t\to\infty} h(t) \tI_{\lambda,\ell}(t)
  = \begin{cases}
     \frac{\langle\be_\ell, \tbu\rangle}{s(\tbB)-2\lambda}
           \bigl(\langle\bu, \EE(\bX_0)\rangle
                 + \frac{\langle\bu, \tBbeta\rangle}{s(\tbB)}\bigr)
      &\text{if \ $\Re(\lambda) \in \bigl(-\infty, \frac{1}{2} s(\tbB)\bigr)$,} \\
     \langle\be_\ell, \tbu\rangle
     \bigl(\langle\bu, \EE(\bX_0)\rangle
           + \frac{\langle\bu, \tBbeta\rangle}{s(\tbB)}\bigr)
      &\text{if \ $\Re(\lambda) = \frac{1}{2} s(\tbB)$ \ and
             \ $\Im(\lambda) = 0$,} \\
     0 &\text{if \ $\Re(\lambda) = \frac{1}{2} s(\tbB)$ \ and
              \ $\Im(\lambda) \ne 0$.}
  \end{cases}
 \end{align}
Further, we have
 \[
   \tI_\lambda(t)
   = \int_0^t \ee^{2\lambda w} \, \dd w
   = \begin{cases}
      t & \text{if \ $\lambda = 0$,} \\
      \frac{1}{2\lambda} (\ee^{2\lambda t} - 1)
       & \text{if \ $\lambda \ne 0$.}
     \end{cases}
 \]
 If \ $\lambda = 0$, \ then
 \[
   \ee^{-s(\tbB)t} \tI_\lambda(t) = t \ee^{-s(\tbB)t} \to 0 \qquad
   \text{as \ $t \to \infty$.}
 \]
If \ $\lambda \in \sigma(\tbB) \setminus \{0\}$ \ with
 \ $\Re(\lambda) \in \bigl(-\infty, \frac{1}{2} s(\tbB)\bigr)$, \ then
 \[
   \ee^{-s(\tbB)t} \tI_\lambda(t)
   = \frac{1}{2\lambda} (\ee^{-(s(\tbB)-2\lambda)t} - \ee^{-s(\tbB)t})
   \to 0 \qquad \text{as \ $t \to \infty$.}
 \]
If \ $\Re(\lambda) = \frac{1}{2} s(\tbB)$, \ then
 \[
   t^{-1} \ee^{-s(\tbB)t} \tI_\lambda(t)
   = \frac{1}{2\lambda t} (\ee^{2\ii\Im(\lambda)t} - \ee^{-s(\tbB)t})
   \to 0 \qquad \text{as \ $t \to \infty$.}
 \]
Consequently,
 \begin{align}\label{help8}
  \lim_{t\to\infty} h(t) \tI_\lambda(t) = 0 .
 \end{align}
Hence, by \eqref{help1}, \eqref{help7} and \eqref{help8}, we have
 \[
   \lim_{t\to\infty} h(t) \EE(\langle\bv,\bX_t\rangle^2) = \tM_{\bv}^{(2)}
 \]
 with
 \[
   \tM_{\bv}^{(2)}
   := \begin{cases}
      \sum_{\ell=1}^d
       \tC_{\bv,\ell} \frac{\langle\be_\ell, \tbu\rangle}{s(\tbB)-2\lambda}
       \biggl(\langle\bu, \EE(\bX_0)\rangle
              + \frac{\langle\bu, \tBbeta\rangle}{s(\tbB)}\biggr)
       &\text{if \ $\Re(\lambda) \in \bigl(-\infty, \frac{1}{2} s(\tbB)\bigr)$,} \\
      \sum_{\ell=1}^d
       \tC_{\bv,\ell} \langle\be_\ell, \tbu\rangle
       \biggl(\langle\bu, \EE(\bX_0)\rangle
              + \frac{\langle\bu, \tBbeta\rangle}{s(\tbB)}\biggr)
       &\text{if \ $\Re(\lambda) = \frac{1}{2} s(\tbB)$ \ and
              \ $\Im(\lambda) = 0$,} \\
      0 &\text{if \ $\Re(\lambda) = \frac{1}{2} s(\tbB)$ \ and
               \ $\Im(\lambda) \ne 0$.}
   \end{cases}
 \]
Using the identity \eqref{ReIm}, then taking the limit as \ $t \to \infty$, \ and
 using Proposition \ref{second_moment_asymptotics_CBI}, we obtain
  \begin{align*}
   h(t)
   \EE\left(\begin{pmatrix}
             \Re(\langle\bv, \bX_t\rangle) \\
             \Im(\langle\bv, \bX_t\rangle)
            \end{pmatrix}
            \begin{pmatrix}
             \Re(\langle\bv, \bX_t\rangle) \\
             \Im(\langle\bv, \bX_t\rangle)
            \end{pmatrix}^\top\right)
   &\to \frac{1}{2} M_{\bv}^{(2)} \bI_2
        + \frac{1}{2}
          \begin{pmatrix}
           \Re(\tM_{\bv}^{(2)}) & \Im(\tM_{\bv}^{(2)}) \\
           \Im(\tM_{\bv}^{(2)}) & - \Re(\tM_{\bv}^{(2)})
          \end{pmatrix} \\
   &= \biggl(\langle\bu, \EE(\bX_0)\rangle
             + \frac{\langle\bu,\tBbeta\rangle}{s(\tbB)}\biggr)
      \bSigma_\bv
 \end{align*}
 as \ $t \to \infty$, \ as desired.
\proofend

\begin{Pro}\label{Pro_Sigma_invertible}
Let \ $(\bX_t)_{t\in\RR_+}$ \ be a supercritical and irreducible multi-type CBI
 process with parameters \ $(d, \bc, \Bbeta, \bB, \nu, \bmu)$ \ such that
 \ $\EE(\|\bX_0\|^2) < \infty$ \ and the moment condition
 \eqref{moment_condition_CBI2} holds.
Let \ $\lambda \in \sigma(\tbB)$ \ with
 \ $\Re(\lambda) \in \bigl(-\infty, \frac{1}{2} s(\tbB)\bigr]$ \ and
 \ $\bv \in \CC^d$ \ be a left-eigenvector of \ $\tbB$ \ corresponding to the
 eigenvalue \ $\lambda$.
\ Then \ $\bSigma_\bv = \bzero$ \ if and only if
 \ $c_\ell \langle\bv, \be_\ell\rangle = 0$ \ and
 \ $\mu_\ell(\{\bz \in \cU_d : \langle\bv, \bz\rangle \ne 0 \}) = 0$ \ for each
 \ $\ell \in \{1, \ldots, d\}$.
\ If, in addition, \ $\Im(\lambda) \ne 0$, \ then \ $\bSigma_\bv$ \ is invertible
 if and only if there exists \ $\ell \in \{1, \ldots, d\}$ \ such that
 \ $c_\ell \langle\bv, \be_\ell\rangle \ne 0$ \ or
 \ $\mu_\ell(\{\bz \in \cU_d : \langle\bv, \bz\rangle \ne 0 \}) > 0$.
\end{Pro}

\noindent {\bf Proof.}
First, suppose that \ $\Re(\lambda) \in \bigl(-\infty, \frac{1}{2} s(\tbB)\bigr)$.
\ If \ $c_\ell \langle\bv, \be_\ell\rangle = 0$ \ and
 \ $\mu_\ell(\{\bz \in \cU_d : \langle\bv, \bz\rangle \ne 0 \}) = 0$ \ for each
 \ $\ell \in \{1, \ldots, d\}$, \ then we have
 \ $C_{\bv,\ell} = \tC_{\bv,\ell} = 0$, \ $\ell \in \{1, \ldots, d\}$, \ yielding
 that \ $\bSigma_\bv = \bzero$.
\ If there exists an \ $\ell \in \{1, \ldots, d\}$ \ such that
 \ $c_\ell \langle\bv, \be_\ell\rangle \ne 0$ \ or
 \ $\mu_\ell(\{\bz \in \cU_d : \langle\bv, \bz\rangle \ne 0 \}) > 0$, \ then we
 check that \ $\bSigma_\bv \ne \bzero$, \ as requested.
On the contrary, let us suppose that \ $\bSigma_\bv = \bzero$.
\ Due to the existence of such an \ $\ell$, \ we have
 \ $C_{\bv,\ell} \in \RR_{++}$ \ and hence
 \ $\sum_{k=1}^d
     \langle \be_k, \tbu\rangle
     \frac{C_{\bv,k}}{s(\tbB)-2\Re(\lambda)}
     \in \RR_{++}$.
\ However, using the notation
 \ $\bSigma_\bv = ((\Sigma_\bv)_{i,j})_{i,j\in\{1,2\}}$, \ since
 \begin{align*}
  &(\Sigma_\bv)_{1,1}
   = \frac{1}{2}
     \sum_{k=1}^d \langle\be_k, \tbu\rangle \frac{C_{\bv,k}}{s(\tbB)-2\Re(\lambda)}
     + \frac{1}{2}
       \sum_{k=1}^d \langle \be_k, \tbu\rangle \Re\left( \frac{\tC_{\bv,k}}{s(\tbB) - 2\lambda}  \right) =0,\\
  &(\Sigma_\bv)_{2,2} =  \frac{1}{2} \sum_{k=1}^d \langle \be_k, \tbu\rangle  \frac{C_{\bv,k}}{s(\tbB) - 2\Re(\lambda)}
                          - \frac{1}{2} \sum_{k=1}^d \langle \be_k, \tbu\rangle \Re\left( \frac{\tC_{\bv,k}}{s(\tbB) - 2\lambda}  \right) =0,
 \end{align*}
 we have \ $\sum_{k=1}^d \langle \be_k, \tbu\rangle  \frac{C_{\bv,k}}{s(\tbB) - 2\Re(\lambda)} =0$, \ yielding us to a contradiction.

Next, suppose that \ $\Re(\lambda) = \frac{1}{2} s(\tbB)$.
\ If \ $c_\ell \langle\bv, \be_\ell\rangle = 0$ \ and \ $\mu_\ell(\{\bz \in \cU_d : \langle\bv, \bz\rangle \ne 0 \}) = 0$ \ for each
 \ $\ell \in \{1, \ldots, d\}$, \ then, as in case of \ $\Re(\lambda)\in(-\infty,\frac{1}{2}s(\tbB))$, \ we have
 \ $\bSigma_\bv = \bzero$.
\ If there exists an \ $\ell \in \{1, \ldots, d\}$ \ such that \ $c_\ell \langle\bv, \be_\ell\rangle \ne 0$ \ or
 \ $\mu_\ell(\{\bz \in \cU_d : \langle\bv, \bz\rangle \ne 0 \}) > 0$, \ then we check that \ $\bSigma_\bv \ne \bzero$, \ as requested.
On the contrary, let us suppose that \ $\bSigma_\bv = \bzero$.
Similarly, as in case of \ $\Re(\lambda)\in(-\infty,\frac{1}{2}s(\tbB))$, \ we have
 \ $\sum_{k=1}^d \langle \be_k, \tbu\rangle C_{\bv,k} =0$, \ yielding us to a contradiction.

Recall that, by \eqref{help13},
 \[
 \bSigma_\bv
  = 2 \sum_{\ell=1}^d
            c_\ell \langle\be_\ell, \tbu\rangle
            \int_0^\infty \Bf(w, \be_\ell) \, \dd w
   + \sum_{\ell=1}^d
            \langle\be_\ell, \tbu\rangle
            \int_0^\infty \int_{\cU_d}
            \Bf(w, \bz) \, \dd w \, \mu_\ell(\dd\bz)
   =: \bSigma_{\bv,1} + \bSigma_{\bv,2},
 \]
 where both \ $\bSigma_{\bv,1}$ \ and \ $\bSigma_{\bv,2}$ \ (and consequently \ $\bSigma_\bv$) \ are symmetric and non-negative definite matrices,
 since \ $\bc\in\RR_+^d$, \ $\tbu\in\RR_{++}^d$, \ and \ $\Bf(w,\bz)$ \ is symmetric and non-negative definite for any \ $w\in\RR_+$ \ and \ $\bz\in\cU_d$.

In what follows, let us assume that \ $\Re(\lambda)\in (-\infty, \frac{1}{2}s(\tbB)]$ \ and \ $\Im(\lambda) \ne 0$.
\ First, let us suppose that \ $\bSigma_\bv$ \ is invertible, and, on the contrary, for each \ $\ell\in\{1,\ldots,d\}$, \ we have
 \ $c_\ell \langle \bv, \be_\ell\rangle = 0$ \ and \ $\mu_\ell(\{ \bz\in\cU_d : \langle \bv, \bz\rangle \ne 0 \}) = 0$.
Then, for each \ $\ell\in\{1,\ldots,d\}$, \ we have \ $C_{\bv,\ell} = \tC_{\bv,\ell}=0$, \ and hence,
 by \eqref{help15_Sigma_v2} and \eqref{help15_Sigma_v}, \ $\bSigma_{\bv}=\bzero$, \ yielding us to a contradiction.

Let us suppose now that there exists \ $\ell\in\{1,\ldots,d\}$ \ such that
 \ $c_\ell \langle \bv, \be_\ell\rangle \ne 0$ \ or \ $\mu_\ell(\{ \bz\in\cU_d : \langle \bv, \bz\rangle \ne 0 \})>0$.
\ Next we show that if \ $\ell\in\{1,\ldots,d\}$ \ is such that \ $c_\ell \langle \bv, \be_\ell\rangle \ne 0$, \
 then \ $\bSigma_{\bv,1}$ \ is strictly positive definite, and that if \ $\ell\in\{1,\ldots,d\}$ \ is such that
 \ $\mu_\ell(\{ \bz\in\cU_d : \langle \bv, \bz\rangle \ne 0 \})>0$, \ then \ $\bSigma_{\bv,2}$ \ is strictly positive definite,
 yielding that \ $\bSigma_{\bv}$ \ is strictly positive definite, and consequently is invertible.
Here for all \ $w\in\RR_+$, \ $\bz\in\cU_d$, \ and \ $a,b\in\RR$, \ we have
 \begin{align*}
   \begin{pmatrix}
     a \\ b
   \end{pmatrix}^\top
    \Bf(w,\bz)
    \begin{pmatrix}
      a \\
      b
    \end{pmatrix}
  & = \Big( a\Re(\ee^{-(s(\tbB) - 2\lambda)w/2} \langle \bv,\bz\rangle)
           + b\Im(\ee^{-(s(\tbB) - 2\lambda)w/2} \langle \bv,\bz\rangle)\Big)^2 \\
  & = \big(\Re((a-\ii b)\ee^{-(s(\tbB) - 2\lambda)w/2} \langle \bv,\bz\rangle ) \big)^2.
 \end{align*}
Consequently, if \ $\ell\in\{1,\ldots,d\}$ \ is such that \ $\langle \bv, \be_\ell\rangle \ne 0$, \
 then for each \ $(a, b)^\top\in\RR^2\setminus \{\bzero\}$, \ we have
 \[
  \begin{pmatrix}
     a \\ b
   \end{pmatrix}^\top
   \int_0^\infty \Bf(w, \be_\ell) \, \dd w
    \begin{pmatrix}
      a \\
      b
    \end{pmatrix}
   = \int_0^\infty
       \big(\Re((a-\ii b)\ee^{-(s(\tbB) - 2\lambda)w/2} \langle \bv,\be_\ell\rangle ) \big)^2 \,\dd w
 \]
 is equal to \ $0$ \ if and only if \ $\Re((a-\ii b)\ee^{-(s(\tbB) - 2\lambda)w/2} \langle \bv,\be_\ell\rangle) = 0$ \
 for every \ $w\in\RR_+$ \ or equivalently \ $\Re(\ee^{\ii \Im(\lambda)w}(a-\ii b)\langle \bv,\be_\ell\rangle) = 0$ \
 for every \ $w\in\RR_+$.
Since \ $(a-\ii b)\langle \bv,\be_\ell\rangle\ne 0$ \ and \ $\Im(\lambda)\ne 0$, \ there exists \ $w\in\RR_+$ \ such that
 \ $\Re(\ee^{\ii \Im(\lambda)w}(a-\ii b)\langle \bv,\be_\ell\rangle) \ne 0$.
\ Indeed, the multiplication by the complex number \ $\ee^{\ii \Im(\lambda)w}$ \ corresponds to a rotation by degree \ $\Im(\lambda)w$.
Hence
  for each \ $(a, b)^\top\in\RR^2\setminus \{\bzero\}$, \ we have
 \[
  \begin{pmatrix}
     a \\ b
   \end{pmatrix}^\top
   \int_0^\infty \Bf(w, \be_\ell) \, \dd w
    \begin{pmatrix}
      a \\
      b
    \end{pmatrix}
   \in\RR_{++}.
 \]
This yields that if \ $\ell\in\{1,\ldots,d\}$ \ is such that \ $c_\ell\langle \bv, \be_\ell\rangle \ne 0$, \
 then for each \ $(a, b)^\top\in\RR^2\setminus \{\bzero\}$,
 \[
  \begin{pmatrix}
     a \\ b
   \end{pmatrix}^\top
   \bSigma_{\bv,1}
    \begin{pmatrix}
      a \\
      b
    \end{pmatrix}
   \in\RR_{++},
 \]
 implying that \ $\bSigma_{\bv,1}$ \ is strictly positive definite.
Further, for each \ $(a, b)^\top\in\RR^2$, \ we have
 \[
  \begin{pmatrix}
     a \\ b
   \end{pmatrix}^\top
   \int_0^\infty \int_{\cU_d} \Bf(w, \bz) \, \dd w\mu_\ell(\dd \bz)
    \begin{pmatrix}
      a \\
      b
    \end{pmatrix}
   = \int_0^\infty \int_{\cU_d}
       \big(\Re((a-\ii b)\ee^{-(s(\tbB) - 2\lambda)w/2} \langle \bv,\bz\rangle ) \big)^2 \,\dd w \mu_\ell(\dd \bz).
 \]
Since \ $\Im(\lambda)\ne0$, \ for each \ $\bz\in\cU_d$ \ with \ $\langle \bv,\bz\rangle \ne 0$ \ and
 \ $(a, b)^\top\in\RR^2\setminus \{\bzero\}$, \ there exists an open subset \ $K_{\bz}$ \ of \ $\RR_+$ \ such that
 \ $\big(\Re((a-\ii b)\ee^{-(s(\tbB) - 2\lambda)w/2} \langle \bv,\bz\rangle ) \big)^2 \in\RR_{++}$ \ for all \ $w\in K_{\bz}$.
\ Consequently, if \ $\ell\in\{1,\ldots,d\}$ \ is such that \ $\mu_\ell(\{ \bz\in\cU_d : \langle \bv, \bz\rangle \ne 0 \})>0$,
 \ then
 \begin{align*}
  &\int_0^\infty \int_{\cU_d}
   \big(\Re((a-\ii b)\ee^{-(s(\tbB) - 2\lambda)w/2} \langle \bv,\bz\rangle ) \big)^2 \,\dd w \mu_\ell(\dd \bz) \\
  &\geq
      \int_{\cU_d}
       \bbone_{\{\langle\bv,\bz\rangle\ne0\}}
      \left(\int_{K_\bz}
            \big(\Re((a-\ii b)\ee^{-(s(\tbB) - 2\lambda)w/2} \langle \bv,\bz\rangle ) \big)^2 \, \dd w \right) \mu_\ell(\dd \bz)
       \in\RR_{++}.
 \end{align*}
This yields that if \ $\ell\in\{1,\ldots,d\}$ \ is such that \ $\mu_\ell(\{ \bz\in\cU_d : \langle \bv, \bz\rangle \ne 0 \})>0$,
 \ then for each \ $(a, b)^\top\in\RR^2\setminus \{\bzero\}$, \ we have
  \[
  \begin{pmatrix}
     a \\ b
   \end{pmatrix}^\top
   \bSigma_{\bv,2}
    \begin{pmatrix}
      a \\
      b
    \end{pmatrix}
   \in\RR_{++},
 \]
 implying that \ $\bSigma_{\bv,2}$ \ is strictly positive definite.
\proofend

\begin{Rem}
Under the conditions of Proposition \ref{Pro_Sigma_invertible}, if
 \ $\lambda \in \sigma(\tbB)$ \ with
 \ $\Re(\lambda) \in \bigl(-\infty, \frac{1}{2} s(\tbB)\bigr]$ \ and
 \ $\Im(\lambda) = 0$ \ and \ $\bv \in \RR^d$ \ is a left eigenvector of \ $\tbB$
 \ corresponding to the eigenvalue \ $\lambda$, \ then \ $\bSigma_\bv$ \ is
 singular.
Indeed, in this case, by \eqref{help15_Sigma_v} and \eqref{help15_Sigma_v2}, we have
 \[
   \bSigma_\bv
   = \begin{cases}
      \begin{pmatrix}
       \frac{1}{s(\tbB)-2\lambda}
       \sum_{\ell=1}^d \langle \be_\ell, \tbu\rangle C_{\bv,\ell} & 0 \\
       0 & 0 \\
      \end{pmatrix}
       & \text{if \ $\Re(\lambda) \in \bigl(-\infty, \frac{1}{2} s(\tbB)\bigr)$,} \\[5mm]
      \begin{pmatrix}
       \sum_{\ell=1}^d \langle \be_\ell, \tbu\rangle C_{\bv,\ell} & 0 \\
       0 & 0 \\
      \end{pmatrix}
       & \text{if \ $\Re(\lambda) = \frac{1}{2} s(\tbB)$.}
     \end{cases}
 \]
Note that if \ $\bv \in \RR^d$ \ is a left eigenvector of \ $\tbB$ \ corresponding
 to an eigenvalue \ $\lambda$ \ of \ $\tbB$, \ then \ $\lambda \in \RR$
 \ necessarily, and hence in case of \ $\lambda\in(-\infty,\frac{1}{2}s(\tbB)]$, \ we have \ $\bSigma_\bv$ \ is not invertible.
However, if \ $\lambda \in \RR$ \ is an eigenvalue of \ $\tbB$ \ and
 \ $\bv \in \CC^d$ \ is a left eigenvector of \ $\tbB$ \ corresponding to
 \ $\lambda$, \ then \ $\Re(\bv) \in \RR^d$ \ and \ $\Im(\bv) \in \RR^d$ \ are also
 left eigenvectors of \ $\tbB$ \ or the zero vector.
\proofend
\end{Rem}

\section{A limit theorem for martingales}
\label{ltm}

The next theorem is about the asymptotic behavior of multivariate martingales.

\begin{Thm}{\bf (Crimaldi and Pratelli \cite[Theorem 2.2]{CriPra})}
\label{THM_Cri_Pra}
Let \ $\bigl(\Omega, \cF, (\cF_t)_{t\in\RR_+}, \PP\bigr)$ \ be a filtered
 probability space satisfying the usual conditions.
Let \ $(\bM_t)_{t\in\RR_+}$ \ be a $d$-dimensional martingale with respect to the
 filtration \ $(\cF_t)_{t\in\RR_+}$ \ such that it has c\`{a}dl\`{a}g sample paths
 almost surely.
Suppose that there exists a function
 \ $\bQ : \RR_+ \to \RR^{d \times d}$ \ such that
 \ $\lim_{t\to\infty} \bQ(t) = \bzero$,
 \begin{equation}\label{ltm_cond1}
  \bQ(t) [\bM]_t \bQ(t)^\top \stoch \bfeta
  \qquad \text{as \ $t \to \infty$,}
 \end{equation}
 where \ $\bfeta$ \ is a \ $d \times d$ \ random (necessarily positive semidefinite) matrix
 and \ $([\bM]_t)_{t\in\RR_+}$ \ denotes the quadratic variation (matrix-valued) process of \ $(\bM_t)_{t\in\RR_+}$,
 \ and
 \begin{equation}\label{ltm_cond2}
   \EE\biggl(\sup_{u\in[0,t]} \|\bQ(t) (\bM_u - \bM_{u-})\|\biggr)
   \to 0 \qquad \text{as \ $t \to \infty$.}
 \end{equation}
Then, for each $\RR^{k\times\ell}$-valued random matrix \ $\bA$ \ defined on
 \ $(\Omega, \cF, \PP)$ \ with some \ $k, \ell \in \NN$, \ we have
 \[
   (\bQ(t) \bM_t, \bA) \distr (\bfeta^{1/2} \bZ, \bA) \qquad
   \text{as \ $t \to \infty$,}
 \]
 where \ $\bZ$ \ is a \ $d$-dimensional random vector with
 \ $\bZ \distre \cN_d(\bzero, \bI_d)$ \ independent of \ $(\bfeta, \bA)$.
\end{Thm}

\section*{Acknowledgements}
We would like to thank the referees for their comments that helped us to improve the paper.
This paper has been revised after the sudden death of Gyula Pap, the third author, in October 2019.

\end{document}